\documentstyle[epsfig]{article}

\newtheorem{prop}{Proposition}
\newcommand{\epr}{\end{prop}}
\newcommand{\bpr}{\begin{prop}}

\newcommand{\ee}{\end{equation}}
\newcommand{\be}{\begin{equation}}
\def\sq#1{{\large\overline{\underline{\left|\matrix{\cr~~ #1~~\cr\cr}\right|}}}}
\def\sqs#1{\overline{\underline{\left|\matrix{\cr~~ #1~~\cr\cr}\right|}}}

\begin{document}  


\centerline{ \huge \bf Tabulation of Painlev\'e 6 Transcendents}

\vskip 0.5 cm
\centerline{{\Large Davide Guzzetti}}

\vskip 0.3 cm 
\centerline{1) International School of Advanced Studies SISSA/ISAS,}
\centerline{Via Bonomea 265, 34136 Trieste, Italy. E-mail: guzzetti@sissa.it}

\vskip 0.3 cm 

\centerline{2) Korea Institute of
    Advanced Study KIAS,}
\centerline{ Hoegiro 87(207-43 Cheongnyangni-dong), Dongdaemun-gu, Seoul 130-722, South Korea.}

\begin{abstract}   

\vskip 0.2 cm
 The paper provides the tables of  the critical behaviours at $x=0,1,\infty$ for the Painlev\'e 6 functions. The connection formulae for the basic solutions are also provided, in parametric form. 
\end{abstract}


\section{Introduction}

In the last few decades, the Painlev\'e equations have emerged as one of the central objects in pure mathematics and mathematical physics. They  define non-linear special functions irreducible to classical ones \cite{Umemura},\cite{Umemura1},\cite{Umemura2}, which find applications in a variety of problems, such as number theory, theory of analytic varieties (like Frobenius structures), random matrix theory, orthogonal polynomials, non linear evolutionary PDEs, combinatorial problems, etc. The properties of the classical special functions have
been organised and tabulated in various classical handbooks.  The scientific community is now engaged in the project of a  
comparable organisation and tabulation of the properties of the Painlev\'e functions.   This paper has been written as a contribution to the tabulation and organization of the Painlev\'e 6 functions. 
                                                                                                              
\vskip 0.2 cm 
The paper provides the tables  (Tables 1, 2 and 3) of the  critical behaviours of the Painlev\'e 6 functions  at the three singular points (critical points) $x=0,1,\infty$ of the equation. It also provides, in parametric form,  the connection formulae relating the critical behaviors of a given function at the three singular points  (Sections \ref{dokidoki} and \ref{CONNE}).  

 Though  the critical behaviours  of the Painlev\'e six functions at the critical points, and their parametrisation (in terms of monodromy data to be introduced), have been extensively studied, a systematic organisation of the results  has been missing. This lack is filled by Tables 1, 2 and 3 and by the parametric formulae of Section \ref{dokidoki}     and \ref{CONNE}.  The material of the paper is based on the results  of   \cite{Jimbo}, \cite{Boalch},  \cite{D4}, \cite{D3},  \cite{D2}, \cite{D1}, \cite{guz2010}.  They are essentially all  that have been obtained  by {\it the monodromy preserving deformations method}. 

Hopefully, Tables 1, 2 and 3 are complete, namely they contain all the possible critical behaviours (Section \ref{congettura}). A partial proof is given in  Appendix B, but no final proof is available yet. 

 The critical behaviours in the tables are given in the form of full expansions. 
When  logarithms appear, no proof of the convergence of the expansions  is available at the moment, thus the logarithmic  expansions are to be considered as asymptotic. For all the other cases tabulated, it is  proved that the expansions converge (\cite{Sh}, \cite{D3}, \cite{kaneko}), thus they define true solutions.   It is still to be verified if an extension of the method of \cite{Sh} and \cite{D3}  is possible in order to prove the convergence of the logarithmic expansions as well.

\vskip 0.2 cm 
As for the local analysis is concerned,   the behaviours (\ref{fullEXP}) and (\ref{atopy}) in Table 1  were also obtained in  \cite{Sh} on the universal covering of a punctured neighbourhood of the critical point, by local analysis of integral equations.  An alternative approach to the local analysis is  the method of power geometry  of \cite{Bruno4} and  \cite{Bruno5}. In  as series of papers  summarised in the review \cite{Bruno7},  expansions of actual and possibly formal solutions of PVI have been  constructed by this method. They are classified in a peculiar way, which resents of the method by which they have been constructed.  In section \ref{esplaintable}, we compare these expansions with our tables, and prove  that they coincide.

\vskip 0.3 cm 
\noindent
{\bf Acknowledgements:}  I thank the Korea Institute of Advanced Study (KIAS, Seoul 130-722, 
Republic of Korea)  for providing  computing resources (Abacus System), an excellent research environment without teaching duties, and  many research funds, which created the conditions for this  work.

\subsection{The Sixth Painlev\'e Equation}

The sixth Painlev\'e equation (denoted PVI or $\hbox{PVI}_{\alpha,\beta,\gamma,\delta}$) is: $$
{d^2y \over dx^2}={1\over 2}\left[ 
{1\over y}+{1\over y-1}+{1\over y-x}
\right]
           \left({dy\over dx}\right)^2
-\left[
{1\over x}+{1\over x-1}+{1\over y-x}
\right]{dy \over dx}
$$
$$
+
{y(y-1)(y-x)\over x^2 (x-1)^2}
\left[
\alpha+\beta {x\over y^2} + \gamma {x-1\over (y-1)^2} +\delta
{x(x-1)\over (y-x)^2}
\right],~~~~~\alpha,\beta,\gamma,\delta\in{\bf C}
$$
The critical points are $x=0,1,\infty$. Following  \cite{Its} (page 8), 
 {\it solving} PVI means that: 
 
{\bf ~i)} We determine 
the {\it explicit} critical behaviour (or asymptotic expansion)  of $y(x)$, by  an
 explicit formula in terms of 
two integration constants.

 {\bf ~ii)} We solve the {\it connection
  problem}, namely we find the {\it explicit} 
relations among couples of integration
constants  at different critical points. 

 If we can solve  i) and ii), then  Painlev\'e transcendents can be
efficiently used in applications 
as special functions.  
  Tables 1, 2 and 3, and the formulae of sections \ref{dokidoki} and \ref{CONNE} give answer to i) and ii) respectively.\footnote{A more restrictive definition of ``solving'' should include 
the determination of the distribution of the movable poles, for which we refer to \cite{Dpoli} and \cite{Br}.}

\vskip 0.2 cm 
\noindent 
{\bf Remark 1:} 
{\it critical behaviour} means the behaviour of a solution at one of the three singularities $x=0,1,\infty$ of the equation. These are called  {\it critical points}.  Please, note that this terminology differs from that of singularity theory, where a critical point of a function  is a point where the first derivative vanishes.


\subsection{Monodromy Data}
According to  \cite{JMU}, 
PVI is the isomonodromy deformation equation of  the $2\times 2$ Fuchsian system: 
\be
   {d\Psi\over d\lambda}=A(x,\lambda)~\Psi,~~~~~
A(x,\lambda):=\left[ {A_0(x)\over \lambda}+{A_x(x) \over \lambda-x}+{A_1(x)
\over
\lambda-1}\right],~~~\lambda\in{\bf C}.
\label{SYSTEM}
\ee
The matrix $A(x,\lambda)$ can be written as
$$
A_{12}(x,\lambda)={g(x)(\lambda-y(x))\over \lambda(\lambda-1)(\lambda-x)}
$$
where $g(x)$ is a certain algebraic function of $x$. 
Therefore, 
\be
\label{ipsilon}
y(x)={x(A_0)_{12}\over x[(A_0)_{12}+(A_1)_{12}]-(A_1)_{12}}
\ee
The  $2\times 2$ matrices  $A_i(x)$  depend 
 on $x$ in such a way that there exists a fundamental matrix solution $\Psi(\lambda,x)$ with monodromy  
independent of  small deformations of $x$. They also depend algebraically on $\alpha,\beta,\gamma,\delta$  according to 
the following relations:  
\begin{equation}
 A_0+A_1+A_x = -{\theta_{\infty}\over 2}
 \sigma_3,~~\theta_\infty\neq 0.~~~~~
\hbox{ Eigenvalues}~( A_i) =\pm {1\over 2} \theta_i, ~~~i=0,1,x;
\label{caffe0}
\end{equation}
 \begin{equation}
    \alpha= {1\over 2} (\theta_{\infty} -1)^2,
~~~-\beta={1\over 2} \theta_0^2, 
~~~ \gamma={1\over 2} \theta_1^2,
~~~ \left({1\over 2} -\delta \right)={1\over 2} \theta_x^2 ,~~~~~\theta_\infty\neq 0
\label{caffe1}
\end{equation}

\begin{figure}
\epsfxsize=8cm
\centerline{\epsffile{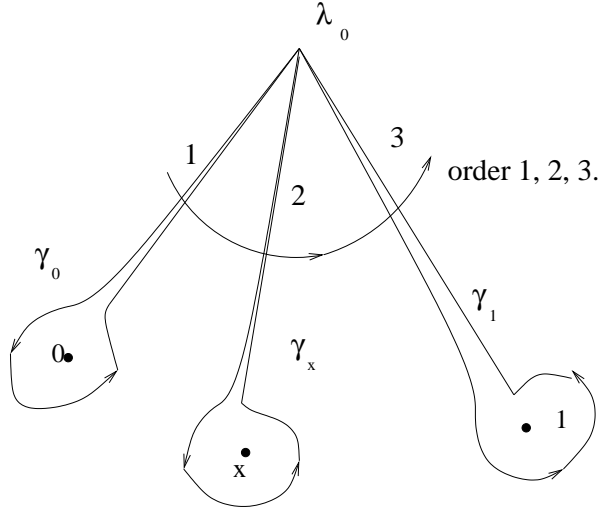}}
\caption{The ordered basis of loops}
\label{figure1}
\end{figure}

In  the ``$\lambda$-plane'' 
${\bf C}\backslash\{0,x,1\}$ we fix a base point $\lambda_0$ 
and  three loops, which are numbered in order 1, 2, 3 according to a 
counter-clockwise order referred to $\lambda_0$.  We choose $0,x,1$ to be the order $1,2,3$. 
We denote the loops by $\gamma_0$, $\gamma_x$, $\gamma_1$. See figure \ref{figure1}.   
  The loop at infinity will be $\gamma_\infty=\gamma_0\gamma_x\gamma_1$. 
  When $\lambda$ goes around a small loop around $\lambda=i$, $i=0,x,1$,  the  fundamental solution transforms like $\Psi \mapsto \Psi M_i$, where $M_0$, $M_x$, $M_1$ are the monodromy matrices w.r.t. the base of loops.

\vskip 0.2 cm
Let $(\theta_0, \theta_1, \theta_x,\theta_\infty)\in{\bf C}^4$ be fixed by PVI, up to the equivalence $\theta_k\mapsto -\theta_k$, $k=0,x,1$, and $\theta_\infty\mapsto 2-\theta_\infty$. Denote $\sim$ the equivalence and let  
$$
\Theta:=\{(\theta_0, \theta_1, \theta_x, \theta_\infty)\in{\bf C}^4~|~\theta_\infty\neq 0\}/ \sim
$$
 be the quotient ($\theta_\infty$ may only be equal to 2). 
  Let $M_{\infty}:=M_1M_xM_0$ be the monodromy at $\lambda=\infty$,  and consider  the set of triples of monodromy matrices, defined up to conjugation $M_i\mapsto CM_i C^{-1}$ ($i=0,x,1$) by an invertible matrix $C$, namely  
$$ 
  M:=\{(M_0,M_x,M_1)~|~\hbox{\rm det}M_0=\hbox{\rm det}M_1=\hbox{\rm det}M_x=1, ~\hbox{Tr}M_\mu= 
2\cos\pi\theta_\mu,~\mu=0,1,x,\infty\}/\hbox{conjugation}
$$

\vskip 0.2 cm 
\noindent
{\bf Definition}: The {\it  monodromy data} of the class of Fuchsian systems (\ref{SYSTEM}), with the basis of loops ordered as  figure \ref{figure1}, is the set ${\cal M}:=\Theta\cup M$. 
\vskip 0.2 cm 
\noindent
 For an equivalent definition, see section 2 of \cite{guz2010}.

\subsection{Critical Behaviours in terms of Monodromy Data -- Parametric Connection Formulae}

We fix branch cuts in the $x$ plane, for example  $-\pi<\arg x<\pi$ and $-\pi<\arg (1-x) <\pi$.  
To every branch $y(x)$, a  system (\ref{SYSTEM}) is associated, and so is a point in ${\cal M}$. Conversely, to a point in ${\cal M}$ a system or a family of systems (\ref{SYSTEM}) is associated through a Riemann-Hilbert problem \cite{AB}, and so is  either one branch $y(x)$ or a family of branches $y(x)$. Let  
\be
\label{monmap}
 f:\{y(x)\hbox{ branch}\} \rightarrow {\cal M}
\ee
be the map from the set $\{y(x)\hbox{ branch}\}$  of all the branches  of all equations $\hbox{PVI}_{\alpha,\beta,\gamma,\delta}$, $(\alpha,\beta,\gamma,\delta)\in{\bf C}^4$,  onto ${\cal M}$, associating to a branch the corresponding monodromy data. 

\bpr  Let the order of loops be fixed. The map (\ref{monmap})  is injective (one-to-one) when restricted to $f^{-1}(\Theta\cup \{(M_0,M_x,M_1)\in M ~|~M_\mu\neq I,~~\forall\mu=0,x,1,\infty\})$ \cite{D1}.

\epr

 \vskip 0.2 cm 
\noindent 
Let 
$$
\sqs{p_{ij}:=\hbox{Tr}(M_iM_j), ~~~j=0,x,1;~~~~~~ p_\mu:=\hbox{Tr}M_\mu=2\cos\pi\theta_\mu, ~~~\mu=0,x,1,\infty. }
$$
Observe that $p_{ij}=p_{ji}$. These are seven invariant functions (w.r.t. conjugation and $\sim$) defined on ${\cal M}$.  They satisfy  the affine cubic Jimbo-Fricke relation (discovered by Jimbo \cite{Jimbo}, named Fricke cubic in \cite{Boalch},  studied in \cite{Iwa}):
$$
p_{0x}^2+p_{01}^2+p_{x1}^2+p_{0x}p_{01}p_{x1}-(p_0p_x+p_1p_\infty)p_{0x}-(p_0p_1+
p_xp_\infty)p_{01}-(p_xp_1+p_0p_\infty)p_{x1}+
$$
\be
\label{FR}
+p_0^2+p_1^2+p_x^2+p_\infty^2+p_0p_xp_1p_\infty -4 =0
\ee
 in agreement with the dimension of ${\cal M}$ 
\footnote{
The above relation follows by taking the trace of the relation $M_1M_xM_0=M_\infty$. 
}
\footnote{
The cubic curve is factorisable when: $p_{0x}=2\cos\pi(\theta_\infty\pm\theta_1)$, or $p_{x1}=2\cos\pi(\theta_\infty\pm\theta_0)$, or $p_{01}=2\cos\pi(\theta_\infty\pm\theta_x)$. These values correspond to reducibility of the subgroups generated by $M_xM_0,M_1$, or $M_1M_x,M_0$, or $M_0M_1,M_x$ respectively. 
}.  Except for special cases, they are local coordinates on ${\cal M}$, according to the following proposition:

\bpr
    $p_0,p_x,p_1,p_\infty$; $p_{0x}, p_{x1}, p_{1x}$     are   coordinates on the subset of ${\cal M}$ where the group generated by  $M_0$, $M_x$, $M_1$ is irreducible  \cite{Iwa}.

\epr

Suppose that  Propositions 1 and 2 hold for a given $\hbox{PVI}_{\alpha,\beta,\gamma,\delta}$. The $p_\mu$'s  $\mu=0,x,1,\infty$ are fixed by $\alpha,\beta,\gamma,\delta$ and only two parameters $p_{ij}$'s are independent. They have the meaning of  integration constants for the associate branch $y(x)$.  Consider now a branch $y(x)$ close to a critical point. It has a critical behaviour 
\be
\label{malamente1}
\sqs{y(x)=y(x,c_1,c_2)}
\ee
 depending on two integration constants $c_1,c_2$. The asymptotic techniques of the method of monodromy preserving deformations provide explicit {\it parametric formulae}  
\be
\label{POA}
\sq{
\left\{
\matrix{ c_1=c_1(\theta_0,\theta_x,\theta_1,\theta_\infty,p_{0x},p_{x1},p_{01})
\cr 
c_2=c_2(\theta_0,\theta_x,\theta_1,\theta_\infty,p_{0x},p_{x1},p_{01})
}
\right.~~~\hbox{ or } ~~~
\left\{
\matrix{
 c_1=c_1(\alpha,\beta,\gamma,\delta,p_{0x},p_{x1},p_{01})
\cr 
c_2=c_2(\alpha,\beta,\gamma,\delta,p_{0x},p_{x1},p_{01})
}
\right.
}
\ee
Conversely, 
the monodromy data are  explicitly computable in terms of the integration constants . 
\be
\label{POA1}
\sq{
\left\{
\matrix{
p_{0x}=p_{0x}(c_1,c_2,\theta_0,\theta_x,\theta_1,\theta_\infty) 
\cr 
p_{x1}=p_{x1}(c_1,c_2,\theta_0,\theta_x,\theta_1,\theta_\infty) 
\cr 
p_{01}=p_{01}(c_1,c_2,\theta_0,\theta_x,\theta_1,\theta_\infty) 
}
\right.
~~~\hbox{ or } ~~
\left\{
\matrix{
p_{0x}=p_{0x}(c_1,c_2,\alpha,\beta,\gamma,\delta) 
\cr 
p_{x1}=p_{x1}(c_1,c_2,\alpha,\beta,\gamma,\delta) 
\cr 
p_{01}=p_{01}(c_1,c_2,\alpha,\beta,\gamma,\delta) 
}
\right.
}
\ee
Explicit means that the formulae are classical functions of their arguments.

\vskip 0.2 cm 
   {\it Three pairs of different  parametric formulae of type (\ref{POA}) can be written at $x=0,1$ and $\infty$ respectively,   
in terms of the same monodromy data, namely for the same branch $y(x)$. Conversely, the monodromy data associated to  a given 
$y(x)$ can be written as in (\ref{POA1})  in three ways, namely in terms of the three couples of integration constants  at $x=0$, 
$1$ and $\infty$ respectively. This  solves the connection problem.} 

For this reason, the formulae  (\ref{POA}) and   (\ref{POA1}) will be refered to as  {\bf parametric connection formulae}. The {\it connection formulae in closed form}, namely the formulae expressing the integration constants appearing in the critical behaviour at one critical point in terms of the integration constants at another critical point (see (\ref{closd}) in Section \ref{CONECLOSE}), are obtained combining those in parametric form at two critical points, with the natural procedure explained in Section \ref{CONECLOSE}. 

\vskip 0.2 cm

Jimbo \cite{Jimbo}  provided for the first time the critical behaviour and its parametrisation in terms of monodromy data for a wide class of two complex parameter solutions
$$
y(x)=
\left\{
\matrix{
a_0 ~x^{1-\sigma_{0x}}(1+O(x^\epsilon)) & x\to 0
\cr 
1+a_1 ~(1-x)^{1-\sigma_{x1}}(1+O((1-x)^{\epsilon})) & x\to 1
\cr 
a_\infty ~x^{\sigma_{01}}(1+O(x^{-\epsilon})) &x\to \infty
}
\right.
$$
where $\epsilon>0$ is small and $0\leq\Re\sigma_{0x},\Re\sigma_{01},\Re\sigma_{x1}<1$. The leading behaviour at each point is governed by two integration constants $(a_0,\sigma_{0x})$,  $(a_1,\sigma_{x1})$, $(a_\infty,\sigma_{01})$. Their explicit expression in terms of monodromy data is computed in \cite{Jimbo} and reproduced here in formula (\ref{aF}) of Section \ref{dokidoki}. The behaviours at $x=0,1$ and $\infty$ are characterised by the {\it monodromy exponents} $\sigma_{0x}$, $\sigma_{x1}$ and $\sigma_{01}$ respectively,  where 
$$
2\cos\pi\sigma_{0x}=p_{0x},~~~~~2\cos\pi\sigma_{x1}=p_{x1},~~~~~2\cos\pi\sigma_{01}=p_{01}
$$
 Jimbo's work  \cite{Jimbo}  is the foundation of all subsequent developments. 
In this paper we use the results of \cite{Jimbo}, \cite{Boalch},  \cite{D4}, \cite{D3},  \cite{D2}, \cite{D1}, \cite{guz2010} to provide behaviours (\ref{malamente1}) in Tables 1, 2 and 3, and the formulae of type (\ref{POA}) and (\ref{POA1}) in Section \ref{dokidoki} .

\subsection{Completeness of  the Tables 1,2, and 3 of Critical behaviors} 
\label{congettura}

 Are Tables 1, 2 and 3  complete, namely  do they provide all the critical behaviours?  We formulate the following conjecture and give strong motivations for it in Appendix B. 

\vskip 0.2 cm 
\noindent
{\bf Conjecture:}
Tables 1, 2 and 3 of Sections \ref{TABLE}, \ref{TABLE2} and \ref{TABLE3}  respectively are complete, namely they include all the  critical behaviours at $x=0,1$ and $\infty$ respectively.

\subsection{Organisation of the Paper}

\noindent
-- In Sections \ref{TABLE}, \ref{TABLE2} and \ref{TABLE3} the critical behaviours (\ref{malamente1}) are provided at $x=0$, 1 and $\infty$ respectively, in the form of  explicit expansions tabulated in Tables 1, 2 and 3 respecively. 

\vskip 0.2 cm 
\noindent
-- Section \ref{dokidoki} provides  the parametric connection formulae of type (\ref{POA}) and (\ref{POA1})  at $x=0$.

\vskip 0.2 cm 
\noindent
-- Section \ref{CONNE} provides:
\vskip 0.2 cm 
~~~ a) 
The procedure to write the parametric connection formulae of type (\ref{POA}) and (\ref{POA1})  at $x=1$ and $\infty$, starting from the same formulae at $x=0$.

\vskip 0.2 cm 
~~~ b) The procedure to write the connection formulae in closed form.

\vskip 0.2 cm \noindent
-- Section \ref{esplaintable} explains how Table 1 is constructed with  references to the results in the literature. 

\vskip 0.2 cm \noindent
-- The conjecture of subsection \ref{congettura} is motivated in Appendix B.


\section{Table 1:  Critical Behaviours when $x\to 0$} 
\label{TABLE}

$\diamond$ The table provides the critical behaviors (\ref{malamente1})  at $x=0$. The branch cut is    $-\pi <
\hbox{arg}  ~x\leq \pi$.  The detailed description of the critical behaviours tabulated  is in section \ref{esplaintable}, to which the numeration of the formulae refers.  The term "basic solutions" is explained in Section \ref{dokidoki}.

\vskip 0.3 cm 
 
\noindent
$\diamond$ The branches in the table may be classified according to their behaviour as follows: 

\vskip 0.2 cm 
- Complex  power behaviours (\ref{fullEXP}),  (\ref{atopy}), (\ref{UUU})  and   (\ref{div}). They are expanded in powers of $x^{n+m\lambda}$, for some $n,m\in{\bf Z}$ and $\lambda\in{\bf C}$. In this case,   $|y(x)|$ may vanish, converge to a constant or diverge when $x\to 0$.

\vskip 0.2 cm 
- Inverse oscillatory behaviours (\ref{lantern1}) and (\ref{TAU}).  $y(x)$ oscillates without vanishing when $x\to 0$, and may have poles in a sector centred at $x=0$.

\vskip 0.2 cm 
- Integer power behaviours, namely Taylor series (\ref{TTLO1}), (\ref{TTLO2}), (\ref{TTLO3}), (\ref{TTLO4}), (\ref{davidekan}) and  (\ref{T1coe}).

\vskip 0.2 cm 
- Logarithmic behaviours,  of type (\ref{log1}),   (\ref{log1zero}), (\ref{logsquare}) and (\ref{LOG12}), namely power expansions  with coefficients which are polynomials of $\ln x$. $y(x)\to$ constant as $x\to 0$, where the constant may be zero or not. 

\vskip 0.2 cm 
- Inverse  logarithmic behaviours of type (\ref{LOG45}) and (\ref{LOG3}), which can be expanded as asymptotic series of $( \ln x)^{-1}$. $y(x)=O(1/\ln x)$, or $O(1/\ln^2 x)$,  as $x\to 0$.

\vskip 0.3 cm 
\noindent 
$\diamond$ 
Except for the case when logarithms appear in the expansions tabulated, the expansions are convergent on suitable domains of the universal covering of a small neighbourhood of $x=0$ with $x=0$  removed (punctured neighbourhood). The proof is based on the results of \cite{Sh} and \cite{D3}.  To be more precise,  for any $\vartheta>0$ there exists  $r(\vartheta)$  small enough (decreasing function of $\vartheta$) such that the expansion of the complex power behaviors converges for  $|\arg(x)|<\vartheta$ and $0<|x|<r(\vartheta)$. If  $\Re\sigma=0$, there is also an  additional constraint $\arg (x)>\varphi_0$, for a suitable $\varphi_0$ fixed by the integration constants.  Also the expansion of the denominator of the inverse oscillatory behaviours is convergent under  the same conditions $|\arg(x)|<\vartheta$, for any choice of a  $\vartheta>0$, and $0<|x|<r(\vartheta)$, plus the additional contraint $\arg (x)<\varphi_0$, for a suitable $\varphi_0$.

 When logarithms appear, no proof of convergence is known to the author, though one may expect that an extension of the method of \cite{Sh} and \cite{D3}  be possible for the power logarithmic expansions. 

Convergence of the Taylor expansions is studied in \cite{kaneko}.

\vskip 0.3 cm 
\noindent 
$\diamond$ 
For the critical behaviors where a ``star'' $\star$ appears close to the numeration of the formula, the parametric connection formulae are given in Section \ref{dokidoki}. In other cases, the formulae can be computed from the cases $\star$ via a birational transformation.  

\vskip 0.3 cm
\noindent 
$\diamond$ 
The table shows also the reducibility of the subgroups $<M_0,M_x>$ and  $<M_xM_0,M_1>$ (here $<A,B,..>$ means group generated by $A,B,...$.)
 
\vskip 0.3 cm
\noindent 
$\diamond$ Notations: 

\vskip 0.2 cm  
-- In the table, $\sigma,\phi$ and $a$ denote complex free parameters (integration constants), and $\nu$ a real free parameter. The coefficients $c_{nm}$, $d_{nm}$ and  $b_n$ are rational functions of $\sqrt{\alpha}$, $\sqrt{\beta}$, $\sqrt{\gamma}$ and $\sqrt{1-2\delta}$. The coefficients $b_n(a)$  are rational functions of $\sqrt{\alpha}$, $\sqrt{\beta}$, $\sqrt{\gamma}$, $\sqrt{1-2\delta}$ and $a$. $P_n(\ln x, a)$ are polynomials in $\ln x$ with coefficients which are rational functions of $\sqrt{\alpha}$, $\sqrt{\beta}$, $\sqrt{\gamma}$, $\sqrt{1-2\delta}$ and $a$. These coefficients can be recursively computed (essentially, by substitution into PVI. See \cite{guz2010}). In the case of (\ref{fullEXP}), the  $c_{nm}$'s  depend rationally also on $\sigma$. In the case of (\ref{lantern1}), the $c_{nm}$'s depend rationally also on $\nu$.

\vskip 0.2 cm 
-- Numbers:    

\be
\label{SeT}
\Sigma_{\beta\delta}^{k}:= \pm 
\Bigl(\sqrt{-2\beta}+(-)^k\sqrt{1-2\delta}\Bigr)  ,~ ~~~~\hbox{if $|\sqrt{-2\beta}+(-)^k\sqrt{1-2\delta}|<1$}.
\ee
otherwise 
\be
\label{iku}
\Sigma_{\beta\delta}^{k}:=
 \Bigl(\sqrt{-2\beta}+(-)^k\sqrt{1-2\delta}\Bigr)\hbox{sgn}\Bigl(\Re(\sqrt{-2\beta}+(-)^k\sqrt{1-2\delta})\Bigr)
,~~~k=1,2,
\ee
In the monodromy deformation parameters, $\Sigma_{\beta\delta}^k
\in\{ (\theta_0+\theta_x),(\theta_0-\theta_x),-(\theta_0+\theta_x),-(\theta_0-\theta_x)\} 
$, or  $\Sigma_{\beta\delta}^k=
 (\theta_0+(-)^k\theta_x)\hbox{sgn}(\Re(\theta_0+(-)^k\theta_x))$. 
\be
\label{SeT1}
\Omega_{\alpha\gamma}^{k} 
:=\pm
\Bigl( \sqrt{2\alpha}+(-)^k\sqrt{2\gamma}\Bigr) 
,~~~~~\hbox{ if $|\Re \{\sqrt{2\alpha}+(-)^k\sqrt{2\gamma} \}|<1$}. 
\ee
otherwise
\be
\label{STaR}
\Omega_{\alpha\gamma}^{k}:=
\Bigl(\sqrt{2\alpha}+(-)^k\sqrt{2\gamma}\Bigl)\hbox{sgn}\Bigl(\Re(\sqrt{2\alpha}+(-)^k\sqrt{2\gamma})\Bigr)
,~~~k=1,2.
\ee
In the monodromy deformation parameters, $\Omega_{\alpha\gamma}
\in\{(\theta_\infty-1+\theta_1), (\theta_\infty-1-\theta_1),-(\theta_\infty-1+\theta_1), -(\theta_\infty-1-\theta_1)\}
$, or 
$\Omega_{\alpha\gamma}
=(\theta_\infty-1+(-)^k\theta_1)\hbox{sgn}(\Re(\theta_\infty-1+(-)^k\theta_1))
$. 

\vskip 0.2 cm 
-- A set: 
\be
\label{infinitapalla}
 {\cal N}_N:=\left\{
\matrix{
\{0,2,4,...,|N|-1\}\cup\{-2,-4,...,-(|N|-1)\} \hbox{ if $N$ is odd}
\cr
 \{1,3,...,|N|-1\}\cup\{-1,-3,...,-(|N|-1)\} \hbox{ if $N$ is even}
}
\right.
\ee

\vskip 0.3 cm 
\noindent
$\diamond$  
{\large\bf How to identify a Critical Behavior from given Monodromy Data:}

\vskip 0.2 cm 
Preliminary, observe that when  $p_{0x}=-2\cos\pi \Omega_{\alpha\gamma}^k$ the Jimbo-Fricke cubic (\ref{FR}) is factorised. Namely, if  $p_{0x}=2\cos\pi(\theta_\infty-\theta_1)$, then (\ref{FR}) is
$$
\left[ p_{01}+p_{x1}e^{i\pi(\theta_\infty-\theta_1)}-2\left(e^{-i\pi\theta_1}\cos\pi\theta_0+e^{i\pi\theta_\infty}\cos\pi\theta_x\right)\right]
\times
$$
\be
\label{era}
\times
\left[
p_{01}+p_{x1}e^{-i\pi(\theta_\infty-\theta_1)}-2\left(e^{i\pi\theta_1}\cos\pi\theta_0+e^{-i\pi\theta_\infty}\cos\pi\theta_x\right)
\right]=0
\ee
If $p_{0x}= 2\cos \pi (\theta_\infty+\theta_1)$, then (\ref{FR}) is: 
$$
\left[p_{01}+p_{x1}e^{i\pi(\theta_\infty+\theta_1)}-2\left(e^{i\pi\theta_1}\cos\pi\theta_0+e^{i\pi\theta_\infty}\cos\pi\theta_x\right)
\right]
\times
$$
\be
\label{era1}
\times
\left[
p_{01}+p_{x1}e^{-i\pi(\theta_\infty+\theta_1)}-2\left(e^{-i\pi\theta_1}\cos\pi\theta_0+e^{-i\pi\theta_\infty}\cos\pi\theta_x\right)
\right]=0
\ee

\vskip 0.3 cm 
 
\noindent
The behaviours at $x=0$ are determined by  $p_{0x}$, which gives the monodromy exponent $\sigma=\sigma_{0x}$ in generic cases (see (\ref{fullEXP}) and (\ref{atopy})), the monodromy data $p_{01}, p_{x1}$, and the values of the numbers  $\Sigma_{\beta\delta}^k$ and $\Omega_{\alpha\gamma}^k$. The following explains how the critical behaviour is singled out:  
$$
\left[\matrix{
p_{0x}<-2 & \Rightarrow & \left\{\matrix{(\ref{lantern1}) \hbox{ if } p_{0x}\neq -2\cos\pi\Omega_{\alpha,\gamma}^k
         \cr (\ref{TAU})  \hbox{ if } p_{0x}= -2\cos\pi\Omega_{\alpha,\gamma}^k
}\right.
\cr
\cr
p_{0x}\neq \left\{\matrix {\pm2,-2\cos\pi\Sigma_{\beta\delta}^k,
\cr
 2\cos\pi\Omega_{\alpha,\gamma}^k
\cr 
\not <-2
}\right.
& 
\Rightarrow 
&
(\ref{fullEXP})
\cr
\cr
p_{0x}=2\cos\pi\Sigma_{\beta\delta}^k\neq \pm 2 
&
\Rightarrow
&
(\ref{atopy}) \hbox{ [or (\ref{davidekan})]}
\cr 
\cr 
p_{0x}=-2 \cos\pi\Omega_{\alpha\gamma}^k\neq \pm 2
&
\Rightarrow
&
\left\{
\matrix{
 \hbox{(\ref{fullEXP}),  if  one factor in (\ref{era}) or (\ref{era1}) $=0$ (*Note)}
\cr 
\hbox{(\ref{UUU}) [or (\ref{T1coe})], if $\alpha\neq 0$ and the other factor in (\ref{era}) or (\ref{era1}) $=0$}
\cr
\hbox{(\ref{div}), if $\alpha = 0$ and  the other factor in (\ref{era}) or (\ref{era1}) $=0$}
}\right.
\cr
\cr
p_{0x}=\pm 2 
&
\Rightarrow
&
\hbox{Taylor or Logarithmic, depending on  $\Sigma_{\beta\delta}^k,\Omega_{\alpha\gamma}^k,\alpha,\beta,\gamma,\delta$
}
}
\right]
$$
*Note: When $p_{0x}=-2 \cos\pi\Omega_{\alpha\gamma}^k\neq \pm 2$, which  of the factors in (\ref{era}) or (\ref{era1}) provides (\ref{fullEXP}), and which provides (\ref{UUU}) [or (\ref{T1coe}), (\ref{div})],  depends on the actual value of $\Omega_{\alpha\gamma}^k$. See Remark 3 in Section \ref{dokidoki} 


$$
\begin{array}{||c c|c|c||} 
\hline\hline  &  &  &  \cr 
 & \hbox{{\bf  Complex power behaviours}} & \matrix{\hbox{Free}\cr
                                                   \hbox{ Param.}
                                                   }
 & \hbox{Other Conditions}    \cr
\hline          
 & & &  \cr
\hbox{(\ref{fullEXP})} \star     
& 
y(x)= \sum_{n=1}^\infty x^n\sum_{m=-n}^n c_{nm}(ax^{\sigma})^m & \matrix{\sigma \cr
\cr  a\neq 0} & 
 \matrix{
0\leq \Re \sigma<1,~~~\sigma\neq \Sigma_{\beta\delta}^k,
\cr
\cr
2\cos\pi\sigma=p_{0x}.
}
\cr 
&
\matrix{
=\left\{
\matrix{ {c_{1,-1}\over a} x^{1-\sigma} (1+O(x^\sigma,x^{1-\sigma})),&  0<\Re \sigma <1
\cr
\cr
x\left[{c_{1,-1}\over a} x^{-\sigma}+c_{10}+ax^\sigma\right]+O(x^2), & \Re \sigma=0
}\right.
\cr 
\cr 
\cr 
c_{11}=1,~~c_{10}={\sigma^2-2\beta+2\delta-1\over 2\sigma^2}
\cr
c_{1,-1}= 
{\Bigl[(\sqrt{-2\beta}-\sqrt{1-2\delta})^2-\sigma^2\Bigr]  \Bigl[(\sqrt{-2\beta}+\sqrt{1-2\delta})^2-\sigma^2\Bigr]\over 16 \sigma^4} 
\cr 
\cr 
\hbox{Basic solutions}
}
&
                                  & 
\matrix{ p_{0x}\neq 2\cos\pi\Sigma_{\beta\delta}^k,\pm 2,
\cr 
p_{0x}\not\in(-\infty,-2].
}
\cr 
& &   &   \cr
\hline
&&&\cr
%
%
\hbox{(\ref{atopy})} \star & 
\matrix{y(x)= \sum_{n=1}^\infty x^n\sum_{m=0}^n c_{nm} (ax^\sigma)^m 
\cr
\cr
\cr 
c_{11}=1,~~
c_{10}={\sqrt{-2\beta}\over\sqrt{-2\beta}+(-)^k\sqrt{1-2\delta}}
\cr 
\cr
\hbox{If $a=0$, $y(x)$ reduces to (\ref{davidekan})}
}
& 
a
&
\matrix{\Re\sigma>-1,~~
 \sigma=\Sigma_{\beta\delta}^{k}\not\in{\bf Z},
\cr
\Sigma_{\beta\delta}^{k} \hbox{ is (\ref{iku}) or (\ref{SeT})}.
\cr 
\cr 
p_{0x}=2\cos\pi\Sigma_{\beta\delta}^k\neq \pm 2.
\cr 
\cr
\cr
<M_0,M_x>\hbox{ reducible.}
}
\cr 
&\hbox{Basic solutions}&&
\cr
&&& \cr 
\hline &&& \cr 
%
%
\hbox{(\ref{UUU})} \star &
\matrix{
y(x)=d_{00} + \sum_{n=1}^\infty x^n \sum_{m=0}^n d_{nm}(\tilde{a}x^{\rho})^m
\cr
\cr
\cr 
d_{00}={\sqrt{\alpha} +(-)^k\sqrt{\gamma} \over \sqrt{\alpha}},~~
\tilde{a}=-a~d_{00}^2,~~~d_{11}=1
\cr 
\cr 
\hbox{If $\tilde{a}=0$, $y(x)$ reduces to (\ref{T1coe})}
}
&
a
&
\matrix{
\alpha\neq 0.
\cr
\cr 
\rho=\Omega_{\alpha\gamma}^{k}-1, ~~\Re\rho>-1,
\cr
\Omega_{\alpha\gamma}^{k}\not\in{\bf Z},~~~
\Omega_{\alpha\gamma}^{k} \hbox{ is (\ref{STaR})}.
\cr 
\cr  
p_{0x}=-2\cos\pi\Omega_{\alpha\gamma}^k\neq \pm 2.
\cr 
\cr
<M_xM_0,M_1>\hbox{ reducible.}
}
\cr
 &&& \cr 
\hline &&& \cr 
%
%
\hbox{(\ref{div})} \star &
y(x)={ 1\over a }x^{-\omega}
\left(
1+\sum_{n=1}^\infty x^n \sum_{m=0}^n d_{nm}(a x^{\omega})^m
\right)
&
a
&
\matrix{\alpha =0,
\cr 
\gamma\not\in(-\infty,0],~~\sqrt{2\gamma}\not \in {\bf Z}.
\cr
\cr 
\omega=\sqrt{2\gamma}~\hbox{sgn}(\Re\sqrt{2\gamma}),
\cr
\Re\omega > 0.
\cr 
\cr
p_{0x}=-2\cos\pi\sqrt{2\gamma}\neq \pm 2.
\cr 
\cr
<M_xM_0,M_1>\hbox{ reducible.}
}
\cr 
&&&
\cr
\hline
\end{array}
$$


\newpage 
$$
\begin{array}{||c c|c|c||} 
\hline\hline
&&&\cr
&\hbox{{\bf Inverse Oscillatory Behaviours}} & \left.\matrix{\hbox{Free}\cr\hbox{Param}}\right. & \hbox{ Other }
\cr
\hline
&&&
\cr
\hbox{(\ref{lantern1})} \star 
& 
\matrix{
y(x)= \Bigl[\sum_{n=0}^\infty x^n\sum_{m=-n-1}^{n+1} c_{nm}
\bigl(e^{i\phi}x^{2i\nu}\bigr)^m
\Bigr]^{-1} 
\cr
\cr 
= \left[A\sin(2\nu \ln x +\phi)+B+ O(x)\right]^{-1}
\cr 
\cr
A=-\sqrt{{\alpha\over 2\nu^2}+B^2},~~
~~~B={2\nu^2+\gamma-\alpha\over 4\nu^2}
}
                          &
 \matrix{\nu 
\cr 
\cr 
\phi}
                             &
\matrix{
\nu\in{\bf R}\backslash\{0\},~~ 2i\nu\neq\Omega_{\alpha\gamma}^k,
\cr 
2\cosh2\pi\nu=-p_{0x}.
\cr 
\cr
\cr 
p_{0x}<-2,
\cr
p_{0x}\neq -2\cos\pi\Omega_{\alpha\gamma}^k.
}
\cr
&&&
\cr 
\hline
&&&\cr
%
%
\hbox{(\ref{TAU})}   \star & 
y(x)=\left[{\sqrt{\alpha} \over \sqrt{\alpha} +(-)^k \sqrt{\gamma}}+ax^{-2i\nu}+\sum_{n=1}^\infty x^n
\sum_{m=0}^{n+1}c_{n+1,m}\bigl(ax^{-2i\nu}\bigr)^m
\right]^{-1}
&
a
&
\left.
\matrix{2i\nu= \Omega_{\alpha\gamma}^{k}\in i {\bf R}\backslash\{0\},
\cr
\Omega_{\alpha\gamma}^{k}  \hbox{ is (\ref{SeT1})}.
\cr 
\cr 
p_{0x}=-2\cos\pi\Omega_{\alpha\gamma}^k<-2.
 \cr 
\cr
<M_xM_0,M_1>\hbox{ reducible.}
}
\right.
\cr&&&\cr
\hline
\end{array}
$$



$$
\begin{array}{||c c|c|c||}
\hline\hline
&&&\cr
 &\hbox{{\bf Taylor expansions}} & \left.\matrix{\hbox{Free}\cr\hbox{Par}}\right. & \hbox{ Other Conditions}
\cr
\hline
&&&
\cr
\hbox{(\ref{davidekan})}&
\matrix{
 y(x)={\sqrt{-2\beta}\over \sqrt{-2\beta}+(-)^k\sqrt{1-2\delta}}x+\sum_{n=2}^\infty b_nx^n &  
\cr 
\cr 
y(x)=0 \hbox{ if } \beta=0
\cr
\cr 
\hbox{This is  (\ref{atopy}) when $a=0$}
}
& &
                  \matrix{ \sqrt{-2\beta}+(-)^k\sqrt{1-2\delta}\not\in{\bf Z}.
\cr 
\cr
 p_{0x}= 2\cos\pi\Sigma_{\beta\delta}^k\neq \pm 2.
\cr 
\cr 
<M_0,M_x>\hbox{ reducible.}
}
\cr
&&& \cr
\hline
%
%
&&&
\cr
\hbox{(\ref{TTLO1})}&
y(x)= \sum_{n=1}^{|N|}b_nx^n+ax^{|N|+1}+\sum_{n= |N|+2}^\infty b_n(a)x^n
&
&
                                        \matrix{\sqrt{-2\beta}+\sqrt{1-2\delta} =N\neq 0
\cr\hbox{ or} \cr
                                                       \sqrt{-2\beta}-\sqrt{1-2\delta}=N\neq 0
}
\cr
                 &               &              &           \cr
&        b_1={\sqrt{-2\beta}\over N}\equiv {\sqrt{-2\beta}\over \sqrt{-2\beta}+\sqrt{1-2\delta}}\hbox{ or } {\sqrt{-2\beta}\over \sqrt{-2\beta}-\sqrt{1-2\delta}}                 
&
 a 
&
                                                        \hbox{\bf and}
\cr
%
%
%
             &       
&
   &
      \matrix{
                           \sqrt{-2\beta}\in \left\{\matrix{
 \{0,-1,-2,...,N\}~~N<0
\cr
    \{0,~1,~2,...,N\}~~~N>0                                                    
                                                      }\right.           
 \cr
\hbox{ or }
\cr
                                  \{ \sqrt{2\alpha}+\sqrt{2\gamma},\sqrt{2\alpha}-\sqrt{2\gamma}\} 
\cap {\cal N}_N \neq \emptyset.
\cr 
\cr 
 p_{0x}=2\cos \pi N=\pm2.
\cr 
\cr
<M_0,M_x>\hbox{ reducible.}
}
\cr
&&&\cr
\hline
&&&\cr
\hbox{(\ref{TTLO2})}\star                
&       
                                                       y(x)= ax+a(a-1)\left(\gamma-\alpha-{1\over 2}\right)x^2+\sum_{n= 3}^\infty b_n(a)x^n
&
        a
& 
\matrix{2\beta=2\delta-1=0.
\cr 
 p_{0x}=2.
\cr
<M_0,M_x>\hbox{ reducible.}
}
\cr
&&&\cr
\hline
&&&\cr
%
%
\hbox{(\ref{T1coe})}
&
\matrix{
y(x)= {\sqrt{\alpha} +(-)^k \sqrt{\gamma} \over \sqrt{\alpha}} 
+
\sum_{n=1}^\infty b_nx^n,~~\hbox{ Basic Taylor}
\cr
\cr
\hbox{This is   (\ref{UUU}) when  $a=0$}  
}
&
&
   \matrix{                            \alpha\neq 0 ,
            ~~\sqrt{2\alpha}+(-)^k\sqrt{2\gamma}\not\in{\bf Z}.
\cr 
\cr 
 p_{0x}=-2\cos\pi\Omega_{\alpha\gamma}^k\neq \pm 2.
\cr 
\cr
<M_xM_0,M_1>\hbox{ reducible.}
}
  \cr
&&&
\cr
\hline
&&&\cr
%
%
\hbox{(\ref{TTLO3})} \star_{N=1}&
y(x)= \sum_{n=0}^{|N|-1}b_n x^n+ax^{|N|}+\sum_{n= |N|+1}^\infty b_n(a)x^n
&
      &
            \left.\matrix{\sqrt{2\alpha}+\sqrt{2\gamma} =N \neq 0
\cr
\hbox{ or}
 \cr
                                 \sqrt{2\alpha}-\sqrt{2\gamma} =N \neq 0
}
\right.
\cr
%
%
&    &      &     \cr
       &
b_0={N \over \sqrt{2\alpha}}    \equiv      {\sqrt{2\alpha}+\sqrt{2\gamma}\over \sqrt{2\alpha}} 
\hbox{ or }
                    {\sqrt{2\alpha}-\sqrt{2\gamma}\over \sqrt{2\alpha}} ,~~\alpha\neq 0
      &
            &
                \hbox{\bf and}
\cr
&      &        a       &       \cr
&\hbox{Basic Taylor solution when $N=1$}
      &
          &
\left.
\matrix{
                           \sqrt{2\alpha}      \in     
\left\{\matrix{
                       \{-1,-2,-3,...,N\}~~N<0
\cr
                       \{~1,~2,~3,...,N\}~~~N>0                                                    
                                                      }\right.  
\cr \hbox{ or}\cr
                     \{\sqrt{-2\beta}+\sqrt{1-2\delta}, \sqrt{-2\beta}-\sqrt{1-2\delta} \}
\cr
                    \cap ~{\cal N}_N\neq\emptyset.
\cr 
\cr
 p_{0x}=-2\cos\pi N=\pm 2.
\cr
\cr 
<M_xM_0,M_1>\hbox{ reducible.}
}
\right.
\cr
&&&\cr
\hline
&&&\cr
\hbox{(\ref{TTLO4})} \star
&
           \matrix{y(x)=a+(1-a)(\delta-\beta)x+\sum_{n= 2}^\infty b_n(a)x^n
\cr 
\cr 
\hbox{ Basic Taylor Solution} 
}
      &
           a
                 &
\matrix{\alpha=\gamma=0.
\cr 
 p_{0x}=-2.
\cr 
<M_xM_0,M_1>\hbox{ reducible.}
}
\cr
&&&\cr
\hline
\end{array}
$$


$$
\begin{array}{||c c|c|c||}
\hline\hline
&&&\cr
&
\hbox{{\bf Logarithmic behaviours}}
 & \left.\matrix{\hbox{Free}\cr\hbox{Par.}}\right. & \hbox{Other Conditions}
\cr
\hline
&&&
\cr
\hbox{(\ref{log1})}& 
\matrix{
y(x)=\Sigma_{n=1}^{|N|}b_nx^n
+\Bigl(a+b_{|N|+1}\ln x\Bigr)x^{|N|+1}+ 
\cr 
\cr 
+\Sigma_{n= |N|+2}^\infty P_n(\ln x;a)x^n
\cr 
\cr 
 b_1= {\sqrt{-2\beta}\over N}   \equiv {\sqrt{-2\beta}\over \sqrt{-2\beta}+\sqrt{1-2\delta}} 
\hbox{ or }
     {\sqrt{-2\beta}\over \sqrt{-2\beta}-\sqrt{1-2\delta}} ,~~N\neq 0.
}
                                                               & 
a
                                                               &
\left.
       \matrix{
\sqrt{-2\beta}+\sqrt{1-2\delta} =N
 \cr
                     \hbox{ or} 
 \cr
\sqrt{-2\beta}-\sqrt{1-2\delta}=N.
\cr 
\cr 
2\beta=2\delta-1\hbox{ if } N=0
}
\right.
\cr
&&&\hbox{\bf and}\cr
&&&\cr
\hbox{(\ref{log1zero})} \star 
   &
        y(x)= \Bigl(a\pm\sqrt{-2\beta}\ln x\Bigr)x+\sum_{n=2}^\infty P_n(\ln x;a)x^n,~~N= 0.
    &
    &    
                                                \matrix{
\sqrt{-2\beta}\neq
\cr
 \left\{\matrix{ 0,-1,-2,...,N, \hbox{ if } N \leq 0
\cr
                                         0,~1,~2,...,N, \hbox{ if } N\geq 0
                                         }\right.
    \cr\cr
                                      \sqrt{2\alpha}\pm\sqrt{2\gamma}\not\in{\cal N}_N.
}                                
\cr
&\hbox{Basic solution when $N=0$}&&
\cr 
&&&\matrix{p_{0x}=2\cos\pi N=\pm 2.
\cr 
\cr
<M_0,M_x>\hbox{ reducible.}
}
\cr
&&&\cr
\hline
&&&\cr
\hbox{(\ref{logsquare})} \star
&
y(x)=\left[{2\beta+1-2\delta\over 4}(a+\ln x)^2+{2\beta\over 2\beta+1-2\delta}\right]x+\sum_{n\geq 2}^\infty P_n(\ln x;a)x^n & a & 
2\beta\neq 2\delta-1.
\cr
&\hbox{Basic solution}
&
&
\matrix{p_{0x}=2.
\cr 
\hbox{ no reduc. subgroups.}
}
\cr
\hline
&&&\cr
(\ref{LOG12}) 
     &
        \matrix{
  y(x)= \sum_{n=0}^{|N|-1}b_n x^n+\Bigl(a+b_N\ln x\Bigr)x^{|N|}+
\cr 
\cr 
+\sum_{n= |N|+1}^\infty P_n(\ln x;a)x^n
\cr
\cr
\cr 
b_0={N\over \sqrt{2\alpha}}\equiv {\sqrt{2\alpha}+\sqrt{2\gamma}\over \sqrt{2\alpha}} \hbox{ or } {\sqrt{2\alpha}-\sqrt{2\gamma}\over \sqrt{2\alpha}} ,~~N\neq 0
}
& 
a
&\left.
                                \matrix{
                        \sqrt{2\alpha}+\sqrt{2\gamma} =N\neq 0
\cr
\hbox{ or} \cr
                           \sqrt{2\alpha}-\sqrt{2\gamma} =N\neq 0
\cr                              
\cr
                                {\bf and}
\cr\cr
                            \sqrt{2\alpha}\neq
\cr
\left\{\matrix{0,-1,...,N,\hbox{ if }N\leq -1
\cr
0,~1,...,N, \hbox{ if } N\geq 1
}\right.
\cr
\cr
                 \sqrt{-2\beta}\pm \sqrt{1-2\delta}\not \in {\cal N}_N.
  }
    \right.
\cr
&&&
\cr
&&&
\matrix{p_{0x}=-2\cos \pi N=\pm 2.
\cr
\cr 
<M_xM_0,M_1>\hbox{ reducible.}
}
\cr
&&&
\cr
\hline
\end{array}
$$


$$
\begin{array}{||c c|c|c||} 
\hline\hline
&&&\cr
&
\hbox{{\bf Inverse logarithmic behaviours}}
& \matrix{\hbox{Free}
\cr 
\hbox{Param.}
} & \hbox{ Other Conditions }
\cr
\hline
&&&
\cr
\hbox{(\ref{LOG45})}&
\left.
\matrix{
y(x)= \left\{a\pm \sqrt{2\alpha}\ln x+\sum_{n=1}^\infty P_n(\ln x;a)x^n\right\}^{-1}
\cr
\cr
= \pm{1\over \sqrt{2\alpha}\ln x}\left[1\mp {a\over \sqrt{2\alpha} \ln x}+O\left({1\over \ln^2 x}\right)\right]
}
\right.
&
a
&\matrix{\alpha=\gamma\neq 0.
\cr 
p_{0x}=-2.
\cr 
<M_xM_0,M_1>\hbox{ reducible.}
}
\cr
&&&\cr
\hline
&&&\cr
(\ref{LOG3}) & 
\left.
\matrix{
y(x)= \left\{{\alpha\over \alpha -\gamma}+ {\gamma-\alpha\over 2} (a+\ln x)^2 +\sum_{n=1}^\infty P_{n+1}(\ln x;a) x^n\right\}^{-1}
\cr\cr
=
{2\over (\gamma-\alpha)\ln^2 x}\left[
1-{2a\over \ln x} +O\left({1\over \ln^2 x}\right)
\right]
}
\right.
&
a
&
\matrix{
\alpha\neq \gamma.
\cr 
p_{0x}=-2 .
\cr 
\hbox{no reduc. subgroups.}
}
\cr
&&&\cr
\hline
\end{array}
$$


\section{Table 2:  Critical Behaviours when  $x\to 1$}
\label{TABLE2}


Table 2 provides the critical 
behaviours (\ref{malamente1}) at $x=1$. The branch cut may be taken to be $-\pi <
\hbox{arg}  ~(1-x)\leq \pi$.

\vskip 0.2 cm 
\noindent
$\diamond$ The table is constructed from the table at $x=0$ by Okamoto's transformation (\ref{onara}), in the following way:

   \vskip 0.2 cm
 a) $\hbox{\rm PVI}_{\alpha,\beta,\gamma,\delta}$  is given,  with coefficients $\alpha,\beta,\gamma,\delta$ (or $\theta_0,\theta_x,\theta_1,\theta_\infty$).

\vskip 0.2 cm 
 b) Take $\hbox{\rm PVI}_{  \alpha^\prime,\beta^\prime,\gamma^\prime,\delta^\prime}$   with coefficients 
$$
\alpha^\prime=\alpha,~~~\beta^\prime=-\gamma,~~~\gamma^\prime=-\beta,~~~\delta^\prime=\delta,~~~\hbox{ (or $\theta_0^\prime=\theta_1,~\theta_x^\prime=\theta_x,~\theta_1^\prime=\theta_0,~\theta_\infty^\prime=\theta_\infty$),}
$$
  and variable $\xi$, and compute the critical behaviours $y_0(\xi)$ for $\xi\to 0$. 
The  critical behaviours  at $x=1$ for $\hbox{\rm PVI}_{\alpha,\beta,\gamma,\delta}$  are  then 
$$
y(x)=1-y_0(1-x),~~~~~x\to 1.
$$
This is why the behaviors in Table 2 are numerated as in Table 1, according to the behavior of $y_0(\xi)$ from which they have been  obtained.

\vskip 0.2 cm 
 c) In the table, the integration constants $\sigma$, $a$ , $\nu$, $\phi$, etc appear  in a way that  their  parametrisation in terms of monodromy data can be obtained   from the parametric formulae of Section \ref{dokidoki} for the behaviors at $x=0$ with the same numeration, by the  substitution (\ref{ricordo1}) of  Section \ref{CONNE}.

\vskip 0.3 cm 
\noindent 
$\diamond$ 
Notations: 
\be
\label{SeTG}
\Sigma_{\gamma\delta}^{k}:=\pm \Bigl(\sqrt{2\gamma}+(-)^k\sqrt{1-2\delta}\Bigr),~~~~~\hbox{ if $|\sqrt{2\gamma}+(-)^k\sqrt{1-2\delta}|<1$.}
\ee
otherwise
\be
\label{ikuG}
\Sigma_{\gamma\delta}^{k}:=
 \Bigl(\sqrt{2\gamma}+(-)^k\sqrt{1-2\delta}\Bigr)\hbox{sgn}\Bigl(\Re(\sqrt{2\gamma}+(-)^k\sqrt{1-2\delta})\Bigr)
,~~~k=1,2,
\ee
In the monodromy deformation parameters, $\Sigma_{\gamma\delta}^k
\in\{ (\theta_1+\theta_x),(\theta_1-\theta_x),-(\theta_1+\theta_x),-(\theta_1-\theta_x)\} 
$, or  $\Sigma_{\gamma\delta}^k=
 (\theta_1+(-)^k\theta_x)\hbox{sgn}(\Re(\theta_1+(-)^k\theta_x))
$. 
\be
\label{SeT1G}
\Omega_{\alpha\beta}^{k} :=\pm
\Bigl( \sqrt{2\alpha}+(-)^k\sqrt{-2\beta}\Bigr),~~~~~\hbox{ if  $|\Re \{\sqrt{2\alpha}+(-)^k\sqrt{-2\beta} \}|<1$}. 
\ee
otherwise
\be
\label{STaRG}
\Omega_{\alpha\beta}^{k}:=
\Bigl(\sqrt{2\alpha}+(-)^k\sqrt{-2\beta}\Bigl)\hbox{sgn}\Bigl(\Re(\sqrt{2\alpha}+(-)^k\sqrt{-2\beta})\Bigr)
,~~~k=1,2
\ee
In the monodromy deformation parameters, $\Omega_{\alpha\beta}^k
\in\{(\theta_\infty-1+\theta_0), (\theta_\infty-1-\theta_0),-(\theta_\infty-1+\theta_0), -(\theta_\infty-1-\theta_0)\}
$, or 
$\Omega_{\alpha\beta}^k
=(\theta_\infty-1+(-)^k\theta_0\Bigr)\hbox{sgn}\Bigl(\Re(\theta_\infty-1+(-)^k\theta_0))
$.

\vskip 0.3 cm 
\noindent
$\diamond$  
{\large\bf How to identify a Critical Behavior from given Monodromy Data:}

\vskip 0.2 cm 
Preliminary, observe that when $p_{x1}=-2\cos\pi \Omega_{\alpha\beta}^k$ the Jimbo-Fricke cubic (\ref{FR}) is factorised. Namely, if  $p_{x1}=2\cos\pi(\theta_\infty-\theta_0)$, then (\ref{FR}) is
$$
\left[ p_{01}+p_{0x}e^{i\pi(\theta_\infty-\theta_0)}-2\left(e^{-i\pi\theta_0}\cos\pi\theta_1+e^{i\pi\theta_\infty}\cos\pi\theta_x\right)\right]
\times
$$
\be
\label{era.1}
\times
\left[
p_{01}+p_{x1}e^{-i\pi(\theta_\infty-\theta_0)}-2\left(e^{i\pi\theta_0}\cos\pi\theta_1+e^{-i\pi\theta_\infty}\cos\pi\theta_x\right)
\right]=0
\ee
If $p_{x1}= 2\cos \pi (\theta_\infty+\theta_0)$, then (\ref{FR}) is: 
$$
\left[p_{01}+p_{0x}e^{i\pi(\theta_\infty+\theta_0)}-2\left(e^{i\pi\theta_0}\cos\pi\theta_1+e^{i\pi\theta_\infty}\cos\pi\theta_x\right)
\right]
\times
$$
\be
\label{era1.1}
\times
\left[
p_{01}+p_{0x}e^{-i\pi(\theta_\infty+\theta_0)}-2\left(e^{-i\pi\theta_0}\cos\pi\theta_1+e^{-i\pi\theta_\infty}\cos\pi\theta_x\right)
\right]=0
\ee

\vskip 0.3 cm 
 
\noindent
 The behaviours at $x=1$ are determined by  $p_{x1}$, which gives the monodromy exponent $\sigma=\sigma_{x1}$ in generic cases (see (\ref{fullEXP}.1) and (\ref{atopy}.1)), the monodromy data $p_{0x}, p_{01}$, and the values of the numbers  $\Sigma_{\gamma\delta}^k$ and $\Omega_{\alpha\beta}^k$. The following explains how the critical behaviour is singled out:  
$$
\left[
\matrix{
p_{x1}<-2 & \Rightarrow & \left\{\matrix{(\ref{lantern1}.1) \hbox{ if } p_{x1}\neq -2\cos\pi\Omega_{\alpha,\beta}^k
         \cr (\ref{TAU}.1)  \hbox{ if } p_{x1}= -2\cos\pi\Omega_{\alpha,\beta}^k
}\right.
\cr
\cr
p_{x1}\neq \left\{\matrix {\pm2,-2\cos\pi\Sigma_{\gamma\delta}^k,
\cr
 2\cos\pi\Omega_{\alpha\beta}^k
\cr 
\not<-2
}\right.
& 
\Rightarrow 
&
(\ref{fullEXP}.1)
\cr
\cr
p_{x1}=2\cos\pi\Sigma_{\gamma\delta}^k\neq \pm 2 
&
\Rightarrow
&
(\ref{atopy}.1) \hbox{ [or (\ref{davidekan}.1)]}
\cr 
\cr 
p_{x1}=-2 \cos\pi\Omega_{\alpha\beta}^k\neq \pm 2
&
\Rightarrow
&
\left\{
\matrix{
 \hbox{(\ref{fullEXP}.1),  if  one factor in (\ref{era.1}) or (\ref{era1.1}) $=0$}
\cr 
\hbox{(\ref{UUU}.1) [or (\ref{T1coe}.1)], if $\alpha\neq 0$ and the other factor in (\ref{era.1}) or (\ref{era1.1}) $=0$}
\cr
\hbox{(\ref{div}.1), if $\alpha = 0$ and  the other factor in (\ref{era.1}) or (\ref{era1.1}) $=0$}
}\right.
\cr
\cr
p_{x1}=\pm 2 
&
\Rightarrow
&
\hbox{Taylor or Logarithmic, depending on  $\Sigma_{\gamma\delta}^k,\Omega_{\alpha\beta}^k,\alpha,\beta,\gamma,\delta$
}
}
\right]
$$
 

\vskip 0.3 cm 
\noindent
$\diamond$ 
In the table, the independent variable is   
$$\xi :=1-x \to 0.
$$ 

$$
\begin{array}{||c c|c|c||} 
\hline\hline  &  &  &  \cr 
 & \hbox{{\bf  Complex power behaviours}} & \hbox{Free Par} & \hbox{Other Conditions}    \cr
\hline           & & &  \cr
\hbox{(\ref{fullEXP}.1) } 
  & 
\matrix{
y(x)=1+ \sum_{n=1}^\infty \xi^n\sum_{m=-n}^n c_{nm}(a\xi^{\sigma})^m \sim
\cr \cr 
\sim
\left\{
\matrix{
1+{c_{1,-1}\over a}\xi^{1-\sigma}, & 0<\Re \sigma <1
\cr
\cr 
1+\xi\left[{c_{1,-1}\over a}\xi^{-\sigma}+c_{10}+c_{11}x^\sigma\right], &\Re \sigma =0
}
\right.
\cr 
\cr 
\cr 
c_{11}=-1,~~~c_{10}={1-2\delta-2\gamma-\sigma^2\over 2\sigma^2}
\cr 
\cr
 c_{1,-1}= 
  {[(\sqrt{2\gamma}-\sqrt{1-2\delta})^2-\sigma^2][(\sqrt{2\gamma}+\sqrt{1-2\delta})^2-\sigma^2]
\over 
 -16 \sigma^4 }
\cr 
\cr
\hbox{Basic solutions}
}
& \matrix{
\sigma 
\cr 
\cr
a\neq 0   }                           
                & 
                   
\matrix{       0\leq \Re \sigma<1,~~~\sigma\neq \Sigma_{\gamma\delta}^k,
\cr
2\cos\pi\sigma=p_{x1}.
\cr
\cr
p_{x1}\neq\ 2\cos\pi\Sigma_{\gamma\delta}^k,~\pm 2,
\cr
p_{0x}\not\in(-\infty,-2].
}
\cr     
& &   &   \cr
\hline
&&&\cr
%
%
\hbox{(\ref{atopy}.1)}
 & 
\matrix{y(x)= 1+\sum_{n=1}^\infty \xi^n\sum_{m=0}^n c_{nm} (a\xi^\sigma)^m 
\cr
\cr
\cr 
c_{11}=-1,~~
c_{10}=-{\sqrt{2\gamma}\over\sqrt{2\gamma}+(-)^k\sqrt{1-2\delta}}
\cr 
\cr 
\cr 
\hbox{If $a=0$, $y(x)$ reduces to (\ref{davidekan}.1)}
}
& 
a
&
\left.
\matrix{\Re\sigma>-1,~~
 \sigma=\Sigma_{\beta\gamma}^{k}\not\in{\bf Z},
\cr
\Sigma_{\gamma\delta}^{k} \hbox{ is (\ref{ikuG}) or (\ref{SeTG})}.
\cr
\cr
p_{x1}=2\cos\pi\Sigma_{\gamma\delta}^{k}\neq \pm 2.
\cr
\cr
<M_1,M_x>\hbox{ reducible.}
}
\right.
\cr 
&&&\cr
&\hbox{Basic solutions}&&
\cr
&&& \cr 
\hline &&& \cr 
%
%
\hbox{(\ref{UUU}.1)}
&
\matrix{
y(x)=(-)^{k+1}{\sqrt{-\beta \over \alpha}} + \sum_{n=1}^\infty \xi^n \sum_{m=0}^n d_{nm}(\tilde{a}\xi^{\rho})^m
\cr
\cr
\cr 
\tilde{a}=a\left[{\sqrt{\alpha}+(-)^k\sqrt{-\beta}
\over 
\sqrt{\alpha}}\right]^2,~d_{11}=1
\cr 
\cr 
\cr
\hbox{If $\tilde{a}=0$, $y(x)$ reduces to (\ref{T1coe}.1)}
}
&
a
&\left.
\matrix{\alpha\neq 0.
\cr
\cr 
\rho=\Omega_{\alpha\beta}^{k}-1,~~\Re\rho>-1,
\cr
\Omega_{\alpha\beta}^{k}\not\in{\bf Z},~~
\Omega_{\alpha\beta}^{k} \hbox{ is (\ref{STaRG})}.
\cr
\cr
p_{x1}=-2\cos\pi\Omega_{\alpha\beta}^k\neq \pm 2.
\cr
\cr
<M_1M_x,M_0>\hbox{ reducible.}
}
\right.
\cr
 &&& \cr 
\hline &&& \cr 
%
%
\hbox{(\ref{div}.1)}
&
y(x)=1-{ 1\over a }\xi^{-\omega}
\left(
1+\sum_{n=1}^\infty \xi^n \sum_{m=0}^n d_{nm}(a \xi^{\omega})^m
\right)
&
a
&\left.
\matrix{\alpha =0,
\cr \beta\not\in[0,+\infty),~~\sqrt{-2\beta}\not \in {\bf Z}.
\cr
\cr
\omega=\sqrt{-2\beta}~\hbox{sgn}(\Re\sqrt{-2\beta}),
\cr 
\Re\omega > 0.
\cr
\cr
p_{x1}=-2\cos\pi\sqrt{-2\beta}\neq \pm 2.
\cr
\cr
<M_1M_x,M_0>\hbox{ reducible.}
}
\right.
\cr 
&&&
\cr
\hline
\end{array}
$$


$$
\begin{array}{||c c|c|c||} 
\hline\hline
&&&\cr
&\hbox{{\bf Inverse Oscillatory Behaviours}} & \left.\matrix{\hbox{Free}\cr\hbox{Param}}\right. & \hbox{ Other }
\cr
\hline
&&&
\cr
\hbox{(\ref{lantern1}.1)}
& 
\matrix{
y(x)= 1-\Bigl[\sum_{n=0}^\infty \xi^n\sum_{m=-n-1}^{n+1} c_{nm}
\bigl(e^{i\phi}\xi^{2i\nu}\bigr)^m
\Bigr]^{-1} 
\cr 
\cr 
=1- \left[A\sin(2\nu \ln \xi +\phi)+B+ O(\xi)\right]^{-1}
\cr
 \cr 
\cr 
A=-\sqrt{{\alpha\over 2\nu^2}+B^2},~~
B={2\nu^2-\beta-\alpha\over 4\nu^2},
\cr
& &&
}
&
 \matrix{
\nu
\cr 
\cr 
\phi 
}
& 
\matrix{
\nu\in{\bf R}\backslash\{0\},
~~ 2i\nu\neq \Omega_{\alpha\beta}^k.
\cr
2\cosh2\pi\nu=-p_{x1}.
\cr
\cr
p_{x1}<-2,
\cr
p_{x1}\neq-2\cos\pi\Omega_{\alpha\beta}^k
}
\cr
&&&
\cr
\hline
&&&\cr
%
%
\hbox{(\ref{TAU}.1)}
&
 y(x)=1-
\left[{\sqrt{\alpha} \over \sqrt{\alpha} +(-)^k \sqrt{-\beta}}+a\xi^{-2i\nu}
+\sum_{n=1}^\infty \xi^n
\sum_{m=0}^{n+1}c_{n+1,m}\bigl(a\xi^{-2i\nu}\bigr)^m
\right]^{-1}
&
a
&
\left.
\matrix{2i\nu= \Omega_{\alpha\beta}^{k}\in i{\bf R}\backslash \{0\},
\cr
\Omega_{\alpha\beta}^{k}  \hbox{ is (\ref{SeT1G})}.
\cr
\cr
p_{x1}=-2\cos\pi\Omega_{\alpha\beta}^k<-2.
\cr
\cr
<M_1M_x,M_0>\hbox{ reducible.}
}
\right.
\cr&&&\cr
\hline
\end{array}
$$



$$
\begin{array}{||c c|c|c||}
\hline\hline
&&&\cr
 &\hbox{{\bf Taylor expansions}} & \left.\matrix{\hbox{Free}\cr\hbox{Par}}\right. & \hbox{ Other Conditions}
\cr
\hline
&&&
\cr
\hbox{(\ref{davidekan}.1)}
&
 \matrix{y(x)=1-{\sqrt{2\gamma}\over \sqrt{2\gamma}\pm\sqrt{1-2\delta}}\xi+\sum_{n=2}^\infty b_n\xi^n 
\cr
\cr
 y(x)=0 \hbox{ if } \gamma=0.
\cr
\cr
\hbox{This is (\ref{atopy}.1) when $a=0$}
}
&  
& 
                \matrix{
   \sqrt{2\gamma}+(-)^k\sqrt{1-2\delta}\not\in{\bf Z}.
\cr
\cr
p_{x1}=2\cos\pi\Sigma_{\gamma\delta}^k\neq \pm 2.
\cr
\cr
<M_1,M_x>\hbox{ reducible.}
}
\cr
&&&\cr
\hline
%
%
&&&
\cr
\hbox{(\ref{TTLO1}.1)}
&
y(x)= 1+\sum_{n=1}^{|N|}b_n\xi^n-a\xi^{|N|+1}+\sum_{n= |N|+2}^\infty b_n(a)\xi^n
&
&
                           \left.
                                        \matrix{\sqrt{2\gamma}+\sqrt{1-2\delta} =N\neq 0
\cr
\hbox{ or} 
\cr
                                                       \sqrt{2\gamma}-\sqrt{1-2\delta}=N\neq 0
}
\right.
\cr
                 &               &              &           \cr
&        b_1=-{\sqrt{2\gamma}\over N}\equiv -{\sqrt{2\gamma}\over \sqrt{2\gamma}+\sqrt{1-2\delta}}
\hbox{ or }
  -{\sqrt{2\gamma}\over \sqrt{2\gamma}-\sqrt{1-2\delta}}              
   &
 a
 &                                                       
 \hbox{\bf and}
\cr
%
%
%
             &       
&
   &
      \left.\matrix{
                           \sqrt{2\gamma}\in \left\{\matrix{
 \{0,-1,-2,...,N\}~~N<0
\cr
    \{0,~1,~2,...,N\}~~~N>0                                                    
                                                      }\right.           
 \cr
\hbox{ or }
\cr
                                  \{ \sqrt{2\alpha}+\sqrt{-2\beta},\sqrt{2\alpha}-\sqrt{-2\beta}\} 
\cap {\cal N}_N \neq \emptyset.
\cr
\cr
p_{x1}=2\cos\pi N=\pm 2.
\cr
\cr
<M_1,M_x>\hbox{ reducible.}
}
\right.
\cr
&&&\cr
\hline
&&&\cr
\hbox{(\ref{TTLO2}.1)}                
&       
                                                       y(x)= 1-a\xi+a(a-1)\left(\beta+\alpha+{1\over 
                                                       2}\right)\xi^2+\sum_{n= 3}^\infty b_n(a)\xi^n
&
        a
& 
\matrix{
2\gamma=1-2\delta=0.
\cr
p_{x1}=2.
\cr
<M_1,M_x>\hbox{ reducible.}
}
\cr
&&&\cr
\hline
&&&\cr
%
%
\hbox{(\ref{T1coe}.1)}
&
\matrix{
y(x)= (-)^{k+1}{\sqrt{-\beta \over \alpha}} 
+
\sum_{n=1}^\infty b_n\xi^n,~~\hbox{ Basic Taylor}
\cr
\cr
\hbox{This is  (\ref{UUU}.1) when  $a=0$}  
}
&
&
\matrix{
                               \alpha\neq 0 ,
            ~~\sqrt{2\alpha}+(-)^k\sqrt{-2\beta}\not\in{\bf Z}.
\cr
p_{x1}=-2\cos\pi\Omega_{\alpha\beta}^k\neq \pm 2.
\cr
<M_1M_x,M_0>\hbox{ reducible.}
}
 
 \cr
&&&
\cr
\hline
&&&\cr
%
%
\hbox{(\ref{TTLO3}.1)}
&
y(x)= \sum_{n=0}^{|N|-1}b_n \xi^n-a\xi^{|N|}+\sum_{n= |N|+1}^\infty b_n(a)\xi^n
&
      &
            \left.\matrix{\sqrt{2\alpha}+\sqrt{-2\beta} =N \neq 0
\cr
\hbox{ or}
 \cr
                                 \sqrt{2\alpha}-\sqrt{-2\beta} =N \neq 0
}
\right.
\cr
%
%
&    &      &     \cr
       &
b_0=1-{N \over \sqrt{2\alpha}}    \equiv      -{\sqrt{-\beta\over \alpha}} 
\hbox{ or }
                    {\sqrt{-\beta\over\alpha}} ,~~\alpha\neq 0
      &
            &
                \hbox{\bf and}
\cr
&      &        a       &       \cr
&\hbox{Basic Taylor solution when $N=1$}
      &
          &
\left.
\matrix{
                           \sqrt{2\alpha}      \in     
\left\{\matrix{
                       \{-1,-2,-3,...,N\}~~N<0
\cr
                       \{~1,~2,~3,...,N\}~~~N>0                                                    
                                                      }\right.  
\cr \hbox{ or}\cr
                     \{\sqrt{2\gamma}+\sqrt{1-2\delta}, \sqrt{2\gamma}-\sqrt{1-2\delta} \}
\cr
                    \cap ~{\cal N}_N\neq\emptyset.
\cr
\cr
p_{x1}=-2\cos\pi N=\pm 2.
\cr
\cr
<M_1M_x,M_0>\hbox{ reducible.}
}
\right.
\cr
&&&\cr
\hline
&&&\cr
\hbox{(\ref{TTLO4}.1)}
&
          \matrix{ y(x)=(1-a)-(1-a)(\delta+\gamma)\xi+\sum_{n= 2}^\infty b_n(a)\xi^n
\cr
\cr
\hbox{ Basic Taylor Solution} 
}
      &
           a
                 &
\matrix{
\alpha=\beta=0.
\cr
p_{x1}=-2.
\cr
<M_1M_x,M_0>\hbox{ reducible.}
}
\cr
&&&\cr
\hline
\end{array}
$$


$$
\begin{array}{||c c|c|c||}
\hline\hline
&&&\cr
&
\hbox{{\bf Logarithmic behaviours}}
 & \left.\matrix{\hbox{Free}\cr\hbox{Par.}}\right. & \hbox{Other Conditions}
\cr
\hline
&&&
\cr
\hbox{(\ref{log1}.1)}
& 
y(x)=1+\Sigma_{n=1}^{|N|}b_n\xi^n
+\Bigl(b_{|N|+1}\ln \xi - a \Bigr)\xi^{|N|+1}+ 
                                                               & 
                                                               &
\left.
       \matrix{
\sqrt{2\gamma}+\sqrt{1-2\delta} =N
 \cr
                     \hbox{ or} 
 \cr
\sqrt{2\gamma}-\sqrt{1-2\delta}=N
}
\right.
                                                                  \cr
&
+\Sigma_{n= |N|+2}^\infty P_n(\ln \xi;a)\xi^n & &
\left.
\matrix{
\hbox{ }
\cr
2\gamma=1-2\delta\hbox{ if } N=0
}
\right.
\cr
&
 &
   &
     \cr
&
        b_1= -{\sqrt{2\gamma}\over N}   \equiv -{\sqrt{2\gamma}\over \sqrt{2\gamma}+\sqrt{1-2\delta}} 
\hbox{ or }
    - {\sqrt{2\gamma}\over \sqrt{2\gamma}-\sqrt{1-2\delta}} ,~~N\neq 0.
   &
a     
 &
       \hbox{\bf and}\cr
&&&\cr
\hbox{(\ref{log1zero}.1)}
   &
        y(x)=1- \Bigl(a\pm\sqrt{2\gamma}\ln \xi\Bigr)\xi+\sum_{n=2}^\infty P_n(\ln \xi;a)\xi^n,~~N= 0.
    &
    &    
                                                \matrix{
\sqrt{2\gamma}\neq
\cr
 \left\{\matrix{ 0,-1,-2,...,N, \hbox{ if } N \leq 0
\cr
                                         0,~1,~2,...,N, \hbox{ if } N\geq 0
                                         }\right.
    \cr\cr
                                      \sqrt{2\alpha}\pm\sqrt{-2\beta}\not\in{\cal N}_N.
\cr
\cr
p_{x1}=2\cos\pi N=\pm 2.
\cr
\cr
<M_1,M_x>\hbox{ reducible.}
}                                
\cr
&\hbox{Basic solution when $N=0$}&&
\cr
&&&\cr
\hline
&&&\cr
\hbox{(\ref{logsquare}.1)}
&
\matrix{
y(x)=1-\left[{1-2\gamma-2\delta\over 4}(a+\ln \xi)^2+{2\gamma\over 2\gamma+2\delta-1}\right]\xi+\sum_{n\geq 2}^\infty P_n(\ln \xi;a)\xi^n 
\cr
\cr
\hbox{Basic solution}
}
& a & 
\matrix{
2\gamma\neq 1-2\delta.
\cr
p_{x1}=2.
\cr
\hbox{no reduc. subrgroups.}
}
\cr
&&&\cr
\hline
&&&\cr
(\ref{LOG12}.1) 
     &
        \matrix{
  y(x)=\sum_{n=0}^{|N|-1}b_n \xi^n+\Bigl(b_N\ln \xi-a\Bigr)x^{|N|}+\sum_{n= |N|+1}^\infty P_n(\ln \xi)\xi^n
\cr
\cr
b_0=1-{N\over \sqrt{2\alpha}}\equiv -\sqrt{-\beta\over \alpha} \hbox{ or } \sqrt{-\beta\over \alpha} 
}
& 
a
&\left.
                                \matrix{
                        \sqrt{2\alpha}+\sqrt{-2\beta} =N\neq 0
\cr
\hbox{ or} \cr
                           \sqrt{2\alpha}-\sqrt{-2\beta} =N\neq 0
\cr                              
\cr
                                {\bf and}
\cr\cr
                            \sqrt{2\alpha}\neq
\cr
\left\{\matrix{0,-1,...,N,\hbox{ if }N\leq -1
\cr
0,~1,...,N, \hbox{ if } N\geq 1
}\right.
\cr
\cr
                 \sqrt{2\gamma}\pm \sqrt{1-2\delta}\not \in {\cal N}_N.
\cr
\cr
p_{x1}=-2\cos\pi N=\pm 2.
\cr
\cr
<M_1M_x,M_0>\hbox{ reducible.}
  }
    \right.
\cr
&&&\cr
\hline
\end{array}
$$


$$
\begin{array}{||c c|c|c||} 
\hline\hline
&&&\cr
&
\hbox{{\bf Inverse  logarithmic behaviours}}
& \hbox{Free Parameters} & \hbox{ Other Conditions }
\cr
\hline
&&&
\cr
\hbox{(\ref{LOG45}.1)}
&
\left.
\matrix{
y(x)=1- \left\{a\pm \sqrt{2\alpha}\ln x+\sum_{n=1}^\infty P_n(\ln x;a)x^n\right\}^{-1}
\cr
\cr
= 1\mp{1\over \sqrt{2\alpha}\ln x}\left[1\mp {a\over \sqrt{2\alpha} \ln x}+O\left({1\over \ln^2 x}\right)\right]
}
\right.
&
a
&\matrix{
\alpha=-\beta\neq0.
\cr
p_{x1}=-2.
\cr
<M_1M_x,M_0>\hbox{ reducible.}
}
\cr
&&&\cr
\hline
&&&\cr
(\ref{LOG3}.1) 
& 
\left.
\matrix{
y(x)=1- \left\{{\alpha\over \alpha +\beta}- {\alpha+\beta\over 2} (a+\ln \xi)^2 +\sum_{n=1}^\infty P_{n+1}(\ln \xi;a) \xi^n\right\}^{-1}
\cr\cr
=
1+{2\over (\alpha+\beta)\ln^2 \xi}\left[
1-{2a\over \ln \xi} +O\left({1\over \ln^2 \xi}\right)
\right]
}
\right.
&
a
&
\matrix{
\alpha\neq -\beta.
\cr
p_{x1}=-2.
\cr
\hbox{no reduc. subgroups.}
}
\cr
&&&\cr
\hline
\end{array}
$$


\section{Table 3:  Critical behaviours when $x\to \infty$}
\label{TABLE3}


Table 3 provides the critical 
behaviours (\ref{malamente1}) at $x=\infty$. The branch cut may be taken to be $-\pi \leq 
\hbox{arg}  x < \pi$.

\vskip 0.2 cm 
\noindent
$\diamond$ 
The table is constructed from the table at $x=0$ by Okamoto's transformation (\ref{onara1}), in the following way:

\vskip 0.2 cm 
a) $\hbox{\rm PVI}_{ \alpha,\beta,\gamma,\delta}$ is given, with coefficients $\alpha,\beta,\gamma,\delta$ (or $\theta_0,\theta_x,\theta_1,\theta_\infty$).

\vskip 0.2 cm 
 b) Take $\hbox{\rm PVI}_{ \alpha^\prime,\beta^\prime,\gamma^\prime,\delta^\prime}$ with coefficients 
$$
\alpha^\prime=\alpha,~~~\beta^\prime=\beta,~~~\gamma^\prime={1\over 2}-\delta,~~~\delta^\prime={1\over 2}-\gamma,
~~~\hbox{ (or $\theta_0^\prime=\theta_0,~\theta_x^\prime=\theta_1,~\theta_1^\prime=\theta_x,~\theta_\infty^\prime=\theta_\infty$ ),}
$$
  and new variable $\xi$, and compute the critical behaviours $y_0(\xi)$ for $\xi\to 0$. 
The   behaviours at $x=\infty$ for $\hbox{\rm PVI}_{ \alpha,\beta,\gamma,\delta}$ are  then 
$$
y(x)=x~y_0\left({1\over x}\right),~~~~~x\to\infty
$$
Therefore, the numeration of the behaviors in Table 3 corresponds to the behavior of $y_0(\xi)$ at $\xi=0$ tabulated in Table 1.

\vskip 0.2 cm 

c) In the table, the integration constants $\sigma$, $a$ , $\nu$, $\phi$, etc appear  in a way that  their  parametrisation in terms of monodromy data can be obtained   from the  parametric formulae at zero of the behaviors with the same numeration, by  the substitution (\ref{ricordo2}) of Section \ref{CONNE}.

\vskip 0.3 cm 
\noindent 
$\diamond$  
Notations:
\be
\label{SeTD}
\Sigma_{\beta\gamma}^{k}:=\pm
\Bigl(\sqrt{-2\beta}+(-)^k\sqrt{2\gamma}\Bigr),~~~~~~\hbox{  if $|\sqrt{-2\beta}+(-)^k\sqrt{2\gamma}|<1$},
\ee
otherwise
\be
\label{ikuD}
\Sigma_{\beta\gamma}^{k}:=
 \Bigl(\sqrt{-2\beta}+(-)^k\sqrt{2\gamma}\Bigr)\hbox{sgn}\Bigl(\Re(\sqrt{-2\beta}+(-)^k\sqrt{2\gamma})\Bigr)
,~~~k=1,2,
\ee
In the monodromy deformation parameters, $\Sigma_{\beta\gamma}^k
\in \{ (\theta_0+\theta_1),(\theta_0-\theta_1),-(\theta_0+\theta_1),-(\theta_0-\theta_1)\} 
$, or  $\Sigma_{\beta\gamma}^k=
 (\theta_0+(-)^k\theta_1)\hbox{sgn}(\Re(\theta_0+(-)^k\theta_1))
$.
\be
\label{SeT1D}
\Omega_{\alpha\delta}^{k}:= \pm 
\Bigl( \sqrt{2\alpha}+(-)^k\sqrt{1-2\delta}\Bigr),~~~~~\hbox{ if  $|\Re \{\sqrt{2\alpha}+(-)^k\sqrt{1-2\delta} \}|<1$}. 
\ee
otherwise
\be
\label{STaRD}
\Omega_{\alpha\delta}^{k}:=
\Bigl(\sqrt{2\alpha}+(-)^k\sqrt{1-2\delta}\Bigl)\hbox{sgn}\Bigl(\Re(\sqrt{2\alpha}+(-)^k\sqrt{1-2\delta})\Bigr)
,~~~k=1,2
\ee
In the monodromy deformation parameters, $\Omega_{\alpha\delta}^k
\in\{(\theta_\infty-1+\theta_x), (\theta_\infty-1-\theta_x),-(\theta_\infty-1+\theta_x), -(\theta_\infty-1-\theta_x)\}
$, or 
$\Omega_{\alpha\delta}^k
=(\theta_\infty-1+(-)^k\theta_x\Bigr)\hbox{sgn}\Bigl(\Re(\theta_\infty-1+(-)^k\theta_x))
$.

\vskip 0.3 cm 
\noindent
$\diamond$  
{\large\bf How to identify a Critical Behavior from given Monodromy Data:}

\vskip 0.2 cm 
Preliminary, observe that when  $p_{01}=-2\cos\pi \Omega_{\alpha\delta}^k$ the Jimbo-Fricke cubic (\ref{FR}) is factorised. Namely, if  $p_{01}=2\cos\pi(\theta_\infty-\theta_x)$, then (\ref{FR}) is
$$
\left[ p_{0x}+p_{x1}e^{i\pi(\theta_x-\theta_\infty)}-2\left(e^{i\pi\theta_x}\cos\pi\theta_0
+e^{-i\pi\theta_\infty}\cos\pi\theta_1\right)\right]
\times
$$
\be
\label{era.in}
\times
\left[
p_{0x}+p_{x1}e^{i\pi(\theta_\infty-\theta_x)}-2\left(e^{-i\pi\theta_x}\cos\pi\theta_0+e^{i\pi\theta_\infty}\cos\pi\theta_1\right)
\right]=0
\ee
If $p_{01}= 2\cos \pi (\theta_\infty+\theta_x)$, then (\ref{FR}) is: 
$$
\left[p_{0x}+p_{x1}e^{-i\pi(\theta_\infty+\theta_x)}-2\left(e^{-i\pi\theta_x}\cos\pi\theta_0+e^{-i\pi\theta_\infty}\cos\pi\theta_1\right)
\right]
\times
$$
\be
\label{era1.in}
\times
\left[
p_{0x}+p_{x1}e^{i\pi(\theta_\infty+\theta_x)}-2\left(e^{i\pi\theta_x}\cos\pi\theta_0+e^{i\pi\theta_\infty}\cos\pi\theta_1\right)
\right]=0
\ee

\vskip 0.3 cm 
 
\noindent
The behaviours at $x=\infty$ are determined by  $p_{01}$, which gives the monodromy exponent $\sigma=\sigma_{01}$ in generic cases (see (\ref{fullEXP}.$\infty$) and (\ref{atopy}.$\infty$)), the monodromy data $p_{0x}, p_{x1}$, and the values of the numbers  $\Sigma_{\beta\gamma}^k$ and $\Omega_{\alpha\delta}^k$. The following explains how the critical behaviour is singled out:  
$$
\left[
\matrix{
p_{01}<-2 & \Rightarrow & \left\{\matrix{(\ref{lantern1}.\infty) \hbox{ if } p_{01}\neq -2\cos\pi\Omega_{\alpha\delta}^k
         \cr (\ref{TAU}.\infty)  \hbox{ if } p_{01}= -2\cos\pi\Omega_{\alpha\delta}^k
}\right.
\cr
\cr
p_{01}\neq \left\{\matrix {\pm2,-2\cos\pi\Sigma_{\beta\gamma}^k,
\cr
 2\cos\pi\Omega_{\alpha\delta}^k
\cr
\not <-2
}\right.
& 
\Rightarrow 
&
(\ref{fullEXP}.\infty)
\cr
\cr
p_{01}=2\cos\pi\Sigma_{\beta\gamma}^k\neq \pm 2 
&
\Rightarrow
&
(\ref{atopy}.\infty) \hbox{ [or (\ref{davidekan}.$\infty$)]}
\cr 
\cr 
p_{01}=-2 \cos\pi\Omega_{\alpha\delta}^k\neq \pm 2
&
\Rightarrow
&
\left\{
\matrix{
 \hbox{(\ref{fullEXP}.$\infty$),  if  one factor in (\ref{era.in}) or (\ref{era1.in}) $=0$}
\cr 
\hbox{(\ref{UUU}.$\infty$) [or (\ref{T1coe}.$\infty$)], if $\alpha\neq 0$ and the other factor in (\ref{era.in}) or (\ref{era1.in}) $=0$}
\cr
\hbox{(\ref{div}.$\infty$), if $\alpha = 0$ and  the other factor in (\ref{era.in}) or (\ref{era1.in}) $=0$}
}\right.
\cr
\cr
p_{01}=\pm 2 
&
\Rightarrow
&
\hbox{Taylor or Logarithmic, depending on  $\Sigma_{\beta\gamma}^k,\Omega_{\alpha\delta}^k,\alpha,\beta,\gamma,\delta$
}
}
\right]
$$
 



$$
\begin{array}{||c c|c|c||} 
\hline\hline   & &  &  \cr 
  &\hbox{{\bf  Complex power behaviours}} & \hbox{Free Par.} & \hbox{Other Conditions}    \cr
\hline           & & &  \cr
\hbox{(\ref{fullEXP}.$\infty$)}
&
 \matrix{y(x)= \sum_{n=1}^\infty x^{1-n}\sum_{m=-n}^n c_{nm}(ax^{-\sigma})^m \sim_{x\to\infty}
\cr
\cr
 \sim \left\{
\matrix{
{c_{1,-1}\over a} x^\sigma, & \Re\sigma>0
\cr 
{c_{1,-1}\over a} x^{-\sigma}+c_{10}+ax^\sigma,& \Re\sigma=0
}
\right.
\cr 
\cr 
\cr 
 c_{11}=1,~~c_{10}={\sigma^2-2\beta-2\gamma\over 2\sigma^2} 
\cr 
c_{1,-1}=
{\Bigl[(\sqrt{-\beta}-\sqrt{\gamma})^2-{\sigma^2\over 2}\Bigr]  \Bigl[(\sqrt{-\beta}+\sqrt{\gamma})^2-
{\sigma^2\over 2}\Bigr]\over 4 \sigma^4 }
}
&
\matrix{
 \sigma 
\cr 
\cr 
a\neq 0 
}
                                  & 
\matrix{
                          0\leq \Re \sigma<1,~~~\sigma\neq \Sigma_{\beta\gamma}^k,
\cr
2\cos\pi\sigma=p_{01}.
\cr
\cr
p_{01}\neq \pm 2,~2\cos\pi\Sigma_{\beta\gamma}^k,
\cr
p_{01}\not\in(-\infty,-2]
}
\cr
&&&\cr
&\hbox{Basic solutions} &   &   \cr
\hline
&&&\cr
%
%
\hbox{(\ref{atopy}.$\infty$)} 
& 
\matrix{
y(x)= \sum_{n=1}^\infty x^{1-n}\sum_{m=0}^n c_{nm} (ax^{-\sigma})^m 
\cr
\cr
\cr
y(x)\sim
\left\{
\matrix{
c_{10},&\Re\sigma>0
\cr
c_{10}+ax^{-\sigma},& \Re\sigma=0
\cr 
ax^{-\sigma},&\Re\sigma<0
}\right\}~~~\beta\neq 0
\cr
y(x)\sim x^{-\sigma} \to 
\left\{
\matrix{
0, & \sigma>0
\cr 
\infty, &-1< \sigma<0 
\cr 
\hbox{oscillates}, &\Re\sigma=0
}
\right\}
~\matrix{
\beta = 0\cr 
\Downarrow
\cr
\sigma= \pm \sqrt{2\gamma}
}
\cr
\cr 
\cr
c_{11}=1,~~
c_{10}={\sqrt{-\beta}\over\sqrt{-\beta}+(-)^k\sqrt{\gamma}}
\cr 
\cr 
\cr 
\hbox{If $a=0$, $y(x)$ reduces to $(\ref{davidekan}.\infty)$}
}
& 
a
&
\matrix{\Re\sigma>-1,~~
 \sigma=\Sigma_{\beta\gamma}^{k}\not\in{\bf Z},
\cr 
\Sigma_{\beta\gamma}^{k} \hbox{ is (\ref{ikuD}) or (\ref{SeTD})}.
\cr
\cr
p_{01}=2\cos\pi \Sigma_{\beta\gamma}^{k}\neq \pm 2.
\cr
\cr
<M_0,M_1>\hbox{ reducible.}
}
\cr 
&&&\cr
&\hbox{Basic solutions}&&
\cr
&&& \cr 
\hline &&& \cr 
%
%
\hbox{(\ref{UUU}.$\infty$)}
&
\matrix{
y(x)=x\left[d_{00} + \sum_{n=1}^\infty x^{-n} \sum_{m=0}^n d_{nm}(\tilde{a}x^{-\rho})^m\right]\sim d_{00}x
\cr
\cr
\cr 
d_{00}={\sqrt{2\alpha} +(-)^k\sqrt{1-2\delta} \over \sqrt{2\alpha}},~~
\tilde{a}=-a~d_{00}^2,~~d_{11}=1
\cr
\cr
\cr
\hbox{If $\tilde{a}=0$, $y(x)$ reduces to $(\ref{T1coe}.\infty)$}
}
&
a
&\left.
\matrix{\alpha\neq 0.
\cr
\cr
\rho=\Omega_{\alpha\delta}^{k}-1,~~
\Re\rho>-1,  
\cr
\Omega_{\alpha\delta}^{k}\not\in{\bf Z},~~
\Omega_{\alpha\delta}^{k} \hbox{ is (\ref{STaRD})}.
\cr
\cr
p_{01}=-2\cos\pi\Omega_{\alpha\delta}^{k}\neq \pm 2.
\cr
\cr
<M_0M_1,M_x>\hbox{ reducible.}
}
\right.
\cr
 &&& \cr 
\hline &&& \cr 
%
%
\hbox{(\ref{div}.$\infty$)}
&
y(x)={ 1\over a }x^{1+\omega}
\left(
1+\sum_{n=1}^\infty x^{-n} \sum_{m=0}^n d_{nm}(a x^{-\omega})^m
\right)
&
a
&\left.
\matrix{\alpha =0, 
\cr \delta\not\in\left[{1\over 2},+\infty\right),~~\sqrt{1-2\delta}\not \in {\bf Z}.
\cr  
\cr 
\omega=\sqrt{1-2\delta}~\hbox{sgn}(\Re\sqrt{1-2\delta}),
\cr
\Re\omega > 0.
\cr
\cr
p_{01}=-2\cos\pi\sqrt{1-2\delta}\neq\pm 2.
\cr
\cr
<M_0M_1,M_x>\hbox{ reducible.}
}
\right.
\cr 
&&&
\cr
\hline
\end{array}
$$


$$
\begin{array}{||c c|c|c||} 
\hline\hline
&&&\cr
&\hbox{{\bf Inverse Oscillatory Behaviours}} & \left.\matrix{\hbox{Free}\cr\hbox{Param}}\right. & \hbox{ Other }
\cr
\hline
&&&
\cr
\hbox{(\ref{lantern1}.$\infty$)}
& 
\matrix{
y(x)= x\Bigl[\sum_{n=0}^\infty x^{-n}\sum_{m=-n-1}^{n+1} c_{nm}
\bigl(e^{i\phi}x^{-2i\nu}\bigr)^m
\Bigr]^{-1} 
\cr 
\cr 
= x\left[B-A\sin(2\nu \ln x -\phi)+ O\left({1\over x}\right)\right]^{-1}
\cr 
\cr 
 A=-\sqrt{{\alpha\over 2\nu^2}+B^2},~~
~~~B={4\nu^2+1-2\delta-2\alpha\over 8\nu^2},
}
&
\matrix{
 \nu
\cr 
\cr
\phi} 
&
\matrix{
\nu\in{\bf R}\backslash\{0\},~~ 2i\nu\neq\Omega_{\alpha\delta}^k,
\cr
2\cosh 2\pi \nu=-p_{01}.
\cr
\cr
p_{01}<-2,
\cr
p_{01}\neq-2\cos\pi\Omega_{\alpha\delta}^k.
}

\cr
&&&
\cr
\hline
&&&\cr
%
%
\hbox{(\ref{TAU}.$\infty$)}
& 
\matrix{
y(x)=x\left[c_0+ax^{2i\nu}+\sum_{n=1}^\infty x^{-n}
\sum_{m=0}^{n+1}c_{n+1,m}\bigl(ax^{2i\nu}\bigr)^m
\right]^{-1}
\cr 
\cr 
c_0={\sqrt{2\alpha} \over \sqrt{2\alpha} +(-)^k \sqrt{1-2\delta}}
}
&
a
&
\left.
\matrix{2i\nu= \Omega_{\alpha\delta}^{k}\in i {\bf R}\backslash\{0\},
\cr
\Omega_{\alpha\delta}^{k}  \hbox{ is (\ref{SeT1D})}
\cr 
\cr 
p_{01}=-2\cos\pi\Omega_{\alpha\delta}^k<-2.
\cr
\cr
<M_0M_1,M_x>\hbox{ reducible.}
}
\right.
\cr&&&\cr
\hline
\end{array}
$$



$$
\begin{array}{||c c|c|c||}
\hline\hline
&&&\cr
 &\hbox{{\bf  Expansions in $x^{-1}$}} & \left.\matrix{\hbox{Free}\cr\hbox{Par}}\right. & \hbox{ Other Conditions}
\cr
\hline
&&&
\cr
\hbox{(\ref{davidekan}.$\infty$)}
&
\matrix{ 
y(x)={\sqrt{-\beta}\over \sqrt{-\beta}+(-)^k\sqrt{\gamma}}+\sum_{n=1}^\infty b_nx^{-n} 
\cr
\cr
y(x)=0 \hbox{ if } \beta=0.
\cr
\cr
\hbox{This is $(\ref{atopy}.\infty)$  when $a=0$.}
}
&
& 
               \matrix{
    \sqrt{-2\beta}+(-)^k\sqrt{2\gamma}\not\in{\bf Z}.
\cr
\cr
p_{01}=2\cos\pi \Sigma_{\beta\gamma}^k\neq\pm 2.
\cr
\cr
<M_0,M_1>\hbox{ reducible.}
}
\cr
&&& \cr
\hline
%
%
&&&
\cr
\hbox{(\ref{TTLO1}.$\infty$)}
&
y(x)= \sum_{n=0}^{|N|-1}b_nx^{-n}+ax^{-|N|}+\sum_{n= |N|+1}^\infty b_n(a)x^{-n}
&
&
                           \left.
                                        \matrix{\sqrt{-2\beta}+\sqrt{2\gamma} =N\neq 0
\cr\hbox{ or} \cr
                                                       \sqrt{-2\beta}-\sqrt{2\gamma}=N\neq 0
}
\right.
\cr
                 &               &              &           \cr
&        b_0={\sqrt{-2\beta}\over N}\equiv {\sqrt{-\beta}\over \sqrt{-\beta}+\sqrt{\gamma}}\hbox{ or } {\sqrt{-\beta}\over \sqrt{-\beta}-\sqrt{\gamma}}                 & a &                                                        \hbox{\bf and}
\cr
%
%
%
             &       
&
   &
      \left.\matrix{
                           \sqrt{-2\beta}\in \left\{\matrix{
 \{0,-1,-2,...,N\}~~N<0
\cr
    \{0,~1,~2,...,N\}~~~N>0                                                    
                                                      }\right.           
 \cr
\hbox{ or }
\cr
                                  \{ \sqrt{2\alpha}+\sqrt{1-2\delta},\sqrt{2\alpha}-\sqrt{1-2\delta}\} 
\cap
\cr
\cap {\cal N}_N \neq \emptyset.
\cr
\cr
p_{01}=2\cos\pi N=\pm 2.
\cr
\cr
<M_0,M_1>\hbox{ reducible.}
}
\right.
\cr
&&&\cr
\hline
&&&\cr
\hbox{(\ref{TTLO2}$.\infty$)}                
&       
                                                       y(x)= 1+a(a-1)\left(\gamma-\alpha-{1\over 2}\right){1\over x}+\sum_{n= 2}^\infty b_n(a)x^{-n}
&
        a
& 
\matrix{
\beta=\gamma=0.
\cr
p_{01}=2.
\cr
<M_0,M_1>\hbox{ reducible.}
}
\cr
&&&\cr
\hline
&&&\cr
%
%
\hbox{(\ref{T1coe}.$\infty$)}
&
\matrix{
y(x)= {\sqrt{2\alpha} +(-)^k \sqrt{1-2\delta} \over \sqrt{2\alpha}}~x 
+
\sum_{n=0}^\infty b_nx^{-n},~~\hbox{ Basic Taylor}
\cr
\cr
\hbox{This is  $(\ref{UUU}.\infty)$ when  $a=0$}  
}
&
&
                      \matrix{
         \alpha\neq 0 ,
            ~~\sqrt{2\alpha}+(-)^k\sqrt{1-2\delta}\not\in{\bf Z}.
\cr
p_{01}=-2\cos\pi \Omega_{\alpha\delta}^k\neq \pm 2.
\cr
<M_0M_1,M_x>\hbox{ reducible.}
}
  \cr
&&&
\cr
\hline
&&&\cr
%
%
\hbox{(\ref{TTLO3}.$\infty$)}
&
y(x)= b_\infty x+\sum_{n=0}^{|N|-2}b_n x^{-n}+ax^{1-|N|}+\sum_{n= |N|}^\infty b_n(a)x^{-n}
&
      &
            \left.\matrix{\sqrt{2\alpha}+\sqrt{1-2\delta} =N \neq 0
\cr
\hbox{ or}
 \cr
                                 \sqrt{2\alpha}-\sqrt{1-2\delta} =N \neq 0
}
\right.
\cr
%
%
&    &      &     \cr
       &
b_\infty={N \over \sqrt{2\alpha}}    \equiv      {\sqrt{2\alpha}+\sqrt{1-2\delta}\over \sqrt{2\alpha}} 
\hbox{ or }
                    {\sqrt{2\alpha}-\sqrt{1-2\delta}\over \sqrt{2\alpha}} ,~~\alpha\neq 0
      &
            &
                \hbox{\bf and}
\cr
&      &        a       &       \cr
&\hbox{Basic Taylor solution when $N=1$}
      &
          &
\left.
\matrix{
                           \sqrt{2\alpha}      \in     
\left\{\matrix{
                       \{-1,-2,-3,...,N\}~~N<0
\cr
                       \{~1,~2,~3,...,N\}~~~N>0                                                    
                                                      }\right.  
\cr \hbox{ or}\cr
                     \{\sqrt{-2\beta}+\sqrt{2\gamma}, \sqrt{-2\beta}-\sqrt{2\gamma} \}
\cap ~{\cal N}_N\neq\emptyset.
\cr
\cr
p_{01}=-2\cos\pi N=\pm 2.
\cr
\cr
<M_0M_1,M_x>\hbox{ reducible.}
}
\right.
\cr
&&&\cr
\hline
&&&\cr
\hbox{(\ref{TTLO4}.$\infty$)}
&
\matrix{
           y(x)=a~x+(1-a)\left({1\over 2}-\gamma-\beta\right)+\sum_{n= 1}^\infty b_n(a)x^{-n}
\cr
\cr
\hbox{ Basic Taylor Solution} .
}
      &
           a
                 &
\matrix{
2\alpha=1-2\delta=0.
\cr
p_{01}=-2.
\cr
<M_0M_1,M_x>\hbox{ reducible.}
}
\cr
&&&\cr
\hline
\end{array}
$$


$$
\begin{array}{||c c|c|c||}
\hline\hline
&&&\cr
&
\hbox{{\bf Logarithmic behaviours}}
 & \left.\matrix{\hbox{Free}\cr\hbox{Param.}}\right. & \hbox{Other Conditions}
\cr
\hline
&&&
\cr
\hbox{(\ref{log1}.$\infty$)}
& 
y(x)=\Sigma_{n=0}^{|N|-1}b_nx^{-n}
+\Bigl(a+b_{|N|+1}\ln x\Bigr)x^{-|N|}+ 
                                                               & 
                                                               &
\left.
       \matrix{
\sqrt{-2\beta}+\sqrt{2\gamma} =N
 \cr
                     \hbox{ or} 
 \cr
\sqrt{-2\beta}-\sqrt{2\gamma}=N.
}
\right.
                                                                  \cr
&
+\Sigma_{n= |N|+1}^\infty P_n(\ln x;a)x^{-n} & &
\left.
\matrix{
\hbox{ }
\cr
\beta=-\gamma\hbox{ if } N=0.
}
\right.
\cr
&
 &
   &
     \cr
&
        b_0= {\sqrt{-2\beta}\over N}   \equiv {\sqrt{-\beta}\over \sqrt{-\beta}+\sqrt{\gamma}} 
\hbox{ or }
     {\sqrt{-\beta}\over \sqrt{-\beta}-\sqrt{\gamma}} ,~~N\neq 0.
   &
a     
 &
        \hbox{\bf and}\cr
&&&\cr
\hbox{(\ref{log1zero}.$\infty$)}
   &
\matrix{       
 y(x)= a\pm\sqrt{-2\beta}\ln x+\sum_{n=1}^\infty P_n(\ln x;a)x^{-n},~~N= 0.
\cr
\cr
\hbox{Basic solution when $N=0$.}
}
    &
    &    
                                                \matrix{
\sqrt{-2\beta}\neq
\cr
 \left\{\matrix{ 0,-1,-2,...,N, \hbox{ if } N \leq 0
\cr
                                         0,~1,~2,...,N, \hbox{ if } N\geq 0,
                                         }\right.
    \cr\cr
                                      \sqrt{2\alpha}\pm\sqrt{1-2\delta}\not\in{\cal N}_N.
\cr
\cr
p_{01}=2\cos\pi N=\pm 2.
\cr
\cr
<M_0,M_1>\hbox{ reducible.}
}                                
\cr
&&&\cr
\hline
&&&\cr
\hbox{(\ref{logsquare}.$\infty$)}
&
\matrix{
y(x)={\beta+\gamma\over 2}(a-\ln x)^2+{\beta\over \beta+\gamma}+\sum_{n=1}^\infty P_n(\ln x;a)x^{-n} 
\cr
\cr
\hbox{Basic solution}
}
& a & 
\matrix{
\beta\neq -\gamma.
\cr
p_{01}=2.
\cr
\hbox{no reduc. subgroups.}
}
\cr
&&&\cr
\hline
&&&
\cr
\hbox{(\ref{LOG12}.$\infty$)} 
     &
        \matrix{
  y(x)= b_\infty x+\sum_{n=0}^{|N|-2}b_n x^{-n}+\Bigl(a+b_N\ln x\Bigr)x^{1-|N|}+
\cr 
\cr
+\sum_{n= |N|}^\infty P_n(\ln x;a)x^{-n}
\cr
\cr
\cr 
b_\infty={N\over \sqrt{2\alpha}}\equiv {\sqrt{2\alpha}+\sqrt{1-2\delta}\over \sqrt{2\alpha}} \hbox{ or } {\sqrt{2\alpha}-\sqrt{1-2\delta}\over \sqrt{2\alpha}} 
}
& 
a
&\left.
                                \matrix{
                        \sqrt{2\alpha}+\sqrt{1-2\delta} =N\neq 0
\cr
\hbox{ or} \cr
                           \sqrt{2\alpha}-\sqrt{1-2\delta} =N\neq 0,
\cr                              
\cr
                                {\bf and}
\cr\cr
                            \sqrt{2\alpha}\neq
\cr
\left\{\matrix{0,-1,...,N,\hbox{ if }N\leq -1
\cr
0,~1,...,N, \hbox{ if } N\geq 1,
}\right.
\cr
\cr
                 \sqrt{-2\beta}\pm \sqrt{2\gamma}\not \in {\cal N}_N.
\cr
\cr
p_{01}=-2\cos\pi N=\pm 2.
\cr
\cr
<M_0M_1,M_x>\hbox{ reducible.}
  }
    \right.
\cr
&&&\cr
\hline
\end{array}
$$


$$
\begin{array}{||c c|c|c||} 
\hline\hline
&&&\cr
&
\hbox{{\bf Inverse  logarithmic behaviours}}
& \matrix{\hbox{Free}
\cr
\hbox{Par.}} & \hbox{ Other Conditions }
\cr
\hline
&&&
\cr
\hbox{(\ref{LOG45}.$\infty$)}
&
\left.
\matrix{
y(x)=x \left\{a\pm \sqrt{2\alpha}\ln x+\sum_{n=1}^\infty P_n(\ln x;a)x^{-n}\right\}^{-1}
\cr
\cr
= \pm{x\over \sqrt{2\alpha}\ln x}\left[1 \pm  {a\over \sqrt{2\alpha} \ln x}+O\left({1\over \ln^2 x}\right)\right]
}
\right.
&
a
&
\matrix{
2\alpha=1-2\delta\neq 0.
\cr
p_{01}=-2.
\cr
<M_0M_1,M_x>\hbox{ reducible.}
}
\cr
&&&\cr
\hline
&&&\cr
\hbox{(\ref{LOG3}.$\infty$)} 
& 
\left.
\matrix{
y(x)= x\left\{{2\alpha\over 2\alpha +2\delta-1}+ {1-2\delta-2\alpha\over 4} (a-\ln x)^2 +\sum_{n=1}^\infty P_{n+1}(\ln x;a) x^{-n}\right\}^{-1}
\cr\cr
=
{4~x\over (1-2\delta-2\alpha)\ln^2 x}\left[
1+{2a\over \ln x} +O\left({1\over \ln^2 x}\right)
\right]
}
\right.
&
a
&
\matrix{
2\alpha\neq 1-2\delta.
\cr
p_{01}=-2
\cr
\hbox{no reduc. subgroups.}
}
\cr
&&&\cr
\hline
\end{array}
$$


\section{Table of Parametric Connection  Formulae of type (\ref{POA}), (\ref{POA1}) at $x=0$}
\label{dokidoki}

\vskip 0.2 cm 
The  integration constants ($\sigma$, $a$, $\nu$, $\phi$, etc) in Tables 1, 2 and 3, are denoted  $c_1, c_2$ in  the right hanside of (\ref{malamente1}), which symbolises a critical behaviour at one of the critical points.  $\sigma$, $a$, $\nu$, $\phi$, are expressed  in terms of the monodromy data by parametric formulae of type (\ref{POA}). 
Conversely, 
the monodromy data  
 as functions of the integration constats are given  by formulae of type (\ref{POA1}).  This section provides formulae of type (\ref{POA}) and (\ref{POA1})  at $x=0$.

Parametric connection formulae are computed in \cite{Jimbo} and \cite{Boalch} (cases (\ref{fullEXP}) and (\ref{atopy})),  and in  \cite{D2}, \cite{D1} and \cite{guz2010} (all  the other cases). 

\subsection{Basic Solutions}

\vskip 0.3 cm 
 It would be a formidable task to give the parametric formulae for all the expansions in the table. Fortunately, we need to consider only a set of {\it basic} ones and Okamoto's bi-rational transformations introduced in  \cite{Okamoto} and reviewed in Appendix A.

\vskip 0.2 cm 
\noindent 
{\bf Definition:} A set of {\it basic solutions } is a set of solutions  that  generate all the other solutions of the equation via a bi rational transformation.  A set of  {\it basic  expansions}  is a set of  expansions that can generate, 
via a bi rational transformation, all the  expansions which formally satisfy   the equation . 

\vskip 0.2 cm 
\bpr
\label{spero}
{\it All the solutions and expansions  in the Table 1 for $x\to 0$  are obtained form  the basic solutions  (\ref{fullEXP}), (\ref{atopy}), and the basic expansions  (\ref{log1zero}) and (\ref{logsquare}), via a bi rational transformation.  Tables 2 and 3, for $x\to 1$ and $x\to \infty$  are obtained from the Table 1  via the bi rational transformations (\ref{onara}) and (\ref{onara1}) respectively. }
\epr
{\it Proof:} See the Appendix A.

\vskip 0.2 cm 
The parametric connection formulae  for  a critical behaviour generated by  a basic one are obtained  by the action of the bi rational transformation on the monodromy data of the basic one (\cite{Jimbo}, \cite{DM}, \cite{D1}, \cite{DM11}).  An example is   (\ref{ricordo3}), associated to the bi rational transformation (\ref{sym2}). Another example is (\ref{ricordo1}) and  (\ref{ricordo2}) to be introduced later. Thus, it is enough to give the parametric connection formulae (\ref{POA}), (\ref{POA1})  for the basic behaviours (\ref{fullEXP}), (\ref{atopy}), (\ref{log1zero}) and (\ref{logsquare}). We will give the formulae also for some other cases, marked with a star $\star$ in Table 1: Namely:

\vskip 0.2 cm
\noindent
--  The basic critical behaviours   (\ref{fullEXP}), (\ref{atopy}),  (\ref{log1zero}) and (\ref{logsquare}).

\vskip 0.2 cm
\noindent
-- The solutions (\ref{lantern1}) and (\ref{TAU}), important for their behaviour with poles  and their applications in quantum cohomology (see for example \cite{Dpoli}). 

\vskip 0.2 cm
\noindent
-- Solutions (\ref{UUU}) and (\ref{div}).

\vskip 0.2 cm
\noindent
-- The Taylor expansion (\ref{TTLO2}), as a limit of  (\ref{log1zero}). 

\vskip 0.2 cm
\noindent
--  The Taylor expansion (\ref{TTLO4}) and one of the two inequivalent expansions (\ref{TTLO3}) when $N=1$, namely: 
$$
y(x)={1\over \sqrt{2\alpha}} +ax+\sum_{n=2}^\infty b_n(a)x^n,~~~\alpha\neq 0,~~\sqrt{2\alpha}\pm\sqrt{2\gamma}=1,~~2\beta=2\delta-1.~~~~~~~\hbox{(\ref{T2coe}) }
$$
(Note: $y(x)=-{1\over \sqrt{2\alpha}} +ax+\sum_{n=2}^\infty b_n(a)x^n,$ occurs for $\sqrt{2\alpha}\pm\sqrt{2\gamma}=-1$; this is equivalent to the choice of the other sign of $\sqrt{2\alpha}$ in (\ref{T2coe})).  
Together with (\ref{T1coe}) (which has no free parameter), they are {\it basic Taylor solutions}, which generate all the other Taylor expansions  via bi rational transformations. 

\vskip 0.5 cm 
\noindent
Let:

 $p_{ij}=$Tr$(M_iM_j)$, $i,j=0,x,1$;~~~ $p_\mu =2\cos\pi\theta_\mu$, $\mu=0,x,1,\infty$;
\vskip 0.2 cm 
 $s(z):=\sin\left({\pi\over 2}z\right)$, ~~~$c(z):=\cos\left({\pi\over 2}z\right)$.

\subsection{Table}

\vskip 0.3 cm 
\noindent
{\large \bf Case   (\ref{fullEXP}) and (\ref{atopy}) .} The two integration constants are $\sigma$ and $a$. Their parametrisation is computed in \cite{Jimbo} and \cite{Boalch}) (note that in \cite{Jimbo}, the last sign in formula (1.8) at the bottom of page 1141 is $\pm \sigma$ instead of $\mp \sigma$):
\be
\label{aF}
\sq{
\sigma= {1\over \pi} \hbox{arcos}\left({p_{0x}\over 2}\right)
},~~~~~
 \sq{
a= {(\theta_x-\theta_0-\sigma)(\theta_x+\theta_0-\sigma)(\theta_\infty+\theta_1-\sigma)\over4 \sigma^2(\theta_\infty+\theta_1+\sigma)}~{1\over  {\bf F}}~ {U\over V} 
}
\ee
where 
$$
{\bf F}:=  {
\Gamma(1+\sigma)^2\Gamma \left({1\over 2}(\theta_0+\theta_x-\sigma)+1
\right) \Gamma\left(  
{1\over 2} ( \theta_x-\theta_0-\sigma)+1
\right)
\over 
\Gamma(1-\sigma)^2 \Gamma \left({1\over2}(\theta_0+\theta_x+\sigma)+1
\right) \Gamma\left(   
{1\over 2} ( \theta_x-\theta_0+\sigma)+1
\right)
}~\times
$$
$$
\times  {
\Gamma \left({1\over 2}(\theta_\infty+\theta_1-\sigma)+1 \right)
 \Gamma\left(  
{1\over 2} ( \theta_1-\theta_\infty-\sigma)+1
\right)
\over
\Gamma \left({1\over 2}(\theta_\infty+\theta_1+\sigma)+1 \right)
 \Gamma\left( 
{1\over 2} ( \theta_1-\theta_\infty+\sigma)+1
\right)
},
$$

$$
   U:={i \over 2}\left(p_{01}+p_{x1}e^{i\pi\sigma}\right)\sin\pi\sigma+
 \cos(\pi\theta_x) \cos(\pi \theta_1) + \cos(\pi\theta_\infty) \cos(\pi \theta_0)~+$$
$$-\left[
 \cos(\pi\theta_x) \cos(\pi \theta_\infty) + \cos(\pi\theta_0) \cos(\pi \theta_1) \right]e^{i\pi\sigma},
$$
$$ 
V:=4 ~s(\theta_0+\theta_x-\sigma) s(\theta_0 - \theta_x+\sigma)~ s (\theta_\infty+\theta_1-\sigma)
 s (\theta_\infty-\theta_1+\sigma).
$$
For both (\ref{fullEXP}) and (\ref{atopy}), we have 
$$
 \hbox{Tr}(M_xM_0)\not\in(-\infty,-2]\cup\{2\}
$$
In particular 
$$
 \hbox{Tr}(M_xM_0)>2 \hbox{ when } \sigma\in i{\bf R} 
$$

\vskip 0.3 cm 
Conversely, the monodromy data in terms of the integration constants are (\cite{Jimbo}, \cite{Boalch} and  \cite{guz2010}):
\be
\label{pF}
\sq{
p_{0x}=2\cos \pi \sigma,~~~p_{x1}={\bf G}_1~a^{-1} ~+ {\bf G}_2+{\bf G}_3 ~a,~~~p_{01}={\bf G}_4~ a^{-1}~+{\bf G}_5+
{\bf G}_6~ a
}
\ee
where (see \cite{guz2010}): 

$$
{\bf G}_2= {2(\Omega~\cos\pi
\theta_x\cos\pi\theta_1~-~\Xi~\sin\pi\theta_x \sin\pi\theta_1)
\over
\sin^2(\pi\sigma)\sin\pi\theta_x\sin\pi\theta_1} ;
$$

$${\bf G}_5= 2\Bigl(\cos\pi\theta_1\cos\pi\theta_0+{\Xi_1\over
  \Omega_1}~\sin\pi\theta_1 \sin\pi\theta_0\Bigr);
$$

$$
{\bf G}_1={\sin\pi\theta_x \sin\pi\theta_1 \over \Omega} ~V_1~{1\over {\cal F}},
~~~~~
{\bf G}_3={\sin\pi\theta_x \sin\pi\theta_1 \over \Omega} ~V~{\cal F}
$$
$$
{\bf G}_4= - e^{i\pi\sigma}{\sin\pi\theta_0\over \sin\pi\theta_x}{\Omega\over \Omega_1} {\bf G}_1,
~~~~~
{\bf G}_6= - e^{-i\pi\sigma}{\sin\pi\theta_0\over \sin\pi\theta_x}{\Omega\over \Omega_1} {\bf G}_3.
$$
 and :

$$
{\cal F}:=  {4\sigma^2~(\theta_\infty+\theta_1+\sigma)\over (\theta_0-\theta_x+\sigma)(\theta_0+\theta_x-\sigma)(\theta_\infty+\theta_1-\sigma)}~{\bf F}
$$
\vskip 0.2 cm 
$$
V_1:= V(\sigma\mapsto-\sigma)
$$

$$
\Xi:=\Bigl(s(\theta_0+\theta_x+\sigma)s(\theta_0-\theta_x-\sigma)
+s(\theta_0-\theta_x+\sigma)s(\theta_0+\theta_x-\sigma)\Bigr)\times
$$
$$
\times
\Bigl(s(\theta_1+\theta_\infty+\sigma)s(\theta_1-\theta_\infty+\sigma)
+s(\theta_1+\theta_\infty-\sigma)s(\theta_1-\theta_\infty-\sigma)\Bigr)
$$

$$
\Xi_1:=\Bigl(s(\theta_0+\theta_x+\sigma)s(\theta_0-\theta_x+\sigma)
+s(\theta_0+\theta_x-\sigma)s(\theta_0-\theta_x-\sigma)\Bigr)\times
$$
$$
\times
\Bigl(s(\theta_1+\theta_\infty+\sigma)s(\theta_1-\theta_\infty+\sigma)
+s(\theta_1+\theta_\infty-\sigma)s(\theta_1-\theta_\infty-\sigma)\Bigr)
$$

$$
\Omega:=\Bigl(-s(\theta_0+\theta_x+\sigma)s(\theta_0-\theta_x-\sigma)
+s(\theta_0-\theta_x+\sigma)s(\theta_0+\theta_x-\sigma)\Bigr)\times
$$
$$
\times
\Bigl(s(\theta_1+\theta_\infty+\sigma)s(\theta_1-\theta_\infty+\sigma)
-s(\theta_1+\theta_\infty-\sigma)s(\theta_1-\theta_\infty-\sigma)\Bigr)
$$

$$
\Omega_1:=\Bigl(s(\theta_0+\theta_x+\sigma)s(\theta_0-\theta_x+\sigma)
-s(\theta_0+\theta_x-\sigma)s(\theta_0-\theta_x-\sigma)\Bigr)\times
$$
$$
\times
\Bigl(s(\theta_1+\theta_\infty+\sigma)s(\theta_1-\theta_\infty+\sigma)
-s(\theta_1+\theta_\infty-\sigma)s(\theta_1-\theta_\infty-\sigma)\Bigr)
$$

Let 
$$
Z(\sigma^2):={16\sigma^4\over (\sigma^2-(\theta_0+\theta_x)^2)(\sigma^2-(\theta_0-\theta_x)^2)},$$
and  observe that
 $$
\Xi(\sigma)=\Xi(-\sigma),~~~ \Xi_1(\sigma)=\Xi_1(-\sigma),~~~\Omega(\sigma)=\Omega(-\sigma),~~~\Omega_1(\sigma)=\Omega_1(-\sigma),
$$
$$F(-\sigma)={1\over F(\sigma)},~~~{\cal F}(-\sigma)={1\over Z(\sigma^2) {\cal F}(\sigma)},
$$
This implies that
$$
G_2(\sigma)=G_2(-\sigma),~~~G_5(\sigma)=G_5(-\sigma),~~~
G_1(-\sigma)=G_3(\sigma)Z(\sigma^2),~~~G_4(-\sigma)=G_6(\sigma)Z(\sigma^2).
$$
and the above, substituted in (\ref{pF}), implies that 
$$
 a(-\sigma)={Z(\sigma^2)\over a(\sigma)}.
$$

\vskip 0.2 cm 
\noindent
{\bf Remark 2:} In  case 
(\ref{atopy}), $a$ is  the limit of (\ref{aF}) for $\sigma$ tending to one of $\pm(\theta_0+\theta_x)$ or $\pm(\theta_0-\theta_x)$. The limit exist finite and non zero. 
\vskip 0.2 cm 

\noindent
{\bf Remark 3:} When $p_{0x}= 2\cos \pi (\theta_\infty-\theta_1)$, the Jimbo-Fricke cubic (\ref{FR}) has factorisation (\ref{era}),  
and $\sigma=\pm(\theta_\infty-\theta_1)+2n_{\pm}$, for  suitable integers $n_+$ and $n_-$  such that  $-1<\Re\sigma<1$ (or $\Re\sigma>-1$ in case of (\ref{atopy})).  

The value of $a$ for $\sigma =(\theta_\infty-\theta_1)+2n_{+}$ is obtained by the limit of (\ref{aF}) when $\sigma \to (\theta_\infty-\theta_1)+2n_{+}$. The limit exists, finite and not zero,   when the monodromy data $p_{x1}$ and $p_{01}$ are such that one of the two factors in (\ref{era}) vanishes.  On the other hand, the limit of $a$ is zero of infinite for the monodromy data $p_{x1}$ and $p_{01}$  such that the other factor of (\ref{era}) vanishes. This is because 
 there are compensations between the vanishing of factors  in the numerator or/and the denominator of  (\ref{aF})  and the behaviour of  $U$, which  is equal to the first factor of (\ref{era}), up to multiplication by $2/[i\sin\pi ((\theta_\infty-\theta_1)]$.  For the monodromy data such that the limit is finite$\neq 0$, the branch $y(x)$ is (\ref{fullEXP}) or  (\ref{atopy}), with $\sigma =(\theta_\infty-\theta_1)+2n_{+}$.   For the other choice of monodromy data, no branch $y(x)$ of type (\ref{fullEXP}),(\ref{atopy})  exists, and the corresponding critical behaviour  is of type (\ref{UUU}) or (\ref{div}). 
The same considerations apply when $\sigma \to -(\theta_\infty-\theta_1)+2n_{-}$, when $U$ is proportional to the second factor of (\ref{era}). 

When $p_{0x}= 2\cos \pi (\theta_\infty+\theta_1)$, the Jimbo-Fricke cubic has factorisation (\ref{era1}) 
and  $\sigma=\pm(\theta_\infty+\theta_1)+2n_{\pm}$  for suitable integers $n_\pm$ such that  $-1<\Re\sigma<1$ (or $\Re\sigma>-1$ in case of (\ref{atopy})). The same considerations as above apply to the limit of $a$ in (\ref{aF}) for $\sigma\to\pm(\theta_\infty+\theta_1)+2n$. There are two choices for the monodoromy data, corresponding to the vanishing of one or the other factor in (\ref{era1}).  For one choice, $a$ exist finite$\neq0$ and the branch is (\ref{fullEXP}) or (\ref{atopy}), for the other choice the limit is zero or infinite, and the critical behaviour is  (\ref{UUU}) or (\ref{div}).

\vskip 0.3 cm
\noindent
{\large \bf Case  (\ref{UUU}) and (\ref{div}).}
The parameter $a$ of (\ref{fullEXP}) and (\ref{atopy}) is expressed in formula (\ref{aF}), in terms of the monodromy data. The parameter $a$ in (\ref{UUU}) and (\ref{div}) is given by the same formula, provided one makes the substitutions 
$(p_{0x},p_{01},p_{x1})\mapsto( -p_{0x},-p_{01},p_{x1})$,  $(\theta_0,\theta_\infty)\mapsto (\theta_\infty-1,\theta_0+1)$, and $\sigma\mapsto 1+\rho$  for (\ref{UUU}),  $\sigma\mapsto \omega$ for (\ref{div}).

\vskip 0.3 cm
\noindent
{\large \bf Case  (\ref{lantern1}) and (\ref{TAU}).} The constants $\nu$ and  $\phi$ of (\ref{lantern1})  in terms of the monodromy data are computed in \cite{guz2010}. 
$$
\sq{2\cosh(2\pi\nu)=-p_{0x},~~~ \phi= i\ln\left(-{r\over \nu A}\right)}
$$
$A$ is the coefficient multiplying the sine in (\ref{lantern1}):
$$
A={1\over 8\nu^2}\sqrt{[4\nu^2+(\theta_\infty-1-\theta_1)^2][4\nu^2+(\theta_\infty-1+\theta_1)^2]}
$$
and 
$$
\sq{
r
=
 {(\theta_\infty-1-\theta_1-2i\nu)(\theta_\infty-1+\theta_1+2i\nu)
(\theta_0+\theta_x+1+2i\nu)\over
 8i\nu(
 2i\nu-\theta_0-\theta_x-1)} ~ {1\over  {\bf F}^*}~ {{\cal U}\over {\cal V}}
}
$$
where
$$
{\bf F}^*:= 
 {
\Gamma(1-2i\nu)^2
\Gamma \left({1\over 2}(\theta_\infty-1+\theta_1+2i\nu)+1\right) 
\Gamma\left(  {1\over 2} ( \theta_1+1-\theta_\infty+2i\nu)+1\right)
\over 
\Gamma(1+2i\nu)^2 
\Gamma \left({1\over2}(\theta_\infty-1+\theta_1-2i\nu)+1\right) 
\Gamma\left( {1\over 2} ( \theta_1+1-\theta_\infty-2i\nu)+1\right)
}~\times
$$
$$
\times  {
\Gamma \left({1\over 2}(\theta_0+\theta_x+1+2i\nu)+1 \right)
\Gamma\left({1\over 2} ( \theta_x-\theta_0-1+2i\nu)\right)
\over
\Gamma \left({1\over 2}(\theta_0+\theta_x+1-2i\nu)+1 \right)
\Gamma\left( {1\over 2} ( \theta_x-\theta_0-1-2i\nu)\right)
},
$$

$$
  {\cal U}:=e^{2\pi\nu}\left[
{1 \over 2}\sinh(2\pi\nu)p_{1x}+ \cos(\pi\theta_x) \cos(\pi \theta_\infty) + \cos(\pi\theta_0) \cos(\pi \theta_1) \right]~+$$
$$
-{1\over 2}
 \sinh(2\pi\nu)p_{01}+ \cos(\pi\theta_x) \cos(\pi \theta_1) + \cos(\pi\theta_\infty) \cos(\pi \theta_0),
$$

$$ 
{\cal V}:=4 ~c(\theta_0+\theta_x+2i\nu)~ c (\theta_0 - \theta_x-2i\nu)~
 c(\theta_\infty+\theta_1+2i\nu)
 ~c (\theta_\infty-\theta_1-2i\nu).
$$

\vskip 0.2 cm 
As for the solution (\ref{TAU}), we have 
$$  
  a=-{r\over 2i\nu}
$$
where $r$ is the above. Note that in this case $2i\nu$ is equal to one of $\pm(\theta_\infty-1+\theta_1)$ or  $\pm(\theta_\infty-1-
\theta_1)$, and the limit of $r$ above exists for these values of $2i\nu$.

\vskip 0.2 cm 
For (\ref{lantern1}) and (\ref{TAU}), we have 
$$
 \hbox{Tr}(M_xM_0)<-2
$$

\vskip 0.2 cm 
Conversely, the monodromy data in terms of the integration constants are  \cite{guz2010}:
$$
\sq{
p_{0x}=-2\cosh (2\pi \nu),~~~p_{x1}={\bf G}_1^{*}~r^{-1} ~+ {\bf G}_2^{*}+{\bf G}_3^{*} ~r,
~~~p_{01}=-\Bigl({\bf G}_4^{*}~ r^{-1}~+{\bf G}_5^{*}+{\bf G}_6^{*}~ r\Bigr)
}
$$
where: 

$$
{\bf G}_2^{*}= -{2(\Omega^{*}~\cos\pi
\theta_x\cos\pi\theta_1~-~\Xi^{*}~\sin\pi\theta_x \sin\pi\theta_1)
\over
\sinh^2(2\pi\nu)\sin\pi\theta_x\sin\pi\theta_1} ;
$$

$${\bf G}_5^{*}=- 2\Bigl(\cos\pi\theta_x\cos\pi\theta_\infty+{\Xi_1^{*}\over
  \Omega_1^{*}}~\sin\pi\theta_x \sin\pi\theta_\infty\Bigr);
$$

$$
{\bf G}_1^{*}=-{\sin\pi\theta_x \sin\pi\theta_1 \over \Omega^{*}} ~{\cal V}_1~{1\over {\cal F}^{*}},
~~~~~
{\bf G}_3^{*}=-{\sin\pi\theta_x \sin\pi\theta_1 \over \Omega^{*}} ~{\cal V}~{\cal F}^{*}
$$
$$
{\bf G}_4^{*}=  e^{2\pi\nu}{\sin\pi\theta_\infty\over \sin\pi\theta_1}{\Omega^{*}\over \Omega_1^{*}} {\bf G}_1^{*},
~~~~~
{\bf G}_6^{*}=  e^{-2\pi\nu}{\sin\pi\theta_\infty\over \sin\pi\theta_1}{\Omega^{*}\over \Omega_1^{*}} {\bf G}_3^{*}.
$$
 and :

$$
{\cal F}:=  {8i\nu~(2i\nu-\theta_0-\theta_x-1)\over (\theta_\infty-\theta_1-1-2i\nu)(\theta_\infty+\theta_1-1+2i\nu)(\theta_0+\theta_x+1+2i\nu)}~{\bf F}^{*}
$$
\vskip 0.2 cm 
$$
{\cal V}_1:= {\cal V}(\nu\mapsto-\nu)
$$

$$
\Xi^{*}:=-\Bigl(c(\theta_\infty+\theta_1-2i\nu)c(\theta_\infty-\theta_1+2i\nu)
+c(\theta_\infty-\theta_1-2i\nu)c(\theta_\infty+\theta_1+2i\nu)\Bigr)\times
$$
$$
\times
\Bigl(c(\theta_x+\theta_0-2i\nu)c(\theta_x-\theta_0-2i\nu)
+c(\theta_x+\theta_0+2i\nu)c(\theta_x-\theta_0+2i\nu)\Bigr)
$$

$$
\Xi_1^{*}:=-\Bigl(c(\theta_\infty+\theta_1-2i\nu)c(\theta_\infty-\theta_1-2i\nu)
+c(\theta_\infty+\theta_1+2i\nu)c(\theta_\infty-\theta_1+2i\nu)\Bigr)\times
$$
$$
\times
\Bigl(c(\theta_x+\theta_0-2i\nu)c(\theta_x-\theta_0-2i\nu)
+c(\theta_x+\theta_0+2i\nu)c(\theta_x-\theta_0+2i\nu)\Bigr)
$$

$$
\Omega^{*}:=-\Bigl(-c(\theta_\infty+\theta_1-2i\nu)c(\theta_\infty-\theta_1+2i\nu)
+c(\theta_\infty-\theta_1-2i\nu)c(\theta_\infty+\theta_1+2i\nu)\Bigr)\times
$$
$$
\times
\Bigl(c(\theta_x+\theta_0-2i\nu)c(\theta_x-\theta_0-2i\nu)
-c(\theta_x+\theta_0+2i\nu)c(\theta_x-\theta_0+2i\nu)\Bigr)
$$

$$
\Omega_1^{*}:=-\Bigl(c(\theta_\infty+\theta_1-2i\nu)c(\theta_\infty-\theta_1-2i\nu)
-c(\theta_\infty+\theta_1+2i\nu)c(\theta_\infty-\theta_1+2i\nu)\Bigr)\times
$$
$$
\times
\Bigl(c(\theta_x+\theta_0-2i\nu)c(\theta_x-\theta_0-2i\nu)
-c(\theta_x+\theta_0+2i\nu)c(\theta_x-\theta_0+2i\nu)\Bigr)
$$

\vskip 0.3 cm \noindent
{\large \bf Case   (\ref{T2coe}) $\subset$ \{(\ref{TTLO3}) with $N=1$\}:}  For the basic Taylor expansion  (\ref{T2coe}), the parametrisation is computed in  \cite{D2}: \footnote{
in the form 
$$
y(x)= {1\over 1-\theta_\infty} +a~x ~+O(x^2),~~~\theta_\infty\pm\theta_1=0,~~~\theta_x\pm \theta_0=0
$$
$$
a={\theta_\infty(2s+\theta_x+1)\over 
2(\theta_\infty-1)},~~~
s={\theta_x\bigl[2\cos(\pi(\theta_\infty+\theta_x))-p_{01}\bigr]\over
 2\bigl[\cos(\pi(\theta_\infty-\theta_x))-
\cos(\pi(\theta_\infty+\theta_x))\bigr]} . 
$$
$$
p_{01}=2\cos(\pi(\theta_\infty+\theta_x))-{4s~\sin(\pi\theta_\infty)\sin(\pi\theta_x)\over \theta_x},~~~p_{01}+p_{1x}=4\cos(\pi\theta_\infty)\cos(\pi\theta_x)
$$
Then, choose $\sqrt{2\alpha}=1-\theta_\infty$}
$$
\sq{a={\sqrt{2\alpha}-1\over 2\sqrt{2\alpha}}\left[
1-{p_{01}+2\cos(\pi\sqrt{2\alpha})~\cos(\pi\sqrt{1-2\delta})\over 2\sin(\pi\sqrt{2\alpha})}{\sqrt{1-2\delta}\over \sin(\pi\sqrt{1-2\delta})}
\right]
}
$$
and conversely
$$
\sq{
\matrix{p_{01}={2\sin(\pi\sqrt{2\alpha})~\sin(\pi\sqrt{1-2\delta})\over (1-\sqrt{2\alpha})\sqrt{1-2\delta}}\left(1+\sqrt{2\alpha}~(2a-1)\right)~-2\cos\left(\pi\sqrt{2\alpha}\right)\cos\left(\pi\sqrt{1-2\delta}\right)
\cr
\cr
p_{0x}=2,~~~~~p_{x1}+p_{01}=-4\cos\left(\pi\sqrt{2\alpha}\right)\cos\left(\pi\sqrt{1-2\delta}\right)
}
}
$$
Parametric formulae for
$$
  y(x)=-{1\over \sqrt{2\alpha}}+ax+O(x^2). 
$$
are obtained from the above, with $\sqrt{2\alpha}\mapsto-\sqrt{2\alpha}$. The monodromy is 
$$
\hbox{Tr}(M_xM_0)=2,~~~<M_0\cdot M_x,M_1> \hbox{ is reducible}
$$
In particular $M_xM_0=I$.

\vskip 0.3 cm \noindent
{\large \bf Case   (\ref{TTLO4}):} 
  For the basic Taylor expansion  (\ref{TTLO4}), with $\alpha=\gamma=0$, the parametrisation is computed in  \cite{D2} (there is a mistake in \cite{D2}, and the formulae below have been computed from (\ref{mac}) below,  via (\ref{sym2}) and (\ref{ricordo3}).). 
$$
 (\ref{TTLO4})~~~~~~~
y(x)~ =~a ~+~{1-a\over 2} (\delta-\beta)~x~+~O(x^2) 
$$
$$
\sq{a= {4\cos\left({\pi\over 2}\left[\sqrt{-2\beta}+\sqrt{1-2\delta}  \right] \right)
     \cos\left({\pi\over 2}\left[\sqrt{-2\beta}-\sqrt{1-2\delta}  \right] \right)    
\over 
   2\cos\Bigl(\pi\sqrt{1-2\delta}\Bigr)+p_{01}
}
}
$$
%
Conversely:
 $$
\sqs{
\matrix{
p_{0x}=-2,
\cr 
\cr 
p_{01}=-2\cos(\pi\sqrt{1-2\delta})+4a^{-1}\cos\left({\pi\over 2}\left[\sqrt{-2\beta}+\sqrt{1-2\delta}  \right] \right)
     \cos\left({\pi\over 2}\left[\sqrt{-2\beta}-\sqrt{1-2\delta}  \right] \right)
\cr 
\cr 
p_{x1}=2\cos(\pi\sqrt{1-2\delta})-4(1-a^{-1})\cos\left({\pi\over 2}\left[\sqrt{-2\beta}+\sqrt{1-2\delta}  \right] \right)
     \cos\left({\pi\over 2}\left[\sqrt{-2\beta}-\sqrt{1-2\delta}  \right] \right)
}
}
$$ 
In  case  (\ref{TTLO4}) we have: 
$$
\hbox{Tr}(M_xM_0)=-2,~~~<M_0\cdot M_x,M_1> \hbox{ is reducible}
$$

\vskip 0.3 cm \noindent
{\large \bf Case   (\ref{log1zero}):} 
 The free parameter of (\ref{log1zero}) is computed  in  \cite{D1}, for $\theta_\infty\pm \theta_1\not\in2{\bf Z}$,  $\theta_x\pm\theta_0=0$. The parametrisation is: 
$$
\sq{
{1\over \sqrt{-2\beta}}~a=
\left\{
{\pi\left[
2\cos\left(\pi\left[\sqrt{-2\beta}+\sqrt{2\gamma}\right]\right)-p_{01}
\right]
\over 4 \cos\left({\pi\over 2}\left[\sqrt{2\alpha}+\sqrt{2\gamma}\right]\right)\cos\left({\pi\over 2}\left[\sqrt{2\alpha}-\sqrt{2\gamma}\right]\right)
\sin(\pi\sqrt{-2\beta})
}+\omega
\right\}e^{i{\pi\over 2}\sqrt{-2\beta}}+\omega_1
}
$$
where
$$
\omega= \Psi_E\left({\sqrt{2\alpha}\over 2} +{\sqrt{2\gamma}\over 2} +{1\over 2} \right)-\Psi_E\left({\sqrt{2\gamma}\over 2} -{\sqrt{2\alpha}\over 2} +{1\over 2} \right)+2\gamma_E,~~~\omega_1:=\gamma_E+\Psi_E(1+\sqrt{-2\beta})+i\pi
$$
 Here $\Psi_E$ is the Euler's psi-function, and $\gamma_E$ is the Euler's gamma constant. Conversely: 
$$
\sqs{
\matrix{p_{0x}=2,~~~~~
p_{01}=2\cos\left(\pi\left[\sqrt{-2\beta}+\sqrt{2\gamma}\right]\right)-4\cos\left({\pi\over 2}\left[\sqrt{2\alpha}
+\sqrt{2\gamma}\right]\right)\cos\left({\pi\over 2}\left[\sqrt{2\alpha}-\sqrt{2\gamma}\right]\right)
\Omega(a)
\cr 
\cr 
p_{x1}= 
2\cos\left(\pi(\sqrt{-2\beta}-\sqrt{2\gamma})\right)
+
4 \cos\left({\pi\over 2}\left[\sqrt{2\alpha}+\sqrt{2\gamma}\right]\right)\cos\left({\pi\over 2}\left[\sqrt{2\alpha}-\sqrt{2\gamma}\right]\right)
\left\{\Omega(a)-e^{-i\pi\sqrt{-2\beta}}
\right\}
}
}
$$
where
$$
\Omega(a):={1\over \pi}\sin(\pi\sqrt{-2\beta})\left[
e^{-{i\over 2}\pi\sqrt{-2\beta}}\left({a\over \sqrt{-2\beta}}-\omega_1\right)-\omega\right]
$$
In this case 
$$
\hbox{Tr}(M_xM_0)=2,~~~<M_0,M_x> \hbox{ is reducible}
$$

\vskip 0.3 cm \noindent
{\large \bf Case (\ref{TTLO2}):} 
When $\theta_0\to 0$, namely $\beta\to 0$, the case (\ref{log1zero})  converges to:
$$
(\ref{TTLO2})~~~~~~~y(x)=ax+a(a-1)\left(
 \gamma-\alpha-{1\over 2}
\right)x^2+O(x^3)
$$ 
\be
\label{mac}
\sq{ a = 
 {2\cos\pi\sqrt{2\gamma}-p_{01}
\over 
4\cos\left({\pi\over 2}\left[\sqrt{2\alpha}+\sqrt{2\gamma}  \right] \right)
     \cos\left({\pi\over 2}\left[\sqrt{2\alpha}-\sqrt{2\gamma}  \right]\right)
}
}
\ee
and 
 $$
\sqs{
p_{0x}=2
\cr
\cr 
p_{01}=2\cos\pi\sqrt{2\gamma}-4a\cos\left({\pi\over 2}\left[\sqrt{2\alpha}+\sqrt{2\gamma}  \right] \right)
     \cos\left({\pi\over 2}\left[\sqrt{2\alpha}-\sqrt{2\gamma}  \right]\right)
\cr 
\cr 
p_{x1}=2\cos\pi\sqrt{2\gamma}+4(a-1)\cos\left({\pi\over 2}\left[\sqrt{2\alpha}+\sqrt{2\gamma}  \right] \right)
     \cos\left({\pi\over 2}\left[\sqrt{2\alpha}-\sqrt{2\gamma}  \right]\right)
}
$$

\vskip 0.3 cm \noindent
{\large \bf Case (\ref{logsquare}):} The free parameter for (\ref{logsquare}) is calculated in  \cite{D1}, when $\theta_0,\theta_x,\theta_1,\theta_\infty\not \in {\bf Z}$, $\theta\pm\theta_x\neq 0$. The parametrisation is: 
$$
\sq{
\matrix{
a={{\bf C}-{\bf c}+p_{01}-p_{x1} \over {\bf b}-{\bf B}}+{2\theta_x\over \theta_0^2-\theta_x^2}+4(\gamma_E-i\pi)+
\cr 
\cr 
+\Psi_E\left({\theta_\infty\over 2}+{\theta_1\over 2}\right)-\Psi_E\left({\theta_1\over 2}-{\theta_\infty\over 2}+1\right)
+\Psi_E\left(-{\theta_x\over 2}-{\theta_0\over 2}\right)
+\Psi_E\left({\theta_x\over 2}-{\theta_0\over 2}+1\right)
}
}
$$
where 
$$
 {\bf b}={4\over \pi} (\sin\pi\theta_1 ~s(\theta_0-\theta_x)s(\theta_0+\theta_x) + \sin\pi\theta_0~s(\theta_\infty+\theta_1)s(\theta_\infty-\theta_1))
$$
$$
 {\bf c}= 2\cos\pi(\theta_0-\theta_1)
$$
$$
{\bf B}= {1\over 2\pi i} \Bigl[
2\cos\pi(\theta_0+\theta_1)+4\cos\pi\theta_x \cos\pi\theta_\infty - 4 e^{i\pi\theta_1}\cos\pi\theta_x - 4 e^{-i\pi\theta_0}\cos\pi\theta_\infty + 3 e^{i\pi(\theta_1-\theta_0)} - e^{i\pi(\theta_0-\theta_1)}
\Bigr]
$$
$$
 {\bf C}= 2 e^{i\pi\theta_1}\cos\pi\theta_x+2e^{-i\pi\theta_0}\cos\pi\theta_\infty -2 e^{i\pi(\theta_1-\theta_0)}
$$
and $s(z):=\sin(\pi z/2)$. Conversely: 
$$
\sqs{
\matrix{
p_{0x}=2 
\cr 
\cr 
p_{01}={\bf A} q(a)^2+({\bf b}-2{\bf A}\omega)q(a)+({\bf c} -{\bf b} \omega +{\bf A}\omega^2)
\cr 
\cr 
p_{x1}={\bf A} q(a)^2+({\bf B}-2{\bf A}\omega)q(a)+({\bf C} -{\bf B} \omega +{\bf A}\omega^2)
}
}
$$
where 
$$
q(a)=a+{2\theta_x\over \theta_x^2-\theta_0^2}-\Psi_E\left(-{\theta_0\over 2}-{\theta_x\over 2}\right)
-\Psi_E\left( {\theta_x\over 2}-{\theta_0\over 2}+1 \right)-2\gamma_E+4\pi i,
$$
$$ 
{\bf A}= {4\over \pi^2} s(\theta_0+\theta_x)s(\theta_0-\theta_x)s(\theta_\infty+\theta_1)s(\theta_\infty-\theta_1),
$$
$$
\omega= \Psi_E\left({\theta_\infty\over 2} +{\theta_1\over 2}  \right)-\Psi_E\left({\theta_1\over 2} -{\theta_\infty\over 2} +1 \right)+2\gamma_E
$$
In this case 
$$
\hbox{Tr}(M_xM_0)=2,~~~<M_0,M_x,M_1>\hbox{ irreducible}
$$


\section{Connection Problem -- Connection Formulae in Closed Form} 
\label{CONNE}
Note: the formulae of this section are proved by means of Appendix A.

\subsection{Parametric Connection Formulae for Table 2 ($x=1$) and and Table 3 ($x=\infty$)} 

{\large \bf -- For Table 1 ($x=0$)},  the  parametric connection formulae of type (\ref{POA}) and (\ref{POA1}) are  the formulae of  Section \ref{dokidoki}.

\vskip 0.2 cm 
Critical behaviours (\ref{malamente1}) in Tables 2 and 3  are labelled with the same numbers than in Table 1. To find the   formulae (\ref{POA}) and (\ref{POA1}) for a critical behaviour of Tables 2 or 3, consider the behaviour  of Table 1  with the same numeration, and take the formulae of Section \ref{dokidoki} for such  behaviour. Then,  do the following procedure:  

\vskip 0.2 cm 
\noindent
 {\large \bf -- For Table 2 ($x=1$).}

\vskip 0.2 cm 
\noindent
 {\large  a) \it Parametric formulae of type (\ref{POA}) at $x=1$.}

\vskip 0.2 cm 
 Take the formulae $ c_i=c_i(\theta_0,\theta_x,\theta_1,\theta_\infty,p_{0x},p_{01},p_{x1})$ of type (\ref{POA}) in Section \ref{dokidoki},  and make  the substitutions 
\be
\label{ricordo1}
\left\{
\matrix{
\theta_0\mapsto \theta_1,~~~\theta_x\mapsto \theta_x,~~~\theta_1 \mapsto\theta_0,~~~\theta_\infty \mapsto \theta_\infty
,~~\hbox{ (or } 
~~
\alpha \mapsto \alpha,~\beta\mapsto -\gamma,~\gamma \mapsto -\beta,~\delta \mapsto \delta)
\cr 
\cr
p_{0x} \mapsto p_{x1} ,~~~
p_{x1} \mapsto p_{0x},~~~p_{01}\mapsto -p_{01}-p_{0x}p_{x1}+p_\infty p_x+p_1p_0.
}
\right.
\ee
where $p_\mu=2\cos\pi\theta_\mu$, $\mu=0,x,1,\infty$. 

\vskip 0.2 cm 
\noindent
 {\large b) \it Formulae of type (\ref{POA1}): Monodromy data in terms of the integration constants $c_1$, $c_2$ at $x=1$.}

\vskip 0.2 cm 
Take the formulae of type (\ref{POA1}) in  Section \ref{dokidoki}, but denote them with a prime  
\be
\label{biop}
p_{ij}^{\prime}=p_{ij}^{\prime}(c_1,c_2,\theta_0^{\prime},\theta_x^{\prime},\theta_1^{\prime},\theta_\infty^{\prime})
\ee
where $c_1$ and $c_2$ are now the constants in the behaviors of  Table 2 at $x=1$. 
Then, the monodromy data corresponding to the behaviours of Table 2 are 
$$ 
\left\{
\matrix{
p_{0x}&=&p_{x1}^{\prime}|_{\theta_0^{\prime}=\theta_1,\theta_x^{\prime}=\theta_x,\theta_1^{\prime}=\theta_0,\theta_\infty^{\prime}=\theta_\infty}~~~~~~~~~~~~~~~~~~~~~~~~~~~~~~~
\cr
\cr
p_{x1}&=&p_{0x}^{\prime}|_{\theta_0^{\prime}=\theta_1,\theta_x^{\prime}=\theta_x,\theta_1^{\prime}=\theta_0,\theta_\infty^{\prime}=\theta_\infty}~~~~~~~~~~~~~~~~~~~~~~~~~~~~~~~
\cr
\cr
p_{01}&=&-\bigl(p_{01}^{\prime}+p_{0x}^{\prime} p_{x1}^{\prime}\bigr)|_{\theta_0^{\prime}=\theta_1,\theta_x^{\prime}=\theta_x,\theta_1^{\prime}=\theta_0,\theta_\infty^{\prime}=\theta_\infty} 
~+p_\infty p_x+p_1 p_0 
}
\right. 
$$

\vskip 0.2 cm 
\noindent
 {\large \bf -- For Table 3 ($x=\infty$).}

\vskip 0.2 cm 
\noindent
 {\large  a) \it Parametric formulae of type (\ref{POA})  at $x=\infty$}:
\vskip 0.2 cm 
 Take  the formulae  $ c_i=c_i(\theta_0,\theta_x,\theta_1,\theta_\infty,p_{0x},p_{01},p_{x1})$  of type (\ref{POA}) in  Section \ref{dokidoki} and make  the substitutions 
\be 
\label{ricordo2}
\left\{
\matrix{
\theta_0 \mapsto \theta_0,~~~\theta_x \mapsto \theta_1,~~~\theta_1\mapsto \theta_x,~~~\theta_\infty\mapsto\theta_\infty,~~
\hbox{ (or }
~~
\alpha\mapsto \alpha,~\beta \mapsto \beta,~\gamma \mapsto {1\over 2}-\delta,~\delta \mapsto {1\over 2}-\gamma)
\cr 
\cr 
p_{0x}\mapsto
 p_{01} ,~~~
p_{x1}\mapsto p_{x1},~~~
p_{01} \mapsto -p_{0x}-p_{x1}p_{01}+p_\infty p_1+p_0p_x
}
\right.
\ee

\vskip 0.2 cm 
\noindent 
 {\large b) \it Formulae of type (\ref{POA1}): Monodromy data in terms of the integration constants $c_1$, $c_2$ at $x=\infty$}. 

\vskip 0.2 cm 
Take the formulae of  type (\ref{POA1}) in  Section \ref{dokidoki} and denote them  as in (\ref{biop}), where $c_1$ and $c_2$ are now the constants in  the behaviors of  Table 3 at $x=\infty$.   
Then, the monodromy data corresponding to the behaviors of Table 3 are 
$$
\left\{\matrix{p_{0x} &= &-\bigl(p_{01}^{\prime}
+p_{0x}^{\prime} p_{x1}^{\prime}\bigr)|_{\theta_0^{\prime}=\theta_0,\theta_x^{\prime}=\theta_1,\theta_1^{\prime}=\theta_x,\theta_\infty^{\prime}=\theta_\infty}
~+p_\infty
p_1+p_0 p_x
\cr\cr
p_{x1}&=&p_{x1}^{\prime}~~~~~~~~~~~~~~~~~~~~~~~~~~~~~~~
\cr\cr
p_{01}&=&p_{0x}^{\prime}~~~~~~~~~~~~~~~~~~~~~~~~~~~~~~~
}
\right.
$$

\subsection{Connection Formulae in Closed Form} 
\label{CONECLOSE}

A branch  $y(x)$ has a critical behavior of type (\ref{malamente1}) at a critical point $x=u$, and another behavior of type (\ref{malamente1}) at another critical point $x=v$, where $u\neq v\in\{0,1,\infty\}$. Namley
$$
y(x)=\left\{
\matrix{ y(x,c_1^{(u)},c_2^{(u)}),~~~x\to u
\cr 
\cr
 y(x,c_1^{(v)},c_2^{(v)}),~~~x\to v
}
\right.
$$
 Explicit formulae of type 
\be
\label{closd}
\left\{
\matrix{c_1^{(v)}=c_1^{(v)}(c_1^{(u)},c_2^{(u)})
\cr 
c_2^{(v)}=c_2^{(v)}(c_1^{(u)},c_2^{(u)})
}
\right. ,~~~\hbox{ and the inverse form} ~ \left\{
\matrix{c_1^{(u)}=c_1^{(u)}(c_1^{(v)},c_2^{(v)})
\cr 
c_2^{(u)}=c_2^{(u)}(c_1^{(v)},c_2^{(v)})
}
\right. 
\ee
 are refered to as {\bf connection formulae in closed form}. 
  To obtain them from the parametric form, we proceed as follows.

\vskip 0.2 cm 
{\bf -- 1:} Extract  $c_1^{(u)},c_2^{(u)}$  from the leading terms of $y(x)=y(x; c_1^{(u)},c_2^{(u)})$, $x\to u$.  

\vskip 0.2 cm 

{\bf -- 2:} Use formulae of type (\ref{POA1}) to compute the monodromy data, starting from the constants $c_1^{(u)},c_2^{(u)}$. 

\vskip 0.2 cm 

{\bf -- 3:} If the critical behaviour at $x=v$ is not known, identify it. Refere to the explanations:  `` $\diamond$ How to identify a Critical Behaviour from given Monodromy Data '' in Sections \ref{TABLE}, \ref{TABLE2}  and \ref{TABLE3}. 

\vskip 0.2 cm 

{\bf -- 4:} Substitute the monodromy data in the formulae of type (\ref{POA}) in order to compute the constants $c_1^{(v)},c_2^{(v)}$.

\vskip 0.3 cm
\noindent
{\it Example:} Let $y(x)=(\ref{fullEXP})$ for $x\to 0$. Suppose we need the critical behaviour at $x=1$. 

\vskip 0.2 cm 
-- 1: We extract $(c_1^{(0)},c_2^{(0)})=(a^{(0)},\sigma^{(0)})$ from 
$$
y(x)\sim {c_{1,-1}\over a^{(0)}} x^{1-\sigma^{(0)}},~~~~~x\to 0.
$$

-- 2: We compute the monodromy data $p_{0x}, p_{x1}, p_{01}$, substituting $a^{(0)}$ and $\sigma^{(0)}$ in (\ref{pF}). Note that  $\theta_0,\theta_x,\theta_1,\theta_\infty$ are given by the coefficients of PVI (arbitrary signs do not affect the formulae). 

\vskip 0.2 cm 

-- 3: We check the value of  $p_{x1}$ and eventually of  $p_{0x}$, $p_{01}$, $\Sigma_{\gamma\delta}^k$ and $\Omega_{\alpha\beta}^k$, to identify the critical behaviour at $x=1$. Suppose for example that $p_{x1}\neq \pm 2$, $2\cos \pi \Sigma _{\gamma\delta}^k$,  $-2\cos\Omega_{\alpha\beta}^k$ and $p_{x1}\not <-2$. Thus the behaviour is of type (\ref{fullEXP}.1), namely: 
$$
y(x)\sim 1+{c_{1,-1}^{(1)}\over a^{(1)}}(1-x)^{1-\sigma^{(1)}},~~~~~x\to 1.
$$
 
-- 4: We compute the integration constants through (\ref{aF}), using the monodromy data obtained at point 2 above plus   the substitution (\ref{ricordo1}), namely: 
$$
  \sigma^{(1)}\equiv \sigma_{x1}={1\over \pi}\hbox{arcos}\left({p_{x1}\over 2}\right)
 ,~~~~~a^{(1)}=a(\sigma^{(1)},p_{0x},~-p_{01}-p_{0x}p_{x1}+p_\infty p_x+p_0p_1,~\theta_1,\theta_x,\theta_0,\theta_\infty)
$$
where $a(\sigma,p_{x1},p_{01},\theta_0,\theta_x,\theta_1,\theta_\infty)$ is the function $a$ given in (\ref{aF}). The explicit formulae of $ a^{(1)}=a^{(1)}(a^{(0)},\sigma^{(0)})$ and $\sigma^{(1)}=\sigma^{(1)}(a^{(0)},\sigma^{(0)})$ are as follows. Let  ${\bf G}_i={\bf G}_i(\sigma^{(0)})={\bf G}_i(\sigma^{(0)},\theta_0,\theta_x,\theta_1,\theta_\infty)$, $i=1,...,6$, be the  long expressions appearing in (\ref{pF}). Then
$$
\sigma^{(1)}={1\over \pi}\hbox{arcos}\left[{1\over 2}\left( {{\bf G}_1(\sigma^{(0)})\over a^{(0)}}+{\bf G_2}(\sigma^{(0)})+{\bf G}_3(\sigma^{(0)})~ a^{(0)}\right)\right],~~~~~0\leq \Re \sigma^{(1)}<1,
$$ 
$$
a^{(1)}= {(\theta_x-\theta_1-\sigma^{(1)})(\theta_x+\theta_1-\sigma^{(1)})(\theta_\infty+\theta_0-\sigma^{(1)})\over4 {\sigma^{(1)}}^2(\theta_\infty+\theta_0+\sigma^{(1)})}~{1\over  {\bf F^{(1)}}}~ {U^{(1)}\over V^{(1)}} 
$$
where 
$$
{\bf F}^{(1)}:=  {
\Gamma(1+\sigma^{(1)})^2\Gamma \left({1\over 2}(\theta_1+\theta_x-\sigma^{(1)})+1
\right) \Gamma\left(  
{1\over 2} ( \theta_x-\theta_1-\sigma^{(1)})+1
\right)
\over 
\Gamma(1-\sigma^{(1)})^2 \Gamma \left({1\over2}(\theta_1+\theta_x+\sigma^{(1)})+1
\right) \Gamma\left(   
{1\over 2} ( \theta_x-\theta_1+\sigma^{(1)})+1
\right)
}~\times
$$
$$
\times  {
\Gamma \left({1\over 2}(\theta_\infty+\theta_0-\sigma^{(1)})+1 \right)
 \Gamma\left(  
{1\over 2} ( \theta_0-\theta_\infty-\sigma^{(1)})+1
\right)
\over
\Gamma \left({1\over 2}(\theta_\infty+\theta_0+\sigma^{(1)})+1 \right)
 \Gamma\left( 
{1\over 2} ( \theta_0-\theta_\infty+\sigma^{(1)})+1
\right)
},
$$

$$ 
V^{(1)}:=4 ~s(\theta_1+\theta_x-\sigma) s(\theta_1 - \theta_x+\sigma)~ s (\theta_\infty+\theta_0-\sigma)
 s (\theta_\infty-\theta_0+\sigma).
$$

$$
  U^{(1)}:=-{i\over 2}\left({{\bf G}_4(\sigma^{(0)})\over a^{(0)}}+{\bf G}_5(\sigma^{(0)})+{\bf G}_6(\sigma^{(0)})~
 a^{(0)}+2\cos\pi\sigma^{(0)}~e^{-i\pi\sigma_1}\right)
\sin\pi\sigma^{(1)}+
$$
$$+{1\over 4}(p_xp_0+p_1p_\infty)-{1\over 4}(p_xp_\infty+p_1p_0)e^{-i\pi\sigma^{(1)}}
$$
  $\Box$

\vskip 0.3 cm 
\noindent 
The example above brings the following remarks: 

1.  In general, it is not possible to  invert by direct computations one of the two systems in (\ref{closd}) and  obtain the inverse form.  
In order to do the inversion, one has to follow the procedure 1. $\to$  4.  explained in this section, making use of the parametric connection formulae. For example,   in the example above the system of closed formulae  
$$\left\{\matrix{
\sigma^{(1)}=\sigma^{(1)}(a^{(0)},\sigma^{(0)})
\cr 
a^{(1)}=a^{(1)}(a^{(0)},\sigma^{(0)}) 
}
\right.
$$ is not directly solvable in order to obtained  $\sigma^{(0)}=\sigma^{(0)}(a^{(1)},\sigma^{(1)})$, $ a^{(0)}=a^{(0)}(a^{(1)},\sigma^{(1)})$.  This is because one should invert several combinations of trigonometric and $\Gamma$ functions. Instead, one can to use the parametric formulae following the steps 1. $\to$  4. above, starting form $a^{(1)},\sigma^{(1)}$ and reaching   $\sigma^{(0)}=\sigma^{(0)}(a^{(1)},\sigma^{(1)})$. $ a^{(0)}=a^{(0)}(a^{(1)},\sigma^{(1)})$.

2. A branch $y(x)$ with a given critical behaviour at $x=u$ may have one of the different types of behaviours at $x=v$, according to the monodromy data. Therefore,  in order to write the connection formulae in closed form,  a big number  of  combinations is to be considered. The size of this paper cannot contain them.  Nevertheless,  all the combinations can be obtained from the parametric connection formulae of Section \ref{dokidoki} and the procedure explained in this section.

 For this reason we have only provided the parametric connection formulae at $x=0$, which contain the necessary tools to construct 
 all the parametric connection formulae at $x=1,\infty$ and the connection formulae in closed form. 


\section{Detailed description of Table 1 ($x\to0$)}
\label{esplaintable}

We explain Table 1, and show how it is obtained from the results of \cite{Jimbo},  \cite{Boalch},  \cite{D4}, \cite{D3},  \cite{D2}, \cite{D1}, \cite{guz2010}. 
In general, substitution into PVI of the leading term of a critical behaviour  provides the entire expansion. An example is the procedure of \cite{guz2010}, section 7. This procedure has the disadvantage that the coefficients of higher order are computed recursively from those of lower order. 
A conjectural close form of the coefficients  at all orders for the $\tau$ function of (\ref{fullEXP})  (and possibly applicable, by a limit procedure, to (\ref{atopy}), (\ref{log1}), (\ref{log1zero}) and (\ref{logsquare})) is given in \cite{Lisovyy}. 

\vskip 0.2 cm 
We also compare the behaviours in Table 1, which is obtained in the framework of the isomonodromy deformations method,    with the  results of power geometry   \cite{Bruno7}. We show that they coincide. See the paragraphs  {\bf {\bf (*)}} and sub-section \ref{COMOPAIR}. 

\vskip 0.3 cm
{\bf 1)} Two-complex parameter ($a$, $\sigma\in{\bf C}$) solutions \cite{Jimbo},\cite{D3},\cite{D2},\cite{D1}, \cite{guz2010}, \cite{Boalch}:
   \be
\label{fullEXP}
\sq{
  y(x)= \sum_{n=1}^\infty x^n\sum_{m=-n}^n c_{nm}
  (ax^\sigma)^m,~~~-1< \Re \sigma<1,~~a\neq 0
}
\ee 
The family is invariant for $\sigma\mapsto -\sigma$, therefore we can restrict to the range 
$$
\sq{
0\leq \Re \sigma <1
}
$$
The $c_{nm}$'s are rational functions of  $\sigma$, and algebraic functions of the coefficients of PVI.  In particular we normalise in such a way that $c_{11}=1$. 
If $\sigma$ is not  (\ref{SeT}) or (\ref{iku}), we have 
\be
\label{c11}
c_{1,-1}={\Bigl[(\sqrt{-2\beta}-\sqrt{1-2\delta})^2-\sigma^2\Bigr]  \Bigl[(\sqrt{-2\beta}+\sqrt{1-2\delta})^2-\sigma^2\Bigr]\over 16 \sigma^4 
}
\ee
\be
\label{c10}
c_{10}={\sigma^2-2\beta+2\delta-1\over 2\sigma^2}
\ee
and $c_{nm}=q_{nm}(\sqrt{\alpha},\sqrt{\beta},\sqrt{\gamma},\sqrt{1-2\delta};\sigma),$ where $q_{nm}$ are rational functions of their arguments. 
If  $\Re\sigma>0$,  the dominant term is 
$$
y(x)={c_{1,-1}\over a} x^{1-\sigma}\Bigl(1+O(x^\sigma,x^{1-\sigma})\Bigr)
$$
If $\sigma=-2i\nu$, $\nu\in {\bf R}$,  solution (\ref{fullEXP})  becomes   a three-real parameter solution 
\be
\label{lantern}
\sq{
y(x)= x\Bigl[A\sin(2\nu \ln x +\phi)+B+ \sum_{n=1}^\infty x^n\sum_{m=-n-1}^{n+1} c_{nm} (e^{i\phi}x^{-2i\nu})^m\Bigr]
} 
\ee
$$
B^2-A^2={\beta\over 2\nu^2},~~~~~B={4\nu^2+2\beta+1-2\delta\over 8\nu^2},
$$
$$
c_{nm}=\tilde{q}_{nm}(\sqrt{\alpha},\sqrt{\beta},\sqrt{\gamma},\sqrt{1-2\delta};\nu)
$$
where $\tilde{q}_{nm}$  are rational functions of their arguments, and the role of $a$  is played by $\phi$ (the log of $a$). 
\vskip 0.2 cm 
\noindent
{\bf Remark 4:} If $\sigma$ is not (\ref{SeT}) or (\ref{iku}), and $\Re\sigma>0$, we can also write (\ref{fullEXP}) as
$$
 y(x)= \sum_{n=0}^\infty x^{2n+1}\sum_{m=-n}^\infty a_{nm} x^{m\rho} +\sum_{n=1}^\infty x^{2n}\sum_{m=-n}^\infty a_{nm} x^{m\rho},~~~~~\rho=1-\sigma,~~a_{nm}\in{\bf C}
$$
\be
\label{ABru}
= v_{\rho}x^\rho \Bigl(1 + \sum_{n=1}^\infty x^{n\rho} +\sum_{m=1}^\infty \sum_{n= n_0}^\infty
v_{nm}x^{m(1-\rho)+n\rho}\Bigr), ~~~v_\rho,~v_{nm}\in {\bf C},~v_\rho\neq 0
\ee
where $ n_0=(m-1)/2$ if $m$ is odd,  $n_0= m/2 -1$ if $m$ is even.

\vskip 0.2 cm 

The leading term of the  expansions (\ref{fullEXP}) was obtained first by Jimbo in \cite{Jimbo}. It can also be obtained with the matching procedure of \cite{D2}, making use of Fuchsian reductions of the system (\ref{SYSTEM}). The full expansion (\ref{fullEXP})  can also be obtained by  local analysis, with the method of  \cite{Sh}, \cite{D3}. This method proves also the convergence of the expansions,  which therefore  define true solutions of PVI, which are {\it basic solutions}. The $c_{nm}$'s are determined by the procedure of \cite{guz2010}, section 7.

\vskip 0.3 cm 
{\bf {\bf (*)}} The class of the solutions (\ref{fullEXP}) when $\Re \sigma\neq 0$ is the  class ${\cal A}_0$ of \cite{Bruno7}. In  \cite{Bruno7} the ${\cal A}_0$-solutions are written as (\ref{ABru}), though the range of $n,m$ is looser, namely $n+m>0$, $n,m\geq 0$.

 \vskip 0.2 cm 
{\bf {\bf (*)}} The image of $\bigl({\cal B}_0^{+}\cup{\cal B}_0^{-}\bigr)\cup\bigl({\cal B}_7^{+}\cup{\cal B}_7^{-}\bigr)$ of \cite{Bruno7} through the symmetry (\ref{sym2})  is in the sub class of solutions (\ref{lantern}).

\vskip 0.3 cm 
\noindent
{\bf $\diamond$------------------------------------------------------------------------------------------------------------------$\diamond$}

\vskip 0.3 cm
{\bf 2)}  If the coefficients are such that  \fbox{$\sqrt{-2\beta}+\sqrt{1-2\delta}\not \in {\bf Z}$} or  \fbox{$\sqrt{-2\beta}-  \sqrt{1-2\delta}\not \in {\bf Z}$} (namely, $\theta_0+\theta_x$ or $\theta_0-\theta_x\not \in {\bf Z}$), PVI admits a one-complex parameter $a\in{\bf C}$ family of solutions (see \cite{D2} and \cite{guz2010})
\be
\label{atopy}
\sq{
y(x)= \sum_{n=1}^\infty x^n\sum_{m=0}^n c_{nm} (ax^{\sigma})^m
=
\sum_{k=0}^\infty y_k(x)~(ax^\sigma)^k,~~~ \sigma\in\Sigma_{\beta\delta},~~ \Re\sigma>-1
}
\ee
 The $c_{nm}$ are functions 
$q_{nm}(\sqrt{\alpha},\sqrt{\beta},\sqrt{\gamma},\sqrt{1-2\delta}),$   like in the case of  (\ref{fullEXP}). In particular $c_{10}= {\sqrt{-2\beta}\over \sqrt{-2\beta}\pm\sqrt{1-2\delta}}$.  There are two possible solutions, corresponding to the signs $\pm$ above, namely to the fact that $\sigma\in \Sigma_{\beta,\delta}^{k}$, $k=1,2$.  

\noindent
The functions  $y_k(x)$ are convergent Taylor expansions $y_l(x)=O(x^l).$ 

\noindent
One of the following two terms is the leading one
$$
 y(x)={\sqrt{-2\beta}\over \sqrt{-2\beta}\pm\sqrt{1-2\delta}}~x+ax^{1+\sigma}+...
$$
If $\beta=0$, the leading term is  $ax^{1+\sigma}$ and 
$$
 y(x)=ax^{1+\sigma}+\sum_{n=2}^\infty x^n\sum_{m=1}^n c_{nm}(ax^{\sigma})^m
$$
The solution (\ref{atopy}) reduces to the following Taylor series  when $a=0$ 
\be
\label{davidekan}
y_0(x)={\sqrt{-2\beta}\over \sqrt{-2\beta}+(-)^k\sqrt{1-2\delta}}x+\sum_{n=2}^\infty q_{n0}(\sqrt{\alpha},\sqrt{\beta},\sqrt{\gamma},\sqrt{1-2\delta})x^n,
\ee
$$
~~~y_0(x)=0 \hbox{ if } \beta=0.
$$

The  local analysis of \cite{Sh} and \cite{D3}  extends the  result of \cite{Jimbo} in  the case $\Re\sigma\geq 1$, $\sigma\in\Sigma_{\beta\delta}$, and proves existence and convergence of solutions (\ref{atopy}). 

 Note that  (\ref{atopy}) is a subcase of (\ref{fullEXP}) when  $-1<\Re\sigma<1$, because the coefficients of $x^{-|m|\sigma}$ in (\ref{fullEXP}) vanish when $\sigma$ is (\ref{SeT}) or (\ref{iku})
 (see  formula (\ref{c11})) provided that $\sqrt{-2\beta}\pm \sqrt{1-2\delta}\not \in {\bf Z}$ (otherwise coefficients diverge). Also note that this happens when we choose $c_{11}\neq 0$ as integration constant. One can choose $c_{1,-1}$ instead, and express $c_{11}$ from (\ref{c11}). In this case, again the coefficients of $x^{|m|\sigma}$ vanish for $\sigma$ in the set (\ref{SeT}) and we again obtain (\ref{atopy}) with $\sigma\mapsto -\sigma$, $\Re(-\sigma)>-1$.

\vskip 0.3 cm 
{\bf {\bf (*)}} The symmetry (\ref{sym2}) transforms the set $\bigl({\cal B}_1(\hbox{non log})\cup
{\cal B}_2(\hbox{non log})\bigr)\cup\bigl({\cal B}_1^{+}\cup{\cal B}_1^{-}\bigr)\cup\bigl({\cal B}_2^{+}\cup{\cal B}_2^{-}\bigr)\cup
{\cal B}_6\cup\bigl({\cal B}_6^{+}\cup{\cal B}_6^{-}\bigr)$ of \cite{Bruno7} into the  solutions (\ref{atopy})  with $\beta \neq 0$. It transforms ${\cal C}_0^\infty \cup{\cal B}_8\cup{\cal B}_9$ of \cite{Bruno7}  into the  solutions (\ref{atopy}) with $\beta =0$. This will be clear after reading {\bf 2.1)} below.  [Here we use ${\cal B}_i$(non log) to denote the subclass of ${\cal B}_i$ such that the coefficients of powers of $x$ are constants, and the power of $x$ have no integer exponent (see {\bf 2.1)} below). We  distinguish it from the subclass ${\cal B}_i$(log) whose coefficients are  polynomials of $\ln x$, and the powers of $x$ have integer exponents (for some value of the coefficients of PVI the coefficients may be constants, and the solutions become Taylor series. See {\bf 3.1)} below). ${\cal B}_i={\cal B}_i$(non log)$\cup{\cal B}_i$(log).]

We conclude that the class   ${\cal A}_0\cup\bigl({\cal B}_0^{+}\cup{\cal B}_0^{-}\bigr)\cup\bigl({\cal B}_1(\hbox{non log})\cup
{\cal B}_2(\hbox{non log})\bigr)\cup\bigl({\cal B}_1^{+}\cup{\cal B}_1^{-}\bigr)\cup\bigl({\cal B}_2^{+}\cup{\cal B}_2^{-}\bigr)\cup
{\cal B}_6\cup\bigl({\cal B}_6^{+}\cup{\cal B}_6^{-}\bigr)\cup\bigl({\cal B}_7^{+}\cup{\cal B}_7^{-}\bigr)\cup{\cal B}_8\cup{\cal B}_9\cup{\cal C}_0^\infty $ of \cite{Bruno7} is  contained in the union of  {\bf 1)} $ \cup $ {\bf 2)} 
$\cup~\{$the image of  {\bf 1)} $ \cup $ {\bf 2)} 
     through the symmetry (\ref{sym2})$\}$.

\vskip 0.3 cm 
\noindent
{\bf $\diamond$------------------------------------------------------------------------------------------------------------------$\diamond$}

\vskip 0.3 cm

{\bf 3)} If \fbox{$   \sqrt{-2\beta}+\sqrt{1-2\delta}=N$} or  \fbox{$\sqrt{-2\beta}-\sqrt{1-2\delta}=N$} 
(namely, $\theta_0+\theta_x=N$ or $\theta_0-\theta_x=N$), $N\in {\bf Z}$, then 
solution (\ref{atopy}) is no longer defined, because some coefficients diverge. We have instead  solutions  
 formally like  (\ref{atopy}), where $\sigma$ is the integer $N$ and the divergent coefficients $c_{nm}$ 
are replaced by  polynomials of $\ln x$.  Namely, from the results of  \cite{D2},\cite{D1}, we see that there exist one-complex parameter ($a\in{\bf C}$) expansions 
\be
\label{log1}
\sq{
  y(x)= \Sigma_{n=1}^{|N|}b_nx^n+\Bigl(a+b_{|N|+1}\ln x\Bigr)x^{|N|+1}+\Sigma_{n= |N|+2}^\infty P_n(\ln x;a)x^n,
~N\neq 0
}
\ee
\be
\label{log1zero}
 \sq{
 y(x)= \Bigl(a\pm\sqrt{-2\beta}\ln x\Bigr)x+\sum_{n=2}^\infty P_n(\ln x;a)x^n,~~\left.\matrix{N= 0\cr 2\beta=2\delta-1}\right. 
}
\ee
where $b_n$ are rational functions of $\sqrt{\alpha},\sqrt{\beta},\sqrt{\gamma}$, and $\sqrt{1-2\delta}$; $a$ is the free complex parameter; $P_n(\ln x;a)$ are polynomials of $\ln x$ of degree $n-|N|$, with coefficients that are  rational functions of  $\sqrt{\alpha},\sqrt{\beta},\sqrt{\gamma},\sqrt{1-2\delta}$ and $a$.
 Solution (\ref{log1zero}) is {\it a basic solution}. 

(\ref{log1}) and (\ref{log1zero})  reduce to  convergent Taylor series  if either $\sqrt{-2\beta}\in \{0,-1,-2,...,N\}$ for $N<0$, $\sqrt{-2\beta}\in \{0,1,2,...,N\}$ for $N>0$; or if $\{\sqrt{2\alpha}+\sqrt{2\gamma},\sqrt{2\alpha}-\sqrt{2\gamma}\}\cap {\cal N}_N\neq \emptyset$, where $
{\cal N}_N$ is (\ref{infinitapalla}). 
The Taylor series  are 
\be
\label{TTLO1}
 \sq{
 y(x)= \sum_{n=1}^{|N|}b_nx^n+ax^{|N|+1}+\sum_{n= |N|+2}^\infty b_n(a)x^n, ~~~N\neq 0
}
\ee
\be
\label{TTLO2}
 \sq{
 y(x)= y(x)=ax+{a(a-1)\over 2}(2\gamma-2\alpha-1)x^2+\sum_{n=3}^\infty b_n(a)x^n, ~~~
\left.
\matrix{N = 0
\cr
2\beta=2\delta-1=0
}
\right.
}
\ee 
$b_n(a)$ are   rational functions of $\sqrt{\alpha},\sqrt{\beta},\sqrt{\gamma},\sqrt{1-2\delta}$ and $a$. 
In particular, (\ref{TTLO2}) is the limit of (\ref{log1zero}) for $\beta \to 0$.

\vskip 0.2 cm 
In \cite{D2}, the  solution (\ref{log1zero}) is obtained with the matching procedure, making use of Fuchsian  reductions of (\ref{SYSTEM}). Convergence of (\ref{log1}) and (\ref{log1zero}) is not proved, but we expect it, the proof being possibly based on an implementation of the method of \cite{Sh} and \cite{D3}. All the coefficients can be recursively determined by direct substitution into PVI. 

Three basic Taylor series (\ref{T1coe}), (\ref{TTLO4}) and   (\ref{T2coe}), given below, are obtained in \cite{D2} by non-Fuchsian reductions of (\ref{SYSTEM}).  
The symmetry (\ref{sym2}) transforms  (\ref{TTLO4}) and   (\ref{T2coe}) respectively  into (\ref{TTLO2}) and the following: 
\be
\label{batsta}
y(x)= \sqrt{-2\beta} ~x~+a~x^2~+\sum_{n=3}^\infty b_n(a;\sqrt{\alpha},\sqrt{\beta})x^n,~~~~~\beta\neq 0,~~ (\sqrt{-2\beta}\pm\sqrt{1-2\delta})^2=1,~~\alpha=\gamma
\ee 
The last solution is included in (\ref{TTLO1}) for $N=1$. 
 Solutions (\ref{TTLO2}), (\ref{batsta}) and (\ref{davidekan})  are {\it basic Taylor solutions}, namely all  the other Taylor  solutions at $x=0$ can be obtained from them by bi-rational transformation.  For example, when $N=1$, solutions (\ref{TTLO1}) contain two nonequivalent families: the family  (\ref{batsta}) and 
$$
y(x)= x+ax^2+{1\over 2}a(\gamma-\alpha)x^3+O(x^4),~~~\delta=-\beta={1\over 2}.
$$
The latter is equivalent to (\ref{TTLO2}) via bi-rational transformation, while (\ref{batsta}) is not. Taylor expansions are also obtained in \cite{kaneko}, and their convergence is proved.

\vskip 0.2 cm 
{\bf {\bf (*)}} The symmetry (\ref{sym2}) transforms the class ${\cal B}_1(\hbox{log})\cup
{\cal B}_2(\hbox{log})$ of \cite{Bruno7} into (\ref{log1}), and ${\cal B}_4\cup{\cal B}_5$ into (\ref{log1zero}).

\vskip 0.3 cm 
\noindent
{\bf $\diamond$------------------------------------------------------------------------------------------------------------------$\diamond$}

\vskip 0.3 cm
{\bf 4)} If \fbox{$2\beta\neq 2\delta-1$} (namely, $\theta_0\pm\theta_x\neq 0$), besides log-solutions (\ref{log1}), there are one-complex parameter ($a\in {\bf C}$) logarithmic solutions ( see \cite{Jimbo}, \cite{D2} and  \cite{D1}), which are also {\it a basic solutions}: 
\be
\label{logsquare}
\sq{
y(x)=\left[{2\beta+1-2\delta\over 4}(a+\ln x)^2+{2\beta\over 2\beta+1-2\delta}\right]x+\sum_{n\geq 2}^\infty P_n(\ln x;a)x^n,~~~a\in{\bf C}.
}
\ee
The above solutions are obtained in \cite{D2}, \cite{D1} by the matching procedure, with Fuchsian reductions of system (\ref{SYSTEM}). As far as convergence is concerned, no proof is published, though we expect that it is convergent and an implementation of the method of \cite{Sh} should allow to prove this. All the coefficients can be recursively determined by direct substitution into PVI. 

\vskip 0.2 cm 
{\bf {\bf (*)}} The symmetry (\ref{sym2}) transforms the class ${\cal B}_3$ of \cite{Bruno7} into (\ref{logsquare}).

\vskip 0.3 cm
\centerline{\bf $\diamond$=================================================$\diamond$}
\vskip 0.3 cm
The next step is the application of the symmetry (\ref{sym2}) to the solutions  {\bf 1)}-{\bf 4)}, to obtain new solutions. In all which follows, it is understood that wherever $\alpha$, $\beta$, $\gamma$ and $\delta$ appear in  formulae related to {\bf 1)}-{\bf 4)} (for example, the coefficients of the expansions), we must substitute  $(\alpha,\beta,\gamma,\delta)\mapsto 
 (-\beta,-\alpha,{1\over 2}-\delta, {1\over 2}-\gamma)
$.

\vskip 0.3 cm 
{\bf 1.1)}  The symmetry (\ref{sym2}) transforms solutions of case 2) into one parameter ($a\in {\bf C}$) solutions,
\vskip 0.3 cm 

{\bf -- 1.1a)} The symmetry (\ref{sym2}) transforms  (\ref{fullEXP}) with $0<\Re \sigma<1$ into
$$
y(x)= {1
\over 
\sum_{m=0}^\infty x^n \sum_{m=-n-1}^{n+1} c_{n+1,m} (ax^{\sigma})^m
}
={x^\sigma/c_{1,-1}
\over 
 1+{c_{10} (ax^\sigma) \over a c_{1,-1}}  +{(ax^{\sigma})^2\over a c_{1,-1}}  +\sum_{n=1}^\infty x^n \sum_{m=-n}^{n+2} {c_{n+1,m-1}\over a c_{1,-1}} (ax^{\sigma})^m}
$$
$$
= {1\over c_{1,-1}} x^\sigma \sum_{n=0}^\infty x^n \sum_{m=-n}^\infty b_{nm} x^{m\sigma},~~~~~b_{nm}\in{\bf C}
$$
Namely
\be
\label{fullEXP1}
y(x)= \sum_{n=1}^\infty x^n \sum_{m=-n}^\infty d_{nm} x^{m\rho},~~~~~\rho=1-\sigma ,~~0<\Re\rho<1,~~~d_{nm}\in{\bf C}
\ee
Explicit computation yields $
d_{nm}=0$ for any $m>n$.  
Note that $1-\rho$ is not (\ref{SeT1}). Thus,  (\ref{fullEXP1}) is a solution (\ref{fullEXP}).

\vskip 0.2 cm 
{\bf {\bf (*)}} We have re-obtained elements in the class ${ \cal A}_0$ of \cite{Bruno7}.

\vskip 0.2 cm 
\noindent
{\it Note:} To verify that $d_{nm}=0$ for $m>n$, one can do the  explicit computation of the $d_{nm}$ by expanding  $ x/y(x)$ and computing the coefficients. Note that to get the correct terms in the truncation  $\sum_{n=1}^N x^n \sum_{m=-n}^n d_{nm} x^{m\rho}$, one needs to start from (\ref{fullEXP}) truncated at  $\sum_{n=1}^{N+3} x^n \sum_{m=-n}^n c_{nm} x^{m\rho}$.  One can also 
conclude directly from the form (\ref{fullEXP1}) that $d_{nm}=0$ for $m>n$. Actually, such a form, substituted into (PVI), implies $d_{nm}=0$ for $m>n$. Namely, we always have that: 
$$
y(x)= \sum_{n=1}^\infty x^n \sum_{m=-\infty}^\infty d_{nm} x^{m\rho}
~~~+~\hbox{ substitution into (PVI) }~~\Longrightarrow~~~
y(x)= \sum_{n=1}^\infty x^n \sum_{m=-n}^n d_{nm} x^{m\rho}
$$

\vskip 0.3 cm

\vskip 0.2 cm

{\bf -- 1.1b)} The symmetry (\ref{sym2}) transforms  (\ref{lantern}) into the following   three-real parameters ($\nu\in{\bf R}$, $\phi\in{\bf C}$) solution (see \cite{guz2010}):
\be
\label{lantern1}
\sq{
y(x)= \Bigl[\sum_{n=0}^\infty x^n\sum_{m=-n-1}^{n+1} c_{nm}~
e^{im\phi}x^{2im\nu}\Bigr]^{-1}= {1\over A\sin(2\nu \ln x +\phi)+B+ O(x)}
}
\ee
where
$$
O(x)= \sum_{n=1}^\infty x^n\sum_{m=-n-1}^{n+1} c_{nm} e^{im\phi}x^{2im\nu},~~~~~
A=-\sqrt{{\alpha\over 2\nu^2}+B^2},~~~B={2\nu^2+\gamma-\alpha\over 4\nu^2}
$$ 
The expansion is convergent. 
All $c_{nm}$, are rational functions of $\nu$, and $     \sqrt{\alpha},\sqrt{\beta},\sqrt{\gamma},\sqrt{1-2\delta}  $. 
If $\arg x$ is bounded (namely, when we considered  $y(x)$ as a  branch), then (\ref{lantern1}) admits two infinite sequences of {\it movable poles}  in the neighbourhood of $x=0$, asymptotically distributed along two rays and accumulating at $x=0$ (see \cite{Dpoli}). One can do a formal (non convergent) expansion, when $x$ is far from the poles, as 
\be
\label{formal}
y(x)=
 {1\over A\sin(2\nu \ln x +\phi)+B}+\delta_\nu x^{2i\nu}
\left\{
1+\sum_{n=0}^\infty x^n\sum_{m=-n}^\infty \delta_{nm}x^{2im\nu}
\right\},~~~\delta_\nu,~\delta_{nm}\in{\bf C}
\ee

\vskip 0.2 cm 
{\bf {\bf (*)}} The class of solutions (\ref{lantern1}) coincide with the union  $\bigl({\cal B}_0^{+}\cup{\cal B}_0^{-}\bigr)\cup\bigl({\cal B}_7^{+}\cup{\cal B}_7^{-}\bigr)$ of \cite{Bruno7}, where a formal expansion,  in a  form equivalent to (\ref{formal}), is given. 

\vskip 0.3 cm 
\noindent
{\bf $\diamond$------------------------------------------------------------------------------------------------------------------$\diamond$}

\vskip 0.3 cm 
{\bf 2.1)}  The symmetry (\ref{sym2}) transforms solutions of case 2) into one parameter ($a\in {\bf C}$) solutions, defined when at least one of  \fbox{$\sqrt{2\alpha}+\sqrt{2\gamma}\not\in{\bf Z}$} or  \fbox{$\sqrt{2\alpha}-\sqrt{2\gamma}\not\in{\bf Z}$ } (namely $\theta_\infty+\theta_1\not\in{\bf Z}$ or  $\theta_\infty-\theta_1\not\in{\bf Z}$). We distinguish the two cases  $\alpha\neq 0$ and $\alpha=0$.

\vskip 0.3 cm 
\noindent 
$\diamond$ Case \fbox{$\alpha\neq 0$}.

\vskip 0.3 cm 

\vskip 0.3 cm 
-- {\bf 2.1a)}  The symmetry (\ref{sym2})  transforms (\ref{atopy}) with $\Re \sigma>0$ into 
$$
y(x)=
{1
\over
c_{10}
\left(
1+{a\over c_{10}} x^\sigma +\sum_{n=1}^\infty x^n \sum_{m=0}^{n+1}{c_{n+1,m}\over c_{10}}(ax^{\sigma})^m
\right)
}
= \sum_{n=0}^\infty x^n\sum_{m=0}^\infty b_{nm}x^{m\sigma}, 
$$
where $b_{01}$ is the  arbitrary constants, and the other coefficients are rational functions of  $\sqrt{\alpha},\sqrt{\beta},\sqrt{\gamma},\sqrt{1-2\delta}$ and $b_{01}$.  In particular $b_{00}= {\sqrt{\alpha} \pm \sqrt{\gamma} \over \sqrt{\alpha}}$. With the substitution $\sigma=\rho+1$, the above  becomes  
\be
\label{UUU}
\sq{
y(x)=d_{00}+ \sum_{n=1}^\infty x^n \sum_{m=0}^n d_{nm}(\tilde{a}x^{\rho})^m=
\sum_{l=0}^\infty y_l(x)(\tilde{a}x^{\rho})^l
}
\ee
$$
  \sqs{\rho+1\in \Omega_{\alpha\gamma}^{k} \hbox{ given by   (\ref{STaR}), }k=1,2, ~~\Rightarrow~~\Re\rho>-1}
$$
$$
d_{00}={\sqrt{\alpha} +(-)^k \sqrt{\gamma} \over \sqrt{\alpha}} ,~~~d_{11}=1,~~~\tilde{a}=-a~(d_{00})^2,~~~~~\alpha\neq 0.
$$
 Here $a$ is the arbitrary constant,  $d_{nm}\in{\bf C}$ are rational functions of $\sqrt{\alpha},\sqrt{\beta},\sqrt{\gamma}$, and $\sqrt{1-2\delta}$. The $y_l(x)=O(x^l)$'s are Taylor series. In particular  
\be
\label{T1coe}
y_0(x)= d_{00}
+
\sum_{n=1}^\infty d_{n0}x^n,
\ee
$$
y_1(x)=-d_{00}x+\sum_{n=2}^\infty d_{n1}x^n
$$
Note that $y_1(x)$  can be directly obtained from (\ref{davidekan}) by the symmetry (\ref{sym2}). Its convergence is proved in \cite{kaneko}. 

\vskip 0.2 cm 
{\bf {\bf (*)}} In terms of the classification of  \cite{Bruno7}, the class of solutions (\ref{UUU}) coincides with the union of the class ${\cal B}_6$  with the subclass of ${\cal B}_1\cup{\cal B}_2$  which does not contain logarithmic coefficients.

\vskip 0.3 cm 

\vskip 0.3 cm

-- {\bf 2.1b)} The symmetry (\ref{sym2})  transforms (\ref{atopy}) with  $\sigma=-2i\nu$, $\nu\in{\bf R}$,  into 
\be
\label{TAU}
\sq{
\matrix{
y(x)= \Bigl\{ \sum_{n=0}^\infty x^n\tau_n(x)\Bigr\}^{-1}, &
 \tau_n(x)= \sum_{m=0}^{n+1}c_{n+1,m}(ax^{-2i\nu})^m, &\nu\in{\bf R} &\nu\neq 0
\cr 
\cr 
&\tau_0(x)= {\sqrt{\alpha} \over \sqrt{\alpha} +(-)^k \sqrt{\gamma}}+ax^{-2i\nu}&&
}
}
\ee 
$$
  \sqs{2i\nu\in\bigl\{ (\sqrt{2\alpha}+(-)^k\sqrt{2\gamma}),-(\sqrt{2\alpha}+(-)^k\sqrt{2\gamma})\bigr\}\cap i{\bf R}\neq \emptyset,~~k=1,2.}
$$
Solutions  with  $2i\nu =\pm(\sqrt{2\alpha}+\sqrt{2\gamma})$ occur when $\sqrt{2\alpha}+\sqrt{2\gamma}$ is imaginary, with   
$2i\nu=\pm(\sqrt{2\alpha}-\sqrt{2\gamma})$ when $\sqrt{2\alpha}-\sqrt{2\gamma}$ is imaginary.   
The coefficients 
$c_{n,+1,m}$'s are rational functions of $\sqrt{\alpha},\sqrt{\beta},\sqrt{\gamma}$ and $\sqrt{1-2\delta}$.  
The denominator $\sum_{n=0}^\infty x^n\tau_n(x)$ may vanish, so $y(x)$ may have movable poles. In a neighbourhood of $x=0$ they are asymptotically distributed along the ray of the zeros of $\tau_0(x)$. The latter are given by 
$$|x|=\exp
\left\{-
{1\over 2\nu}
\left(
\hbox{arg}{c_{10}\over a}+(2l+1)\pi
\right)
\right\}
,
~~~~~\hbox{arg} ~x= {1\over 2\nu} \ln \left|{c_{10}\over a}\right|
$$
$$
c_{10}= {\sqrt{\alpha} \over \sqrt{\alpha} \pm \sqrt{\gamma}},~~~l\in{\bf Z}, ~~l\geq l_0>0
$$
Far from the line of the poles, we can write 
\be
\label{TAU1}
 y(x)= {1\over \tau_0(x)}\Bigl[1-{\tau_1(x)\over \tau_0(x)}+O(x^2)\Bigr],~~~{\tau_1(x)\over \tau_0(x)}=O(x)
\ee

\vskip 0.2 cm 
{\bf {\bf (*)}} The class (\ref{TAU}),  expanded as (\ref{TAU1}),  coincides with the union  $\bigl({\cal B}_1^{+}\cup   {\cal B}_1^{-}\bigr)\cup  \bigl({\cal B}_2^{+}\cup  {\cal B}_2^{-}\bigr)\cup  \bigl({\cal B}_6^{+}\cup  {\cal B}_6^{-}\bigr)$ of \cite{Bruno7}.

\vskip 0.3 cm 

\vskip 0.3 cm

--{\bf 2.1c)}  The symmetry (\ref{sym2})  transforms (\ref{atopy}) with $-1<\Re\sigma<0$  into
$$
y(x)= {{1\over a} x^{-\sigma}\over 1+{c_{10}\over a}x^{-\sigma}+\sum_{n=1}^\infty x^n \sum_{m=-1}^n {c_{n+1,m+1}} (ax^\sigma)^m}
$$
$$
= {1\over a}x^{-\sigma} \sum_{n=0}^\infty x^n \sum_{m=-n}^\infty b_{nm} x^{-m\sigma}
= \sum_{n=1}^\infty x^n \sum_{m=-n}^\infty d_{nm}x^{m(1+\sigma)}  
$$
This is a solution in the class (\ref{fullEXP}),  with new $\tilde{\sigma} =1+\sigma$, $0<\Re\tilde{\sigma}<1$.
\vskip 0.2 cm 
{\bf {\bf (*)}} These  solutions belong to ${\cal A}_0$ of \cite{Bruno7}.

\vskip 0.5 cm 

\vskip 0.3 cm
\noindent
$\diamond$ {\bf Case} \fbox{$\alpha=0$}.

\vskip 0.3 cm 

\vskip 0.3 cm

-- {\bf 2.1d)}  The symmetry (\ref{sym2}) transforms solutions (\ref{atopy})  with $\beta=0$ into the following solutions  with $\alpha=0$:
\be
\label{prediv}
y(x)={1\over 
ax^\sigma+\sum_{n=1}^\infty x^n \sum_{m=1}^{n+1} c_{n+1,m}(ax^{\sigma})^m}
=\left.
{1\over a x^\sigma
\left(
1+\sum_{n=1}^\infty x^n \sum _{m=0}^n {c_{n+1,m+1}} (ax^{\sigma})^m
\right)
}\right. 
\ee
We rename $\sigma$ with the letter $\omega$, namely
\be
\label{div}
\sq{
y(x)={1\over a}  x^{-\omega} 
\left(
1+\sum_{n=1}^\infty x^n \sum_{m=0}^n d_{nm}(a x^\omega)^m
\right)={1\over a} x^{-\omega} \sum_{k=0}^\infty y_k(x)(ax^{\omega})^k}
\ee
$$
\omega=\sqrt{2\gamma}~\hbox{sgn}(\Re\sqrt{2\gamma}) ~~\hbox{ or }~~ 
\omega \in \{\sqrt{2\gamma},-\sqrt{2\gamma}\} ~\hbox{ if } ~ -1<\Re\sqrt{2\gamma}<1; ~~~\sqs{\sqrt{2\gamma}\not\in{\bf Z}} 
$$ 
$a$ is the free parameter, $\Re \omega>-1$,  $d_{nm}$ are rational functions of
 $\sqrt{\alpha},\sqrt{\beta},\sqrt{\gamma}$ and $\sqrt{1-2\delta}$. $y_l(x)=O(x^l)$ are Taylor series, 
being $y_0(x)=1+O(x)$ .
 When  $\omega =\pm \sqrt{2\gamma}\in i{\bf R}$, (\ref{prediv}) can be formally expanded as (\ref{div}) only far from the poles, namely the zeros of the denominator of   (\ref{prediv}), which possibly exist for any small value of $|x|$.   This case falls into the case    (\ref{TAU}),  $\alpha=0$.  
   When $-1<\Re\omega<0$, the solutions (\ref{div}) are in the class (\ref{fullEXP}) and vanish for $x\to 0$. 
When $ \Re\omega> 0 $, $|y(x)|$  diverges for $x\to 0$. In  Table 1 we have tabulated (\ref{div})  for  the divergent case $\Re \omega>0$, because the other cases are already included in  (\ref{TAU}) and (\ref{fullEXP}). This explains the condition  
\be
\label{etardi}
\sqs{
\omega=\sqrt{2\gamma}~\hbox{sgn}(\Re\sqrt{2\gamma})\not\in{\bf Z},~~~\gamma\not\in (-\infty,0]
}
\ee
\vskip 0.2 cm 
 
{\bf {\bf (*)}} The class of solution (\ref{div}), with restriction (\ref{etardi}),  coincides with ${\cal C}_0^\infty\cup {\cal B}_8\cup{\cal B}_9$ 
of \cite{Bruno7}.

\vskip 0.3 cm 
\noindent
{\bf $\diamond$------------------------------------------------------------------------------------------------------------------$\diamond$}

\vskip 0.3 cm
{\bf 3.1)} The symmetry (\ref{sym2}) applied to case 3),  gives  solutions with \fbox{$\sqrt{2\alpha}+\sqrt{2\gamma}=N$} or \fbox{$ \sqrt{2\alpha}-\sqrt{2\gamma}=N$}  (namely $\theta_\infty-1+\theta_1=N$, or $\theta_\infty-1-\theta_1=N$), $N\in{\bf Z}$: 
\be
\label{LOG12}
\sq{
  y(x)= \sum_{n=0}^{|N|-1}b_n x^n+\Bigl(a+b_N\ln x\Bigr)x^{|N|}+\sum_{n= |N|+1}^\infty P_n(\ln x;a)x^n,~N\neq 0
}
\ee
\be
\label{LOG45}
\sq{
y(x)= {1\over (a\pm \sqrt{2\alpha}\ln x)+\sum_{n=1}^\infty P_n(\ln x;a)x^n}=\pm{1\over \sqrt{2\alpha}\ln x}
\left(
1\mp {1\over \sqrt{2\alpha}\ln x }+
O\left(
{1\over \ln^2 x}
\right)
\right),\left.\matrix{N=0 \cr\alpha=\gamma\neq 0
}\right.
}
\ee
 The $b_n$'s are rational functions of $\sqrt{\alpha},\sqrt{\beta},\sqrt{\gamma}$ and $\sqrt{1-2\delta}$; $a$ is free complex parameter; $P_n(\ln x;a)$ are polynomials of $\ln x$ of degree $n-|N|+1$, with coefficients which are rational functions of $\sqrt{\alpha},\sqrt{\beta},\sqrt{\gamma},\sqrt{1-2\delta}$ and $a$. The  solutions become convergent Taylor expansions, either if $\{\sqrt{-2\beta}+\sqrt{1-2\delta}, \sqrt{-2\beta}-\sqrt{1-2\delta} \}\cap{\cal N}_N\neq\emptyset$, 
namely $
   \{ \theta_0+\theta_x,\theta_0-\theta_x\}\cap{\cal N}_N\neq\emptyset
$,  
 or if $\sqrt{2\alpha}\in \{-1,...,N\}$, if $N<0$,  $\sqrt{2\alpha}\in\{1,...,N\}$, if $N>0$.  
The expansions are 
\be
\label{TTLO3}
 \sq{
 y(x)= \sum_{n=0}^{|N|-1}b_n x^n+ax^{|N|}+\sum_{n= |N|+1}^\infty b_n(a)x^n,~~~N\neq 0
}
\ee
\be
\label{TTLO4}
\sq{
y(x)=a+(1-a)(\delta-\beta)x+\sum_{n= 2}^\infty b_n(a)x^n,~~N= 0,~~\alpha=\gamma=0
}
\ee
We observe that the case (\ref{TTLO3}) with $N=1$ gives two subcases:
\be
\label{T2coe}
y(x)={1\over \sqrt{2\alpha}} +ax+O(x^2),~~~~~\alpha\neq 0,~~(\sqrt{2\alpha}\pm\sqrt{2\gamma})^2=1,~~1-2\delta+2\beta=0
\ee
 and  
$$
y(x)=1+ax+{1\over 4} a(1+4a+2\beta-2\delta)x^2+O(x^3),~~~\alpha={1\over 2},~~\gamma=0
$$
The latter  is equivalent to (\ref{TTLO4}) via a bi-rational transformation. 

Solutions (\ref{TTLO4}), (\ref{T2coe}) and (\ref{T1coe}) are basic Taylor solutions, which generate all the other Taylor expansions by the action of the bi-rational transformations. They are obtained in \cite{D2}  with non Fuchsian  reductions  of (\ref{SYSTEM}). They are also derived in \cite{kaneko}. 

\vskip 0.2 cm
{\bf {\bf (*)}} In terms of the classification of \cite{Bruno7}, the class of the solutions (\ref{LOG12}) is the subclass of ${\cal B}_1\cup {\cal B}_2$ with logarithmic coefficients. The class of the solutions (\ref{LOG45}) coincides with the class ${\cal B}_4\cup {\cal B}_5$. 
The sub case (\ref{TTLO4})  gives the class ${\cal B}_{10}$ of \cite{Bruno7}.

\vskip 0.3 cm 
\noindent
{\bf $\diamond$------------------------------------------------------------------------------------------------------------------$\diamond$}

\vskip 0.3 cm
{\bf 4.1)}  Symmetry (\ref{sym2}) applied to 4) gives the following solutions defined for \fbox{$\alpha\neq \gamma$} (namely $\theta_\infty \pm \theta_1\neq 1$): 
\be
\label{LOG3}
\sq{
y(x)= {1\over\left[ {\gamma-\alpha\over 2} (a+\ln x)^2 +{\alpha\over \alpha -\gamma}\right] +\sum_{n=1}^\infty P_{n+1}(\ln x) x^n}
 =
{2\over (\gamma-\alpha)\ln^2 x}\left[
1-{2a\over \ln x} +O\left({1\over \ln^2 x}\right)
\right],~~~\alpha\neq \gamma
}
\ee

{\bf {\bf (*)}} In terms of the classification of \cite{Bruno7}, this is the class ${\cal B}_3$. 


\subsection{Comparison with the Results of Power Geometry} 
\label{COMOPAIR}
 Summarising the paragraphs {\bf {\bf (*)}},  we see that the results of the monodromy preserving deformations method and of power geometry coincide, as follows:
$$
\left.
\matrix{
\hbox{\bf Table 1 } & & \hbox{\bf 39 families of \cite{Bruno7}}
\cr
\cr
{\bf 1)} \hbox{ behaviour } (\ref{fullEXP}), \hbox{ with $\Re \sigma\neq 0$} & &
\cr
\hbox{ [ This includes cases {\bf 1.1a)},}& =& {\cal A}_0
\cr
 \hbox{ {\bf 2)} with $-1<\Re\sigma<1$, and {\bf 2.1c)}]}  && 
\cr
\cr
{\bf 1.1b)}\hbox{ behaviour  (\ref{lantern1}), and {\bf 1)} behaviour (\ref{lantern})}&=&  {\cal B}_0^{+}\cup {\cal B}_0^{-}
\cup {\cal B}_7^{+}\cup {\cal B}_7^{-} \hbox{ and image}^1
\cr
\cr
   {\bf 2.1a)}\hbox{ behaviour (\ref{UUU}), and {\bf 2)} behaviour (\ref{atopy}) if $\Re\sigma>0$ }&=& {\cal B}_1\hbox{(non log)}\cup{\cal B}_2\hbox{(non log)}\cup {\cal B}_6\hbox{ and image}
\cr
\cr
   {\bf 2.1b)}\hbox{ behaviour (\ref{TAU}), and {\bf 2)} behaviour (\ref{atopy}) if $\Re\sigma=0$ }&= &{\cal B}_1^{+}\cup{\cal B}_1^{-} \cup{\cal B}_2^{+}\cup{\cal B}_2^{-} \cup{\cal B}_6^{+}\cup{\cal B}_6^{-} \hbox{ and image}
\cr
\cr
   {\bf 2.1d)}\hbox{ behaviour (\ref{div}), and {\bf 2)} behaviour (\ref{atopy}) if $\beta=0$ }&=&{\cal C}_0^\infty \cup {\cal B}_8\cup{\cal B}_9  \hbox{ and image}
\cr
\cr
{\bf 3.1)}, \hbox{ and {\bf 3)}}&=& {\cal B}_1\hbox{(log)}\cup{\cal B}_2\hbox{(log)}\cup {\cal B}_4\cup{\cal B}_5\cup {\cal B}_{10}  \hbox{ and image}
\cr
\cr
{\bf 4.1)},\hbox{ and {\bf 4)}} &=& {\cal B}_3  \hbox{ and image}
}
\right.
$$

\vskip 0.2 cm 
\noindent
${}^1$``image'' means the image of the families  through (\ref{sym2}).

\vskip 0.2 cm
\noindent 
In \cite{Bruno7} a classification of actual, asymptotic or formal expansions which satisfy an ODE is given. It terms of it, we have that: 

\noindent
- Solutions (\ref{fullEXP}) and (\ref{atopy}) for $\Re\sigma\neq 0$, (\ref{UUU}) for $\Re \rho \neq 0$, (\ref{div}) for $\Re \omega\neq 0$,  and the Taylor series   (\ref{TTLO1}), (\ref{TTLO2}), (\ref{TTLO3}), (\ref{TTLO4}), (\ref{davidekan}) and  (\ref{T1coe})  have {\it power expansions}. 

\noindent
- Solutions (\ref{log1}),   (\ref{log1zero}), (\ref{logsquare}) and (\ref{LOG12})  have  {\it power logarithmic expansions}.

\noindent
- Solutions (\ref{LOG45}) and (\ref{LOG3}), when expanded in powers of $(\ln x)^{-1}$, have {\it complicated expansions}. 

\noindent
- Solutions (\ref{lantern}) (namely, (\ref{fullEXP}) for $\Re\sigma= 0$) and (\ref{atopy}) for $\Re\sigma= 0$, (\ref{UUU}) for $\Re\rho=0$, and  (\ref{div}) for $\Re \omega= 0$ have {\it semi exotic expansions}.
 
\noindent
- Solutions (\ref{lantern1}) and (\ref{TAU}),  when formally expanded (like in  (\ref{formal})), 
have {\it exotic expansions}.




\section{Appendix A: Bi-rational Transformations}
\label{simmetries}

{Okamoto's bi-rational transformations} \cite{Okamoto} are
 symmetries of (PVI), namely invertible transformations: 
\be
\label{OKAvino}
y^{\prime}(x)={P\left(x,y(x),{dy(x)\over dx}\right)\over Q\left(x,y(x),{dy(x)\over dx}\right)}, ~~~x^{\prime}={p(x)\over q(x)},~~~~~
(\theta_0,\theta_x,\theta_1,\theta_\infty)\mapsto 
(\theta_0^{\prime},\theta_x^{\prime},\theta_1^{\prime},\theta_\infty^{\prime})
\ee
such that $y(x)$ satisfies (PVI) with coefficients $\theta_0,\theta_x,\theta_1,\theta_\infty$ and variable $x$,  if and only if $y^{\prime}(x^{\prime})$ satisfies (PVI) with coefficients $\theta_0^{\prime},\theta_x^{\prime},\theta_1^{\prime},\theta_\infty^{\prime}$ and variable $x^\prime$. The functions 
 $P,Q$ are polynomials; $p,q$ are linear; the transformation of the $\theta_\mu$'s
  is an element of a linear representation of one of the following groups.

\vskip 0.2 cm
* Permutation group:
\be
\theta_0^{\prime}=\theta_1,~~\theta_x^{\prime}=\theta_x,~~\theta_1^{\prime}=\theta_0,~~\theta_\infty^{\prime}=\theta_\infty;~~~y^{\prime}(x^\prime)=1-y(x),~~~x^\prime=1-x. 
\label{onara}
\ee
 
\be
\label{onara1}
\theta_x^\prime=\theta_1,~~\theta_1^\prime=\theta_x;~~~~~~~
\theta_0^\prime=\theta_0,~~\theta_\infty^\prime=\theta_\infty;~~~~~~~~~~
y^\prime(x^\prime)={1\over x}y(x),~~~x^\prime={1\over x}.
\ee

\be
\label{sym2}
\theta_0^\prime=\theta_\infty-1,~~\theta_x^\prime=\theta_1,~~\theta_1^\prime=\theta_x,~~\theta_\infty^\prime=\theta_0+1;~~~~~~~y^\prime(x^\prime)={x\over y(x)},~~~x=x^\prime.
\ee 
In terms of the coefficients of PVI, the transformation (\ref{sym2}) is $
\alpha^\prime= -\beta$, $\beta^\prime=-\alpha$, $\gamma^\prime={1\over 2}-\delta$, $\delta^\prime={1\over 2}-\gamma
$.

\vskip 0.3 cm
* Weyl Group  of the root system $D_4$: 
\vskip 0.2 cm
\noindent
$w_1$:
$$\theta_1^\prime=-\theta_1;~~~~~\theta_0^\prime=\theta_0,~~ \theta_x^\prime= \theta_x,~~ \theta_\infty^\prime =\theta_\infty.
$$
$w_2$:
$$
\theta_0^\prime={\theta_0+\theta_1+\theta_x+\theta_\infty\over 2} -1,~~~
\theta_1^\prime={\theta_0+\theta_1-\theta_x-\theta_{\infty}\over 2}+ 1,
$$
$$
\theta_x^\prime={\theta_0-\theta_1+\theta_x-\theta_{\infty}\over 2} +1,~~~
\theta_{\infty}^\prime={\theta_0-\theta_1-\theta_x+\theta_\infty\over 2}+1
$$
$
w_3$:
$$
\theta_\infty^\prime=2-\theta_\infty;~~~~~\theta_0^\prime=\theta_0,~~ \theta_x^\prime= \theta_x,~~ \theta_1^\prime =\theta_1.
$$
$
w_4$:
$$
\theta_\infty^\prime=2-\theta_\infty;~~\theta_x^\prime=2-\theta_x;~~~~~ \theta_0^\prime= \theta_0,~~ \theta_1^\prime =\theta_1.
$$
\vskip 0.2 cm
\noindent
The variable $x^\prime=x$, but $y^{\prime}(x)$ is quite complicated
and will not be given here.

\vskip 0.3 cm

* Shift $l_j: v=(v_1,v_2,v_3,v_4)\mapsto 
v+e_j$, $j=1,2,3,4$, where  $e_1=(1,0,0,0)$, ..., $e_4=(0,0,0,1)$:
$$
l_1:~~~~~~~~~~\theta_0^\prime=\theta_0+1,~~\theta_1^\prime=\theta_1+1;~~~~~
 \theta_x^\prime=\theta_x,~~\theta_\infty^\prime=\theta_\infty.
$$
$$
l_2:~~~~~~~~~~\theta_0^\prime=\theta_0+1,~~\theta_1^\prime=\theta_1-1;~~~~~
 \theta_x^\prime=\theta_x,~~\theta_\infty^\prime=\theta_\infty.
$$
$$
l_3:~~~~~~~~~~\theta_x^\prime=\theta_x+1,~~\theta_\infty^\prime=\theta_\infty+1;~~~~~
 \theta_0^\prime=\theta_0,~~\theta_1^\prime=\theta_1.
$$
$$
l_4:~~~~~~~~~~\theta_x^\prime=\theta_x+1,~~\theta_\infty^\prime=\theta_\infty-1;~~~~~
 \theta_0^\prime=\theta_0,~~\theta_1^\prime=\theta_1.
$$
The variable $x^\prime=x$, but $y^{\prime}(x)$ is quite complicated and will not be given here.

\vskip 0.2 cm 
 Transformations $w_1,w_2,w_3,w_4$ and $l_3$ are sufficient to generate the other $l_j$'s, and they give a representation of the affine Weyl group  of $D_4$. 
The transformations of $y(x)$ corresponding to $w_1,w_2,w_3$ and $w_4$ are obtained from formula (2.10) at page 355 of \cite{Okamoto}, while the transformation corresponding to $l_3$ is obtained from formulae (2.5), (2.6) at page 354 of \cite{Okamoto}, provided that $h$  is substituted with $h^{+}$ of (1.13), page 352.

\vskip 0.3 cm 
 The effect of the bi-rational transformation on the other monodromy data can be found in \cite{Jimbo}, \cite{DM}, \cite{D1}, and more extensively and generally in \cite{DM11} (see also \cite{Marta}). In particular, we have: 
\vskip 0.2 cm  

-- For (\ref{onara}): 
\be
\label{ricordo11}
\left\{
\matrix{
p_{01}^\prime
&=&-p_{01}-p_{0x}p_{x1}+p_\infty p_x+p_1p_0 
\cr
\cr
p_{0x}^\prime&=&p_{x1}~~~~~~~~~~~~~~~~~~~~~~~~~~~~~~~                                     
\cr
\cr
p_{x1}^\prime&=&p_{0x}~~~~~~~~~~~~~~~~~~~~~~~~~~~~~~~~
}
\right.
~\hbox{ or }
\left\{
\matrix{
p_{01}&=&-p_{01}^\prime-p_{0x}^\prime p_{x1}^\prime+p_\infty^\prime
p_x^\prime+p_1^\prime p_0^\prime 
\cr
\cr
p_{0x}&=&p_{1x}^\prime~~~~~~~~~~~~~~~~~~~~~~~~~~~~~~~
\cr
\cr
p_{x1}&=&p_{0x}^\prime~~~~~~~~~~~~~~~~~~~~~~~~~~~~~~~
}
\right.
\ee

-- For (\ref{onara1}): 
\be
\label{ricordo22}
\left\{
\matrix{p_{0x}^\prime &=&
 -p_{01}-p_{0x}p_{x1} +p_\infty p_x +p_0p_1 
\cr
\cr
p_{01}^\prime&=&p_{0x}~~~~~~~~~~~~~~~~~~~~~~~~~~~~~~~
\cr\cr
p_{1x}^\prime&=&p_{1x}~~~~~~~~~~~~~~~~~~~~~~~~~~~~~~~
}
\right.
~\hbox{ or }
\left\{\matrix{p_{01} &= &-p_{0x}^\prime-p_{01}^\prime p_{x1}^\prime+p_\infty^\prime
p_1^\prime+p_0^\prime p_x^\prime
\cr\cr
p_{0x}&=&p_{01}^\prime~~~~~~~~~~~~~~~~~~~~~~~~~~~~~~~
\cr\cr
p_{1x}&=&p_{1x}^\prime~~~~~~~~~~~~~~~~~~~~~~~~~~~~~~~
}
\right.
\ee

-- For  (\ref{sym2}):
\be
\label{ricordo3}
p_{0x}^\prime=-p_{0x},~~~p_{01}^\prime=-p_{01},~~~p_{x1}^\prime=p_{x1}
\ee
 The parametrisation for the solutions in cases {\bf N.1)},  N=1,2,3,4, of section \ref{esplaintable} are obtained from that of the equivalent  cases  {\bf N)}  by means of the above (\ref{ricordo3}) [note: $\sigma^\prime=1-\sigma$ when $\sigma$ appears].

\vskip 0.3 cm 

{\bf Poof of Proposition \ref{spero}}:  Tables 2 and 3 are obtained form Table 1 via (\ref{onara}) and (\ref{onara1}) rsp. We  prove that the expansions of Table 1  are obtained form the set of basic expansions   (\ref{fullEXP}), (\ref{atopy}),  (\ref{log1zero}) and (\ref{logsquare}).

1st Step. The complex power behaviours, the Taylor expansions and the power logarithmic expansions given in the table are the only allowed forms of these kinds. This is proved by direct substitution into PVI. 

 2nd Step. The complex power behaviours (\ref{UUU}) and (\ref{div}), and the inverse oscillatory behaviours are obtained from  (\ref{fullEXP}) and  (\ref{atopy}) via (\ref{sym2}). The inverse power logarithmic behaviours and the power logarithmic behaviours (\ref{LOG12}) are obtained from (\ref{log1}), (\ref{log1zero}) and (\ref{logsquare}) via (\ref{sym2}). Details are explained in section \ref{esplaintable}.  

 3rd Step.  All solutions  (\ref{log1})  are obtained from (\ref{log1zero}), by a  bi-rational transformation belonging to the affine Weil group of $D_4$, because these transformations  change $\theta_0\pm \theta_x$ into $\theta_0\pm \theta_x+N$, $N\in {\bf Z}$.  It is a formidable task to obtain $y(x)$ by repeated applications of bi-rational transformations, but one can avoid this by observing that a transformation (\ref{OKAvino}) applied to a power logarithmic behaviour necessarily produces a power logarithmic behaviour or an inverse power logarithmic behaviour. Then, the conclusion comes  from the fact that the expansions in the table exhaust all possible forms of such behaviours, being relations like $\theta_\mu\pm \theta_\nu\in{\bf Z}$, or $\not\in{\bf Z}$, preserved by bi-rational transformations of the affine $D_4$.

4th Step.  All Taylor expansions are obtained from  (\ref{log1zero}). For example,  when we apply the bi-rational transformation $l_3$ (or $l_3^{-1}$)  to (\ref{log1zero}), we obtain the case  $N=1$ of solution   (\ref{log1}), namely:
\be
\label{uiui}
 y(x)=\pm\sqrt{-2\beta}~x+\left(
a-{1\over 2}\sqrt{-2\beta}\left(\sqrt{-2\beta}\mp 1\right)\left(\alpha-\gamma\right)\ln x
\right) x^2+\sum_{n=2}^\infty x^n\sum_{m=0}^{n-1}b_{nm}(a)\ln^m x,
\ee
$$
\sqrt{-2\beta}+\sqrt{1-2\delta}=\pm 1 \hbox{ or }\sqrt{-2\beta}-\sqrt{1-2\delta}=\pm 1 $$ 
Then, let $\beta=0$ in (\ref{log1zero}) and $\alpha=\gamma$ in (\ref{uiui}). We obtain respectively the two basic Taylor solutions: 
$$
\left.
\matrix{
y(x)=a+{a(a-1)\over 2}(2\gamma-2\alpha-1)x^2+\sum_{n=3}^\infty b_n(a) x^n, & \hbox{this is  (\ref{TTLO2}) in the table}
\cr
\cr
y(x)=\pm\sqrt{-2\beta}~x+ax^2+\sum_{n=3}^\infty b_n(a) x^n, & \hbox{ this is (\ref{batsta}), section \ref{esplaintable}}
}
\right.
$$
These are transformed  by (\ref{sym2}) into (\ref{TTLO4}) and (\ref{T2coe}), which may also be used as basic Taylor solutions. Then,  bi-rational transformations belonging to the affine Weil group of $D_4$ applied the basic Taylor give all the other Taylor solutions at $x=0$. $\Box$


\section{Appendix B: Completeness of Table 1 -- Partial Proof}

 We analyse the conjecture  of Section  \ref{congettura}.   
It is enough to prove that Table 1 is complete, being Table 2 and 3 constructed from 1. The conjecture  will be proved when we are able to prove that  the  restriction to  Table 1 of the  map  (\ref{monmap}), namely
  $$ 
  f:\{y(x) \hbox{ branches in  Table 1}\}\rightarrow {\cal M}
$$
is surjective. Let us  write the disjoint union ${\cal M}={\cal M}_r\cup{\cal M}_i$, defining  ${\cal M}_r$ and ${\cal M}_i$ to be the subsets where the group 
$<M_0,M_x,M_1>$ is reducible or irreducible respectively (it is understood that the elements are defined up to conjugation).  
 Because $f:f^{-1}({\cal M}_i)\to {\cal M}_i$ is one to one, we should prove that that:

\vskip 0.2 cm 
{\bf I)} $f^{-1}({\cal M}_r)\subset \{y(x)$ branches of in  table$\}$, 

{\bf II)} $f(\{y(x)$ in the table such that the corresponding monodromy group is irreducible$\})={\cal M}_i$.

\vskip 0.3 cm 
\noindent
{\bf I)}  When  $
<M_0,M_x,M_1>$  is reducible, all the corresponding solutions of (PVI) 
 are  equivalent by bi-rational canonical transformations 
 to the  one-parameter family of solutions \cite{hitchin} 
\be 
y(x)=
{\theta_1+\theta_\infty-1+x(1+\theta_x)\over \theta_\infty-1}-
{1\over \theta_\infty-1}~{x~(1-x)\over u(x;a)}~{du(x;a)\over dx},~~~\theta_\infty+\theta_1+\theta_0+\theta_x=0
\label{sabishii}
\ee
where $u(x;a)=u_1(x)+a u_2(x)$; $a \in {\bf C}$,  $u_1(x)$ and
$u_2(x)$ 
are linear independent solutions of the hyper-geometric equation: 
$$
x(1-x) {d^2 u\over dx^2} +
\left\{
[2-(\theta_\infty+\theta_1)]-(4-\theta_\infty+\theta_x)x
\right\}{du\over dx} 
-(2-\theta_\infty)(1+\theta_x)u=0
$$
Their behaviour is  logarithmic, Taylor or (\ref{UUU}), but the parametrisation in terms of $p_{0x}$, $p_{x1}$, $p_{01}$ does not apply\footnote{The rational solutions fall into the above case, see \cite{mazzoccoratio}.}. Thus, I) is true. 
%
%

\vskip 0.3 cm 
\noindent
{\bf II)} We know  the monodromy group  $<M_0,M_x,M_1>$ associated to the solutions of the table, when this group is irreducible. This is precisely the case when the parametric formulae in terms of $p_{ij}$ and $\theta_\mu$ are computed.  Based on  \cite{Jimbo},\cite{Boalch},\cite{D3},\cite{D2},\cite{D1},\cite{guz2010}  and keeping into account the effect of bi-rational transformations on the monodromy data, we conclude the following:

\vskip 0.3 cm 
\noindent
{\bf II,~a)} The generic case  of complex power (\ref{fullEXP}) and inverse oscillatory behaviours (\ref{lantern1}). 
\vskip 0.3 cm 
\noindent
$\diamond$ 
 $
f(\{y(x)= (\ref{fullEXP}),(\ref{lantern1})\}) 
$

\noindent
$=
\left\{ <M_0,M_x> \hbox{ irreducible},~\hbox{tr}(M_xM_0)\not\in\bigl\{-2,~2,~2\cos\bigl(\pi(\theta_0+\theta_x)\bigr),
~2\cos\bigl(\pi(\theta_0-\theta_x)\bigr)\bigr\}\right\}
$

\vskip 0.2 cm 
\noindent
Note that $\hbox{tr}(M_xM_0)\not \in\bigl\{2\cos(\pi(\theta_0+\theta_x)),2\cos(\pi(\theta_0-\theta_x))\bigr\}$ implies that $<M_0,M_x>$ is irreducible, but does not exclude reducibility of other subgroups, like $<M_xM_0,M_1>$ (see below).

\vskip 0.3 cm 
\noindent
{\bf II,~b)} Cases when the subgroup $<M_0,M_x>$ is reducible. In this case $\hbox{tr}(M_xM_0)\in\bigl\{2\cos(\pi(\theta_0+\theta_x)),2\cos(\pi(\theta_0-\theta_x))\bigr\}$.    

\vskip 0.3 cm 
\noindent
$\diamond$ 
$
f(\{y=(\ref{atopy}),(\ref{davidekan})\})
$

\noindent
$=\bigcup_{k=1,2}\{ \hbox{tr}(M_xM_0)\neq \pm2,~<M_0,M_x>\hbox{ reducible}, 
~\theta_0+(-)^k\theta_x\not\in{\bf Z}\}
$

\vskip 0.2 cm 
\noindent
 The two signs in $\theta_0\pm \theta_x$ correspond to $\pm$ in the first coefficient ${\theta_0\over \theta_0\pm \theta_x}$ in (\ref{atopy}). Case (\ref{davidekan}) is a subcase of (\ref{atopy}) for $a=0$.   

\vskip 0.3 cm 
\noindent
$\diamond$ 
$
f(\{y=(\ref{log1}),(\ref{log1zero}),(\ref{TTLO1}), (\ref{TTLO2})\})
$

\noindent
$
=\bigcup_{k=1,2}\bigcup_{N\in{\bf Z}}\Bigl\{ \hbox{tr}(M_xM_0)=2\cos \pi N=\pm 2, ~<M_0,M_x>\hbox{ reducible},~ \theta_0+(-)^k\theta_x=N\}
$

\vskip 0.2 cm 
\noindent
The branches (\ref{log1}) occur when $N\neq 0$, and become (\ref{TTLO1}) if either $\theta_0\in\{0,\pm 1,...,\pm N\}$ or one of $(\theta_\infty-1)\pm \theta_1\in{\cal N}_N$. (\ref{log1zero}) occur when $N=0$. (\ref{TTLO2}) occur when $N=0=\theta_0=\theta_x$.

\vskip 0.3 cm 
\noindent
{\bf II,~c)} Cases when $<M_0,M_x>$ is irreducible and $<M_xM_0,M_1>$ is  reducible, so that   $\hbox{tr}(M_xM_0)\in \bigl\{2\cos(\pi(\theta_\infty+\theta_1)),~2\cos(\pi(\theta_\infty-\theta_1))\bigr\}$.

\vskip 0.3 cm 
\noindent
$\diamond$ 
$
f(\{y=(\ref{UUU}),(\ref{T1coe}),(\ref{div}),(\ref{TAU})\})
$

\noindent
$
=\{\hbox{tr}(M_xM_0)\neq \pm 2,~<M_0,M_x>\hbox{ irreducible}, ~ <M_xM_0,M_1> \hbox{ reducible },~ (\theta_\infty-1)\pm \theta_1\not\in{\bf Z}\}
$

\vskip 0.2 cm
\noindent
 The two signs in $(\theta_\infty-1)\pm \theta_1$ correspond to $\pm$ in the first coefficient of (\ref{UUU}). Note that (\ref{T1coe}) is a subcase of (\ref{UUU}), and that (\ref{div}) occurs for $\theta_\infty=1$. 

\vskip 0.3 cm 
\noindent
$\diamond$ 
$
f(\{y=(\ref{LOG12}),(\ref{LOG45}),(\ref{TTLO3}),(\ref{TTLO4})\})
$

\noindent
$=
\bigcup_{k=1,2}\bigcup_{N\in{\bf Z}}\Bigl\{\hbox{tr}(M_xM_0)=-2\cos\pi N=\pm 2,~<M_0,M_x>\hbox{ irreducible},~ <M_xM_0,M_1> \hbox{ reducible},~(\theta_\infty-1)+(-)^k\theta_1=N\}
$

\vskip 0.2 cm 
\noindent
The branches (\ref{LOG12}) occur for $N\neq 0$, and become (\ref{TTLO3}) if either $\theta_\infty-1\in\{\pm1,...,\pm N\}$ or one of $\theta_0\pm \theta_x \in{\cal N}_N$.  (\ref{LOG45}) occur for $N=0$ and $\theta_\infty-1\neq 0, \theta_1\neq 0$.   (\ref{TTLO4}) occur if $N=0=\theta_\infty-1=\theta_1$. 

\vskip 0.3 cm 
\noindent
$\diamond$ 
 The subset of $f(\{y(x)= (\ref{fullEXP})\})$ s.t. $<M_xM_0,M_1>$ is reducible. It contains the image of (\ref{atopy}) through the bi-rational transformation (\ref{sym2})  [see case 2.1c) and 2.1d)  in Section \ref{esplaintable}].

\vskip 0.3 cm 
\noindent
{\bf II, d)} Cases when both $<M_0,M_x>$ and $<M_xM_0,M_1>$ are irreducible. 

\vskip 0.3 cm 
\noindent
$\diamond$ 
The  subset  of  $
f(\{y(x)= (\ref{fullEXP}),(\ref{lantern1})\}) 
$ s.t. $<M_xM_0,M_1>$ is irreducible. 

\vskip 0.3 cm 
\noindent
$\diamond$ 
$
f(\{y=(\ref{logsquare})\})
$

\noindent
$=
\{\hbox{tr}(M_xM_0)=2,~<M_0,M_x>\hbox{ and } <M_xM_0,M_1> \hbox{ irreducible},~\theta_0\pm\theta_x\neq 0\}
$

\vskip 0.2 cm 
\noindent 
 If one of $\theta_0\pm\theta_x\to 0$, then $y=(\ref{logsquare})\to (\ref{TTLO2})$ and $<M_0,M_x>$ becomes reducible. 

\vskip 0.3 cm 
\noindent
$\diamond$ 
$
f(\{y=(\ref{LOG3})\})
$

\noindent
$
=\{\hbox{tr}(M_xM_0)=-2,~<M_0,M_x>\hbox{ and } <M_xM_0,M_1> \hbox{ irreducible},~(\theta_\infty-1)\pm \theta_1\neq 0\}
$

\vskip 0.2 cm 
\noindent 
 If one of $(\theta_\infty-1)\pm \theta_1\to 0$, then $y=(\ref{LOG3})\to  (\ref{TTLO4})$ and $<M_xM_0,M_1>$ becomes reducible. This case is considered above.

\vskip 0.2 cm 
\noindent
The images through $f$ obtained in cases II, a), II, b), II, c) and II, d) is contained in  ${\cal M}_i$ and {\it it is  almost the whole}  ${\cal M}_i$, though we are not able to prove that {\it it is} ${\cal M}_i$. 
$\Box$


\end{document}